\documentclass[10pt,aps,onecolumn,superscriptaddress]{revtex4-2}

\usepackage[T1]{fontenc}
\usepackage[utf8]{inputenc}
\usepackage[english]{babel}

\usepackage{amsfonts,amssymb,amsmath}
\usepackage{graphicx,float}

\usepackage{hyperref}
\hypersetup{
    colorlinks=true,
    linkcolor=blue,
    filecolor=magenta,      
    urlcolor=cyan,
    citecolor=blue,
}

\graphicspath{{figures/}}

\begin{document}
	
	\title{Painting the Phase Space of Dissipative Systems \\ with Lagrangian Descriptors}
	
	\author{V\'ictor J. Garc\'ia-Garrido}
	\email{vjose.garcia@uah.es}
	\affiliation{Departamento de F\'isica y Matem\'aticas, \\ Universidad de Alcal\'a, Madrid, 28871, Spain. \\[.2cm]}
	\author{Julia Garc\'{i}a-Luengo}
	\email{julia.gluengo@upm.es}
	\affiliation{Departamento de Matem\'{a}tica Aplicada a las TIC, \\ Universidad Polit\'{e}cnica de Madrid, Madrid, 28031, Spain.}

	\begin{abstract}
	In this paper we apply the method of Lagrangian descriptors to explore the geometrical structures in phase space that govern the dynamics of dissipative systems. We demonstrate through many classical examples taken from the nonlinear dynamics literature that this tool can provide valuable information and insights to develop a more general and detailed understanding of the global behavior and underlying geometry of these systems. In order to achieve this goal, we analyze systems that display dynamical features such as hyperbolic points with different expansion and contraction rates, limit cycles, slow manifolds and strange attractors. Furthermore, we study how this technique can be used to detect transition ellipsoids that arise in Hamiltonian systems subject to dissipative forces, and which play a crucial role in characterizing trajectories that evolve across an index-1 saddle point of the underlying potential energy surface.
	\end{abstract}
	
	\maketitle


\section{Introduction}
\label{sec:intro}

Dissipative dynamical systems are ubiquitous in Nature. For example, in mechanical systems they occur when energy is pumped into the system through an external driving force acting on it, or energy is drained by means of mechanisms such as friction. The process by which dissipation manifests itself in phase space is through the contraction and/or expansion of volumes of ensembles of initial conditions as they evolve in time. These systems display a rich variety of dynamical behavior governed by objects such as limit cycles, attractors, slow manifolds, etc., and their understanding is of paramount importance to unravel the underlying dynamics at play \cite{guck1983,nayfeh1995,strogatz,kuehn2015,meiss2017}. 

The goal we pursue in this work is to illustrate how the method of Lagrangian descriptors (LDs) can be used to explore the phase space of dynamical systems subject dissipation. We believe that this technique would be of great help to gain insight into the dynamical mechanisms that take place in such systems, and provide some guidance in order to carry out a detailed analysis of any system of interest. The advantage of applying LDs is that it is a scalar-based trajectory diagnostic tool that is very simple to implement computationally. This technique was first developed a decade ago for the study of Lagrangian transport and mixing in geophysical flows \cite{madrid2009,mendoza2010}, and in this context it has been applied for instance to the real-time assessment of oil spills \cite{gg2016,guillermo2021} and also to the path planning of transoceanic missions of autonomous underwater vehicle missions \cite{ramos2018}. The identification of hyperbolic trajectories and their stable and unstable manifolds in the time-dependent vector fields defined by ocean currents plays a crucial role in the analysis of these problems, and it was the main focus of these studies. However, other geometrical objects (limit cycles, attractors and slow manifolds) that are relevant for dissipative systems have not been analyzed in the literature with this tool. The contribution of this work is aimed towards filling this gap.

The capability of LDs to reveal phase space structure is addressed here through many classical examples from nonlinear dynamics. The contents of this paper are outlined as follows. In Section \ref{sec:sec1} we provide a brief introduction on the method of Lagrangian descriptors, and describe the  numerical setup we have used for this tool throughout this work when exploring the phase sapce of dissipative systems. Section \ref{sec:sec2} is devoted to presenting the results obtained from our analysis of several well-known model dynamical systems by means of LDs. Firstly, in Sec. \ref{subsec:sec1} we look at simple saddle systems that display different expansion and contraction rates, and thus their dynamics are characterized by two distinct time scales. Then, in Sec. \ref{subsec:sec2} we move on to illustrate how LDs can be used to detect bifurcations, and in particular, Andronov-Hopf bifurcations that give rise to the birth of limit cycles. Moreover, we demonstrate that LDs are capable of revealing the limit cycle structures in the classical van der Pol oscillator. Section \ref{subsec:sec3} is focused on the analysis of systems with slow manifolds, and Sec. \ref{subsec:sec4} describes how this tools succeeds in highlighting the intricate and fractal structure of attractors and strange attractors in the Duffing oscillator system. We extend our study of dissipative dynamics further in Sec. \ref{subsec:sec5}, where we illustrate the potential that LDs brings for the analysis of transition ellipsoids in Hamiltonian systems subject to dissipative forces. The development of techniques such as LDs for identifying these geometrical objects in phase space is of paramount importance, since they greatly facilitate the dynamical study of trajectories evolving across an index-1 saddle of the underlying potential energy surface that describes the model system. Finally, in Sec. \ref{sec:conc} we present the conclusions of this work and outlook for future research.

\section{Lagrangian Descriptors}
\label{sec:sec1}

This section briefly introduces the method of Lagrangian descriptors (LDs) that we have used in this work to analyze the phase space structures that govern the dynamics of dissipative systems. This technique is a scalar-based diagnostic tool that was originally introduced in the context of fluid mechanics to study transport and mixing in geophysical flows \cite{madrid2009,mendoza2010,mancho2013lagrangian}, but in the past years it has been found to provide useful information for the analysis of the high dimensional phase space of chemical systems, see e.g. \cite{Agaoglou2019,ldbook2020} and references therein. Although the method was initially developed to analyze continuous-time dynamical systems \cite{mancho2013lagrangian,lopesino2017} with general time dependence, it has also been applied to maps \cite{carlos,GG2019b}, stochastic dynamical systems \cite{balibrea2016} and holomorphic dynamics  \cite{GG2020a}.

The idea behind the method of Lagrangian descriptors is very simple. Take any initial condition and accumulate along its trajectory the values attained by a positive scalar function that depends on the phase space variables. This calculation is carried out both forward and backward in time as the system evolves. Once this computation is done for a grid of initial conditions, the scalar output obtained from the method will highlight the location of the invariant stable and unstable manifolds intersecting this slice. These manifolds will appear as 'singular features' of the scalar field, that is, they are detected at points where the values of LDs display an abrupt change, which indicates distinct dynamical behavior. Forward integration of trajectories detects stable manifolds while backward evolution does the same for unstable manifolds. It is also important to remark here that this technique also reveals the structure of KAM tori, but in this case the tori regions correspond to smooth values of the LD function. In fact, tori can be visualized by computing long-term time averages of LDs as discussed in \cite{lopesino2017,GG2019b}.

Consider a general time-dependent dynamical system given by the evolution law:
\begin{equation}
\dfrac{d\mathbf{x}}{dt} = \mathbf{f}(\mathbf{x},t) \;,\quad \mathbf{x} \in \mathbb{R}^{n} \;,\; t \in \mathbb{R} \;,
\end{equation}
where the vector field defining flow satisfies $\mathbf{f}(\mathbf{x},t) \in C^{r}\left(\mathbb{R}^n\right) \times C(\mathbb{R})$, with $r \geq 1$. The method of LDs was originally defined in \cite{madrid2009} as the computation of the arclength of trajectories starting at the initial condition $\mathbf{x}(t_0)$ when they evolve forward and backward for a fixed integration time $\tau > 0$. In this work we will use the definition of LDs based on the $p$-norm, which was first presented in \cite{lopesino2017} to provide a rigorous mathematical foundation for this technique. Given a value $p \in (0,1]$, fix  forward and backward integration times, $\tau_f > 0$ and $\tau_b > 0$ respectively, the diagnostic is given by:
\begin{equation}
\mathcal{L}_p(\mathbf{x}_{0},t_0,\tau_f,\tau_b) = \int^{t_0+\tau_f}_{t_0-\tau_b} \sum_{k=1}^{n} |f_{k}(\mathbf{x}(t;\mathbf{x}_0),t)|^p \; dt  \;, 
\end{equation}
where $\mathbf{x}_0 = \mathbf{x}(t_0)$ is any inital condition. Notice that it can be decomposed into its backward and forward components:
\begin{equation}
\mathcal{L}_p^{(b)}(\mathbf{x}_{0},t_0,\tau_b) = \int^{t_0}_{t_0-\tau_b} \sum_{k=1}^{n}  |f_{k}(\mathbf{x}(t;\mathbf{x}_0),t)|^p \; dt \quad,\quad \mathcal{L}_p^{(f)}(\mathbf{x}_{0},t_0,\tau_f) = \int^{t_0+\tau_f}_{t_0} \sum_{k=1}^{n} |f_{k}(\mathbf{x}(t;\mathbf{x}_0),t)|^p \; dt \;,
\end{equation}

Although this alternative definition of LDs does not have such an intuitive physical interpretation as that of arclength, it has been shown to provide many advantages. For example, it allows for a rigorous mathematical analysis of the notion of ``singular structures'' and to establish a mathematical connection of this notion to invariant stable and unstable manifolds in phase space. In \cite{lopesino2017} it was shown that forward integration reveals stable manifolds, while backward evolution unveils unstable manifolds, at points where the scalar field $\mathcal{L}_p$ is non-differentiable. Therefore its gradient displays jumps \cite{lopesino2017,demian2017,naik2019a}. This property that the method detects manifolds where the function is non-differentiable, allows for an easy extraction of these geometrical objects by means of edge detection algorithms similar to those used in image processing, such as Sobel filters and many others. Here, we have applied for this purpose the simple approach of computing the norm of the gradient vector of the scalar field, $||\nabla \mathcal{L}_p||$, and also its laplacian $\Delta \mathcal{L}_p$.

We finish this overview of the method by pointing out that there exists many systems where initial conditions may escape to infinity at a very fast rate, or even in finite time (blow up) when we integrate them forward or backward. In these situations, one must adapt LDs so that the underlying geometry of the phase space can still be recovered when applying this tool. To do so, it suffices to implement the strategy of accumulating the value of LDs along a trajectory for the whole integration time chosen, or until the trajectory leaves a fixed region of the phase space, whatever happens first. This procedure has been proved to successfully circumvent the issue caused by escaping trajectories on the values attained by the LD scalar field \cite{GG2019a,GG2019b,katsanikas2020b,katsanikas2020c}.


\section{Results}
\label{sec:sec2}


\subsection{Dissipative Saddle Systems}
\label{subsec:sec1}

\subsubsection{The Linear Saddle}

Let us begin our study of the phase space of dissipative systems by analyzing the simplest case of a linear saddle with different expansion and contraction rates. Consider the dynamical system:
\begin{equation}
	\begin{cases}
		\dot{x} = \lambda x \\[.2cm]
		\dot{y} = -\mu y
	\end{cases}
	\label{saddle}
\end{equation}
where $\lambda > 0$ is the parameter that measures expansion about the origin, and $\mu > 0$ determines the contraction rate. This system is dissipative, since the divergence of the vector field is $\lambda-\mu$, which is not zero when $\lambda \neq \mu$. The analytical solution of this system has the expression:
\begin{equation}
x(t) = x_0 \, e^{\lambda t} \quad,\quad y(t) = y_0 \, e^{-\mu t}
\label{anal_sol}
\end{equation}
where $(x_0,y_0)$ is the initial condition of the trajectory at time $t_0 = 0$. For this system, the origin is a hyperbolic equilibrium point and its stable and unstable manifolds are given by:
\begin{equation}
	\mathcal{W}^{s}(0,0) = \left\{(x,y) \in \mathbb{R}^2 \; \Big| \; x = 0 \right\} \quad,\quad \mathcal{W}^{u}(0,0) = \left\{(x,y) \in \mathbb{R}^2 \; \Big| \; y = 0 \right\}
\end{equation}
If we apply the $p$-norm definition of LDs to Eq. \eqref{saddle} it is straightforward to show that the manifolds of the system are revealed at points where the scalar field $\mathcal{L}_p$ is non-differentiable. This was shown for the first time in \cite{lopesino2017}. Consider a grid of initial conditions $(x_0,y_0)$ in the plane at ime $t_0 = 0$ and integrate them forward and backward for a time interval $\tau_f$ and $\tau_b$ respectively. Then, LDs yields:
\begin{equation}
	\mathcal{L}_p(x_0,y_0,t_0=0,\tau_f,\tau_b) = \int_{-\tau_b}^{\tau_f} |\dot{x}|^p + |\dot{y}|^p \, dt = \dfrac{\lambda^{p-1} |x_0|^p}{p}  \left(e^{p\lambda \tau_f} - e^{-p\lambda \tau_b}\right) - \dfrac{\mu^{p-1} |y_0|^p}{p} \left(e^{-p\mu \tau_f} - e^{p\mu \tau_b}\right)
	\label{estim}
\end{equation}
and this expression is non-differentiable at points $(0,y_0)$, which correspond to the stable manifold, and also at $(x_0,0)$, which is the unstable manifold. Recall that the forward component of LDs highlights the stable manifold whereas the backward contribution reveals the unstable manifold. Notice however that both coordinate axis are singularities, even if we only integrate forward or backward:
\begin{equation}
	\mathcal{L}_p^{(f)}(x_0,y_0,t_0=0,\tau_f,0) = \int_{0}^{\tau_f} |\dot{x}|^p + |\dot{y}|^p \, dt = \dfrac{\lambda^{p-1} |x_0|^p}{p}  \left(e^{p\lambda \tau_f} - 1\right) - \dfrac{\mu^{p-1} |y_0|^p}{p}  \left(e^{-p\mu \tau_f} - 1\right)
\end{equation}
\begin{equation}
	\mathcal{L}_p^{(b)}(x_0,y_0,t_0=0,0,\tau_b) = \int_{-\tau_b}^{0} |\dot{x}|^p + |\dot{y}|^p \, dt = \dfrac{\lambda^{p-1} |x_0|^p}{p}  \left(1 - e^{-p\lambda \tau_b}\right) - \dfrac{\mu^{p-1} |y_0|^p}{p}  \left(1 - e^{p\mu \tau_b}\right)
\end{equation}
If the integration times are large enough, i.e. $\tau_f \gg 1$ and $\tau_b \gg 1$ we get:
\begin{equation}
	\mathcal{L}_p^{(f)}(x_0,y_0,t_0=0,\tau_f,0) \approx \dfrac{\lambda^{p-1} |x_0|^p}{p} e^{p\lambda \tau_f} + \dfrac{\mu^{p-1} |y_0|^p}{p} \approx \dfrac{\lambda^{p-1} |x_0|^p}{p} e^{p\lambda \tau_f}
\end{equation}
\begin{equation}
	 \mathcal{L}_p^{(b)}(x_0,y_0,t_0=0,0,\tau_b) \approx \dfrac{\lambda^{p-1} |x_0|^p}{p} + \dfrac{\mu^{p-1} |y_0|^p}{p} e^{p\mu \tau_b} \approx \dfrac{\mu^{p-1} |y_0|^p}{p} e^{p\mu \tau_b}
\end{equation}
and this result illustrates that, in forward integration, the singularity that dominates is that corresponding to the stable manifold, while backward integration enhances the unstable manifold.

\begin{figure}[htbp]
	\begin{center}
	A)\includegraphics[scale=0.22]{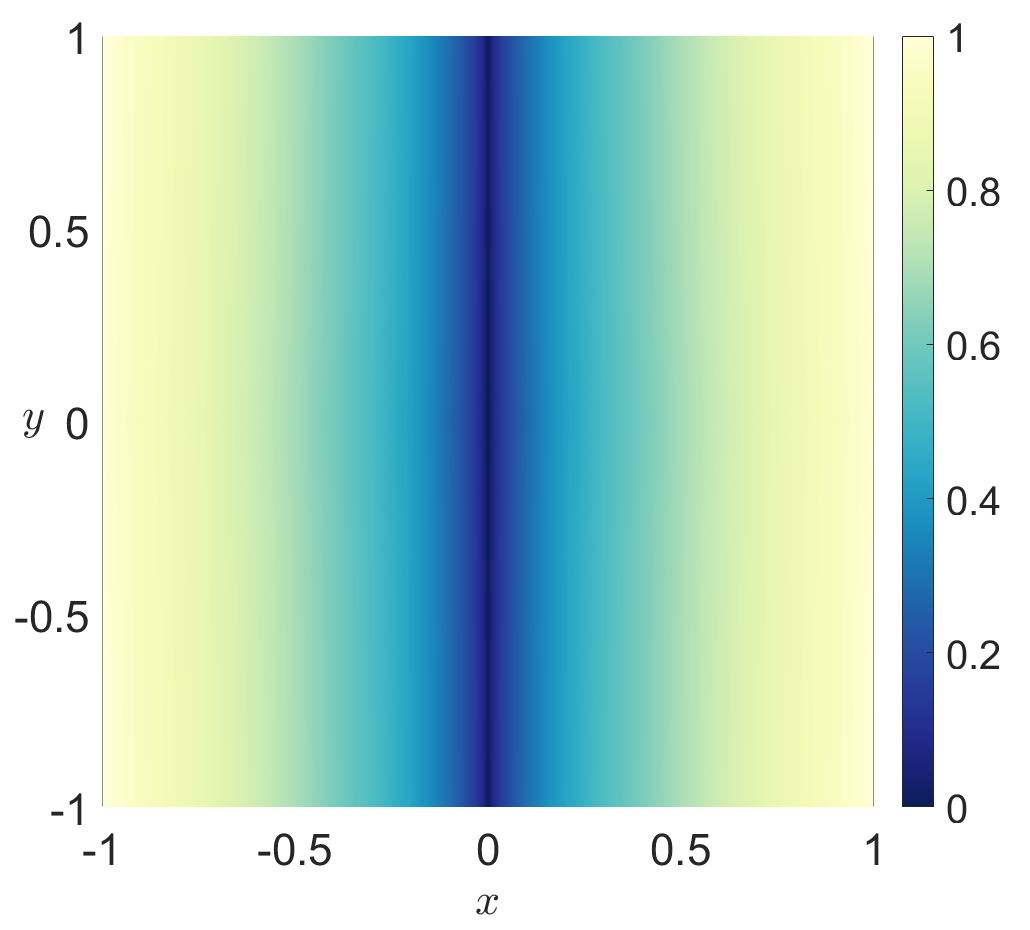} 
	B)\includegraphics[scale=0.22]{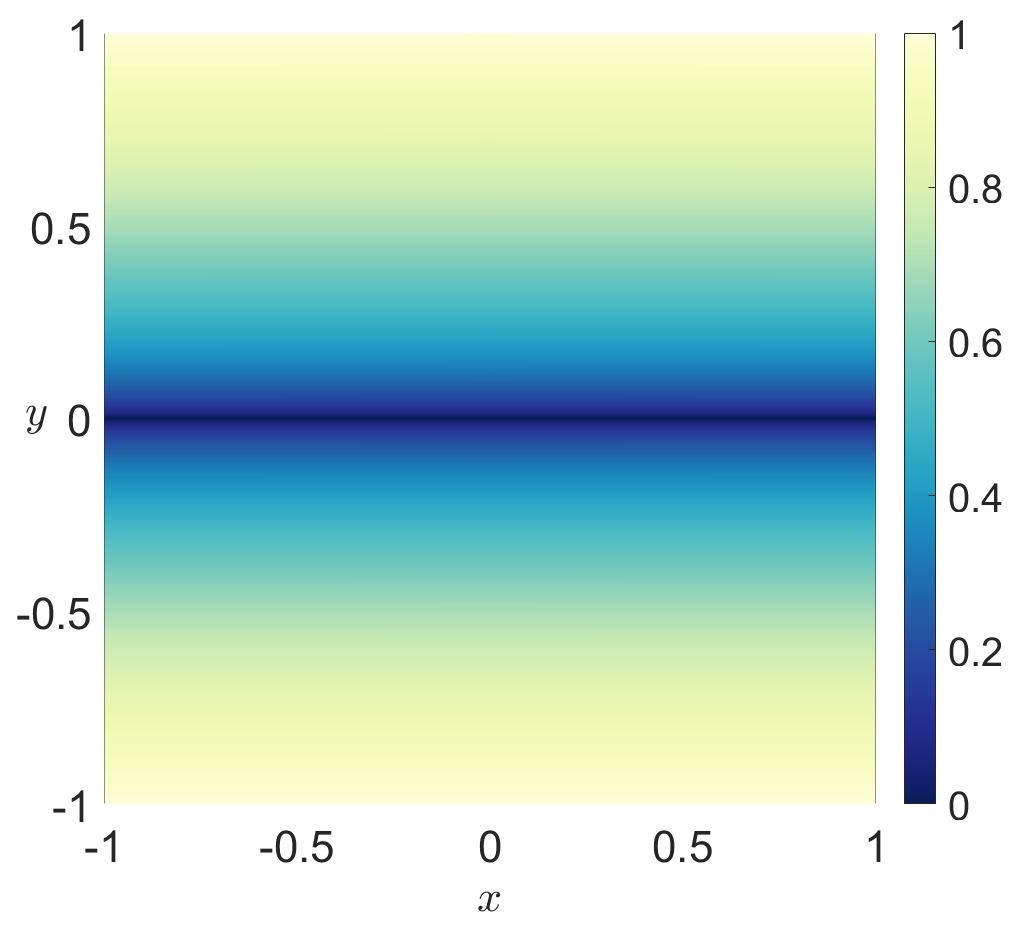} 
	C)\includegraphics[scale=0.22]{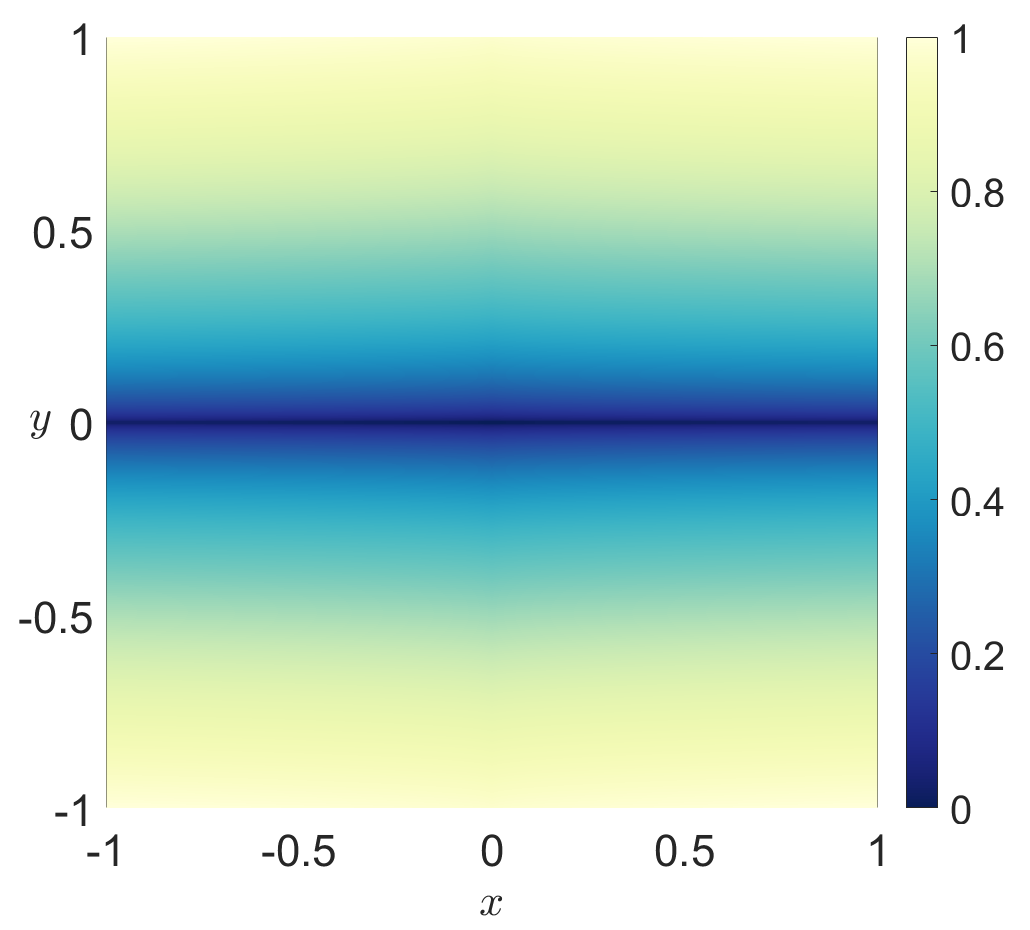} 
	D)\includegraphics[scale=0.22]{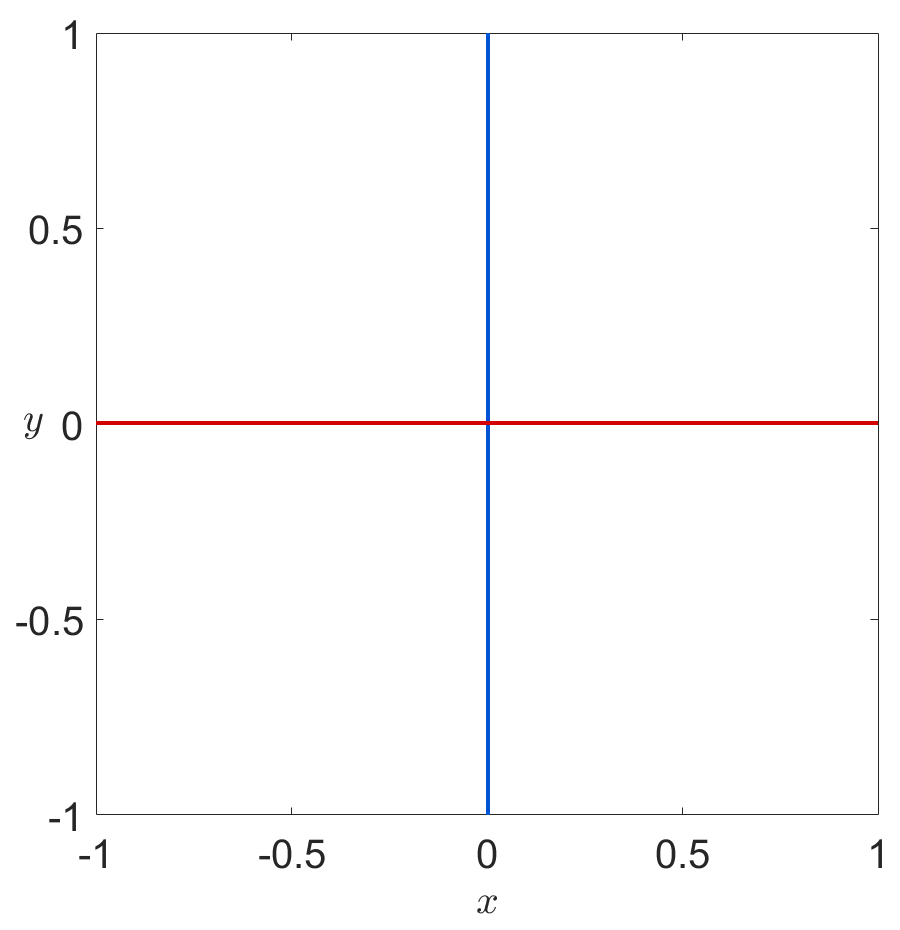}
	\end{center}
	\caption{Phase space of system in Eq. \eqref{saddle} for the model parameters $\lambda =1$, $\mu = 2$ as revealed by LDs using the $p$-norm with $p = 1/2$. A) Forward component of LDs calculated for $\tau_f = 8$. B) Backward component of LDs calculated for $\tau_b = 8$. C) Total LD value obtained by adding the forward and backward contributions. D) Invariant stable (blue) and unstable (red) manifolds extracted from the scalar output of LDs by means of applying the laplacian $\Delta \mathcal{L}_p$.}
	\label{ld_saddle_sameTime}
\end{figure}

However, an issue arises in this system when it comes to visualizing both manifolds simultaneously with LDs in the same picture, since there are two different timescales at play when $\lambda \neq \mu$. There is a competition between the exponential growth, described by $e^{\lambda t}$, and the decay rate identified with $e^{-\mu t}$, and this can obscure phase space structure. We illustrate this behavior in Fig. \ref{ld_saddle_sameTime}, where we have calculated LDs for $p = 1/2$ and the system parameters $\lambda = 1$ and $\mu = 2$, using an integration time $\tau_f = \tau_b = 8$. The values of LDs are normalized by subtracting the minimum and dividing by the difference between the maximum and minimum. In panels A) and B) we depict $\mathcal{L}_{p}^{(f)}$ and $\mathcal{L}_{p}^{(b)}$ respectively, and we can see how the method nicely highlights the location of the stable and unstable manifolds respectively. However, when we plot $\mathcal{L}_{p} = \mathcal{L}_{p}^{(f)} + \mathcal{L}_{p}^{(b)}$ in C), only the unstable manifold is visible and we would like to get both manifolds in the same picture. This happens because $\lambda < \mu$, and therefore, when we compute LD backward in time for an initial condition $(0,z_0)$ on the stable manifold using $\tau_b = \tau$, and integrate a symmetric initial condition $(z_0,0)$ on the unstable manifold forward for $\tau_f = \tau$, this yields:
\begin{equation}
	\mathcal{L}_p^{(f)}(z_0,0,t_0=0,\tau,0) = \dfrac{\lambda^{p-1} |z_0|^p}{p}  \left(e^{p\lambda \tau} - 1\right) \quad,\quad \mathcal{L}_p^{(b)}(0,z_0,t_0=0,0,\tau) = \dfrac{\mu^{p-1} |z_0|^p}{p}  \left(e^{p\mu \tau} - 1\right)
\end{equation}
If we divide these two quantities we can compare their magnitude:
\begin{equation}
	\dfrac{\mathcal{L}_p^{(f)}(z_0,0,t_0=0,\tau,0)}{\mathcal{L}_p^{(b)}(0,z_0,t_0=0,0,\tau)} = \left(\dfrac{\lambda}{\mu}\right)^{p-1} \dfrac{e^{p\lambda \tau} - 1}{e^{p\mu \tau} - 1}
\end{equation}
and, whenever $\tau \gg 1$, this gives:
\begin{equation}
	\dfrac{\mathcal{L}_p^{(f)}(z_0,0,t_0=0,\tau,0)}{\mathcal{L}_p^{(b)}(0,z_0,t_0=0,0,\tau)} \approx \left(\dfrac{\lambda}{\mu}\right)^{p-1} e^{p\left(\lambda-\mu\right) \tau}
\end{equation}
Therefore, if $\lambda < \mu$, it is clear that $\mathcal{L}_p^{(b)}$ dominates, while if $\lambda > \mu$ then $\mathcal{L}_p^{(f)}$ becomes the dominating term. This explains why in Fig. \ref{ld_saddle_sameTime} C) the values of LDs on the unstable manifold dominate those attained on the stable manifold, and hence we cannot see the stable manifold when we plot the quantity $\mathcal{L}_p^{(f)} + \mathcal{L}_p^{(b)}$. Nevertheless, we would like to point out that we can easily extract the location of the points corresponding to the stable and unstable manifolds by means of computing the gradient of $\mathcal{L}^(f)$ and  $\mathcal{L}^(b)$ respectively, and overlay both results in the same plot. This is shown in Fig. \ref{ld_saddle_sameTime} D) where we have used the standard coloring scheme widely applied in the literature, that is, red for the unstable manifold and blue for the stable manifold.

\begin{figure}[htbp]
	\begin{center}
		A)\includegraphics[scale=0.2]{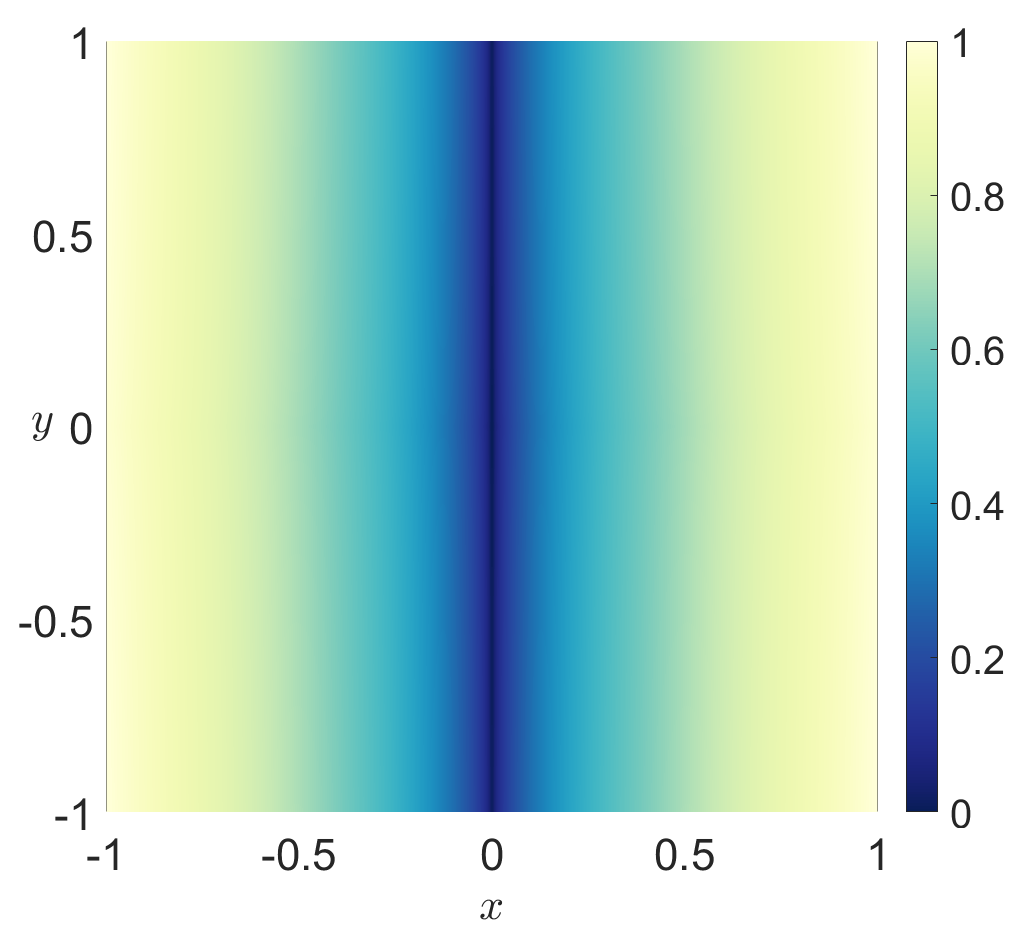} 
		B)\includegraphics[scale=0.2]{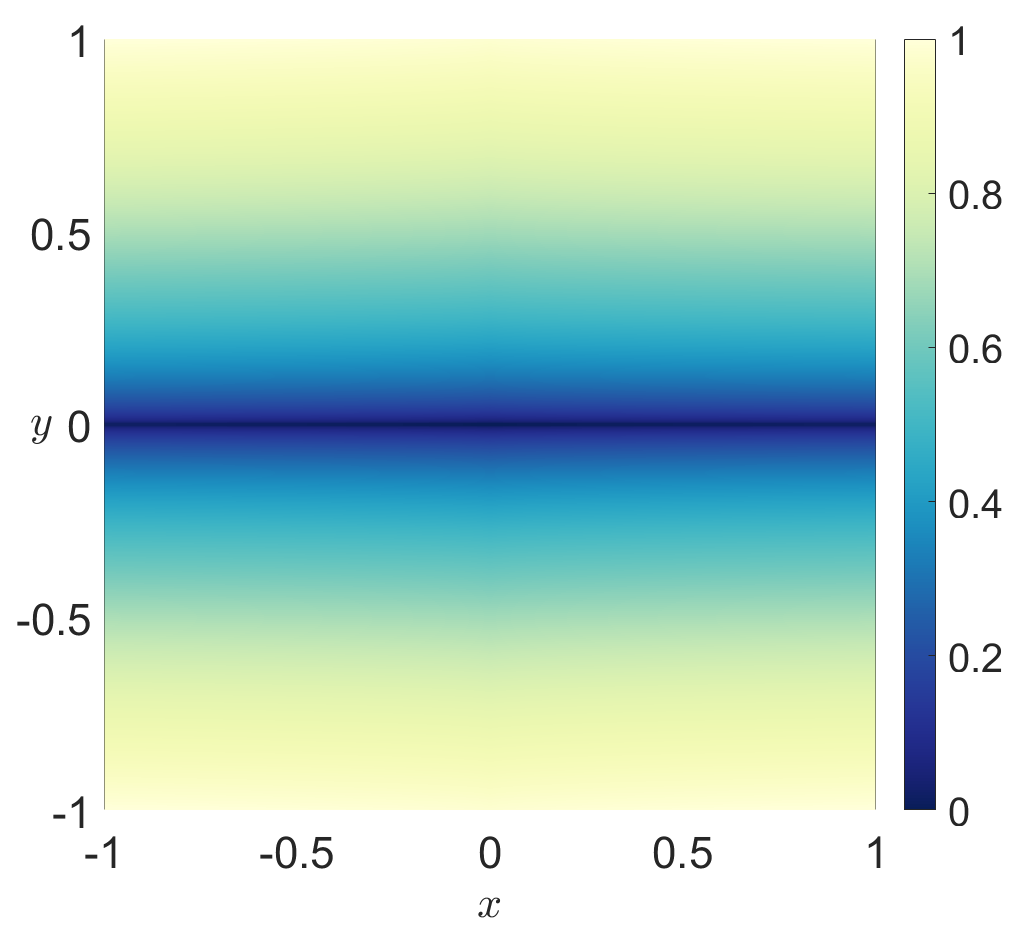} 
		C)\includegraphics[scale=0.2]{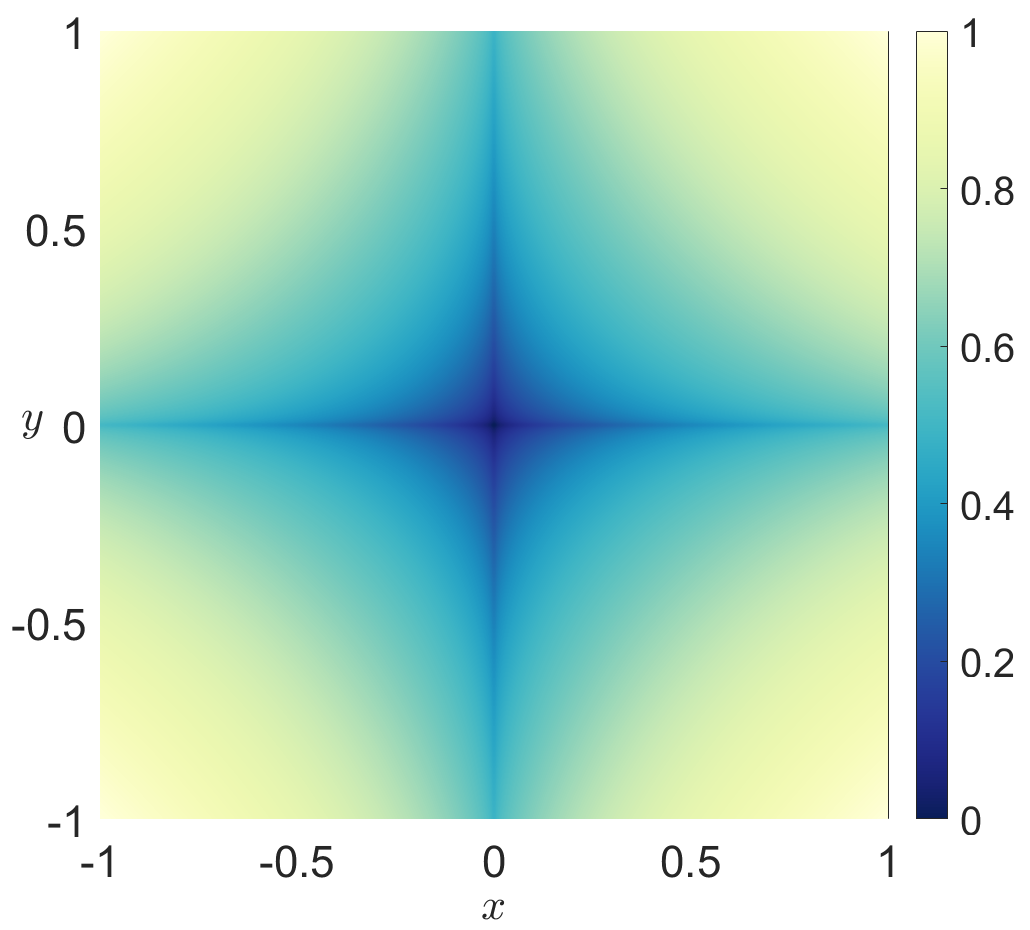} 
	\end{center}
	\caption{Phase space of system in Eq. \eqref{saddle} for the model parameters $\lambda = 1$, $\mu = 2$ as revealed by LDs using the $p$-norm with $p = 1/2$. A) Forward component of LDs calculated for $\tau_f = 8$. B) Backward component of LDs calculated for $\tau_b = 4.346$. C) Total LD value obtained by adding the forward and backward contributions.}
	\label{ld_saddle_difTime}
\end{figure}

This visualization problem in the total LD scalar field was first discussed in detail in \cite{lopesino2017}, where it was shown that one can circumvent this issue by changing the value of $p$ that defines the $p$-norm used to compute LDs. The argument that was developed in that paper assumed that $\tau_f = \tau_b$, that is, all initial conditions are integrated for the same time forwards and backwards. However, it seems more natural to proceed as follows. Fix a value $p \in (0,1]$ to define the $p$-norm of LDs, and then adjust the forward and backward integration times $\tau_f$ and $\tau_b$ so that the time scales of the problem compensate and the contributions of the forward and backward components of LDs become comparable. We will determine here a condition for $\tau_f$ and $\tau_b$ that solves the issue of depicting all the phase space structures simultaneously in the same simulation. Take the initial conditions $(z_0,0)$ and $(0,z_0)$ on the unstable and stable manifolds respectively, and integrate the one on $\mathcal{W}^{u}$ forward in time for $\tau_f$, and the one located on $\mathcal{W}^s$ backward in time for $\tau_b$. Impose that the values of LDs obtained are comparable, that is:
\begin{equation}
	\mathcal{L}_p^{(f)}(z_0,0,t_0=0,\tau_f,0) \approx \mathcal{L}_p^{(b)}(0,z_0,t_0=0,0,\tau_b) \quad \Leftrightarrow \quad \dfrac{\lambda^{p-1} |z_0|^p}{p} \left(e^{p\lambda \tau_f} - 1\right) \approx \dfrac{\mu^{p-1} |z_0|^p}{p}  \left(e^{p\mu \tau_b} - 1\right)
\end{equation}
Taking $\tau_f,\tau_b \gg 1$ and simplifying, this gives:
\begin{equation}
	\left(\dfrac{\mu}{\lambda}\right)^{p-1} \approx e^{p\left(\lambda\tau_f - \mu\tau_b\right)} \quad \Leftrightarrow \quad \tau_b \approx \dfrac{\lambda}{\mu} \tau_f + \dfrac{1-p}{\mu p} \ln\left(\dfrac{\mu}{\lambda}\right)
	\label{intTimes_rel}
\end{equation}
Therefore, we have found a simple formula that determines how to set the forward and backward integration times in terms of the Lyapunov exponents of the saddle system. In particular, for the case $\lambda = \mu$ the formula indicates that forward and backward integration times should be chosen to be equal. In order to validate that this approach helps improving the visualization of the total LD field, we calculate LDs for the model parameters $\lambda = 1$, $\mu = 2$ and use a forward integration time of $\tau_f = 8$, and a corresponding backward integration time of $\tau_b = 4.346$ determined from the formula in Eq. \eqref{intTimes_rel}. We display the total LD function in Fig. \ref{ld_saddle_difTime} C), and this confirms that we acn adjust the forward and backward integration accordingly to obtain a complete picture of the phase space of the system, despite the different expansion and contraction rates, and hence the timescale difference is captured by the method.

\subsubsection{A Nonlinear Saddle Example}

Next, we move on to analyze the nonlinear system with a saddle point at the origin:
\begin{equation}
	\begin{cases}
		\dot{x} = \mu x \\[.2cm]
		\dot{y} = \lambda \left(y - x^2\right)
	\end{cases}
	\label{nonlin_saddle}
\end{equation}
and focus our analysis in the case where $\lambda = -2$ and $\mu = 1$. It is a simple exercise to show that this system can be reduced to a non-autonomous (in $x$) first order linear ODE that can be solved by means of an integrating factor. In fact, the solution is:
\begin{equation}
   x(t) = x_0 e^{t} \quad,\quad y(t) = \frac{1}{2}x_0^{2}e^{2t}+\left(y_0-\frac{1}{2}x_0^{2}\right)e^{-2t}
\end{equation}
where $(x_0,y_0)$ is the initial condition of the trajectory at time $t_0 = 0$. This system has an unstable equilibrium point at the origin and its stable and unstable manifolds are given by:
\begin{equation}
    \mathcal{W}^{s}(0,0) = \left\{(x,y) \in \mathbb{R}^2 \; \Big| \; x = 0\right\} \quad,\quad 	
	 \mathcal{W}^{u}(0,0) = \left\{(x,y) \in \mathbb{R}^2 \; \Big| \; y = \dfrac{x^2}{2}\right\}
	\label{nonlin_saddle_mani}
\end{equation}

\begin{figure}[htbp]
\begin{center}
A)\includegraphics[scale=0.2]{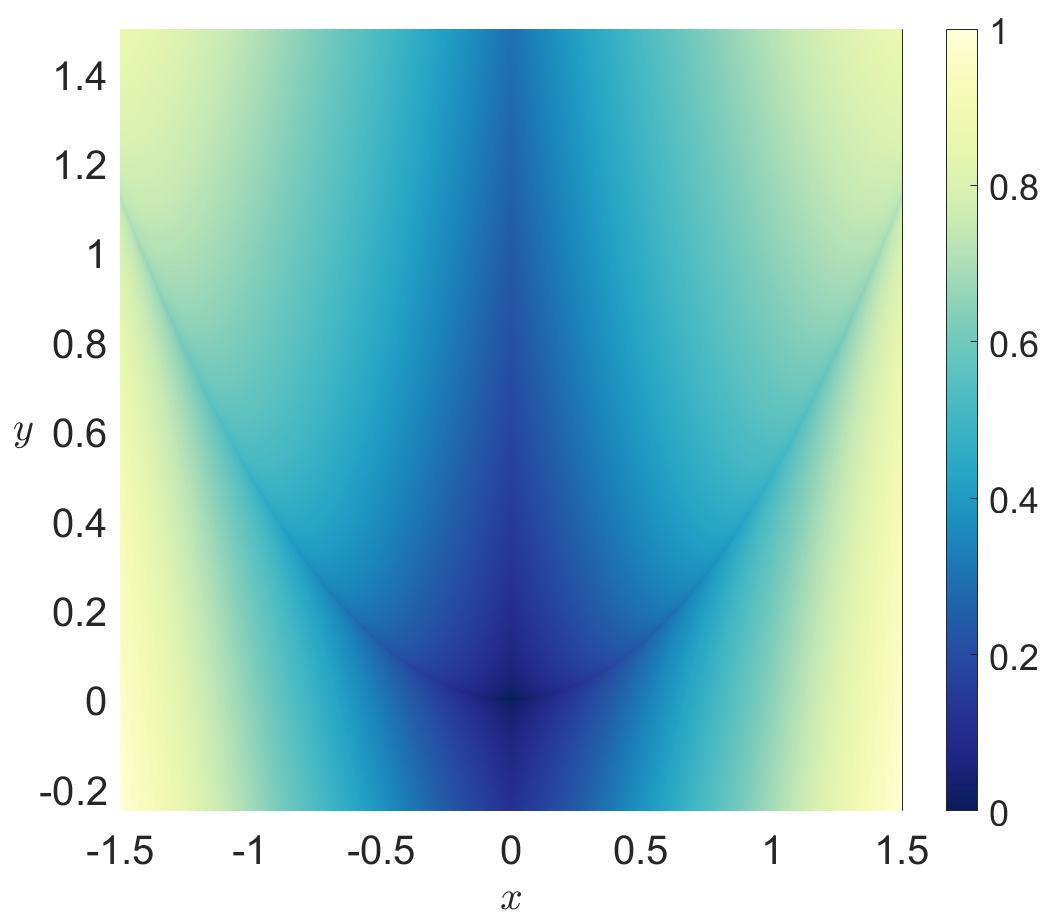} 
B)\includegraphics[scale=0.2]{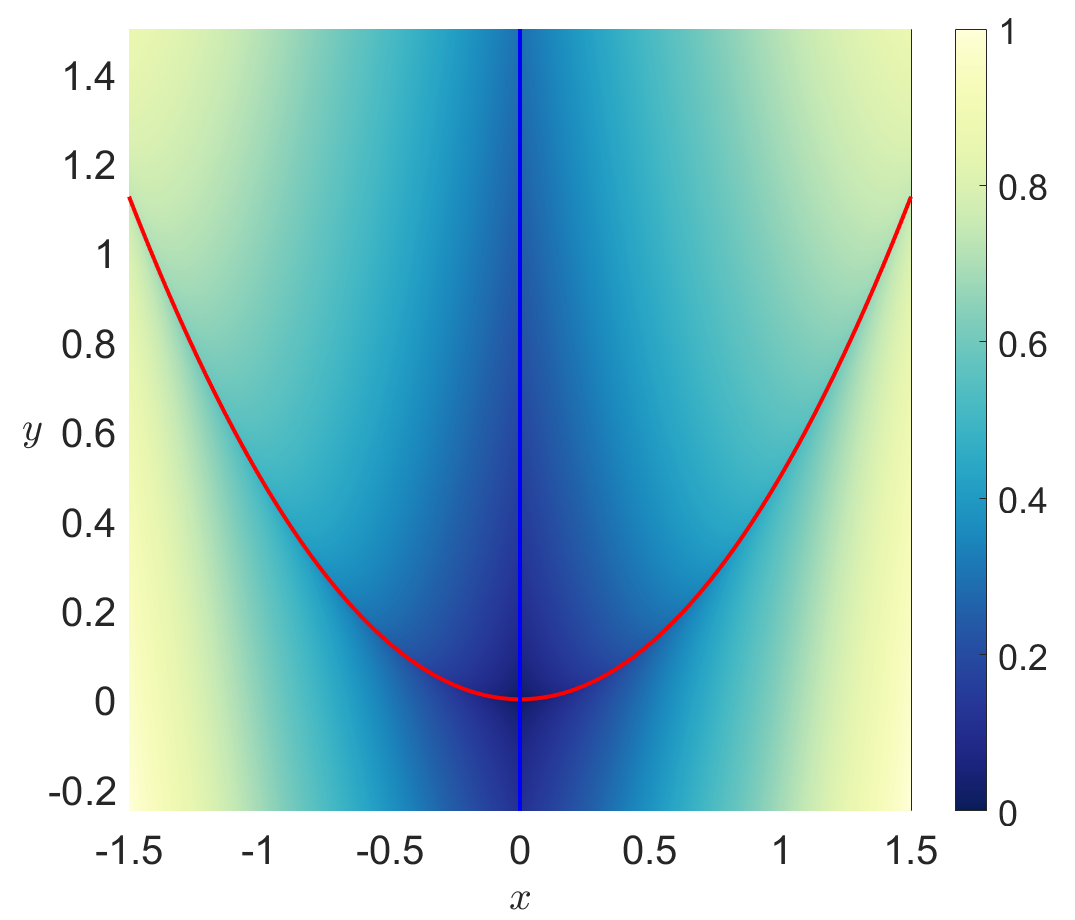} 
C)\includegraphics[scale=0.2]{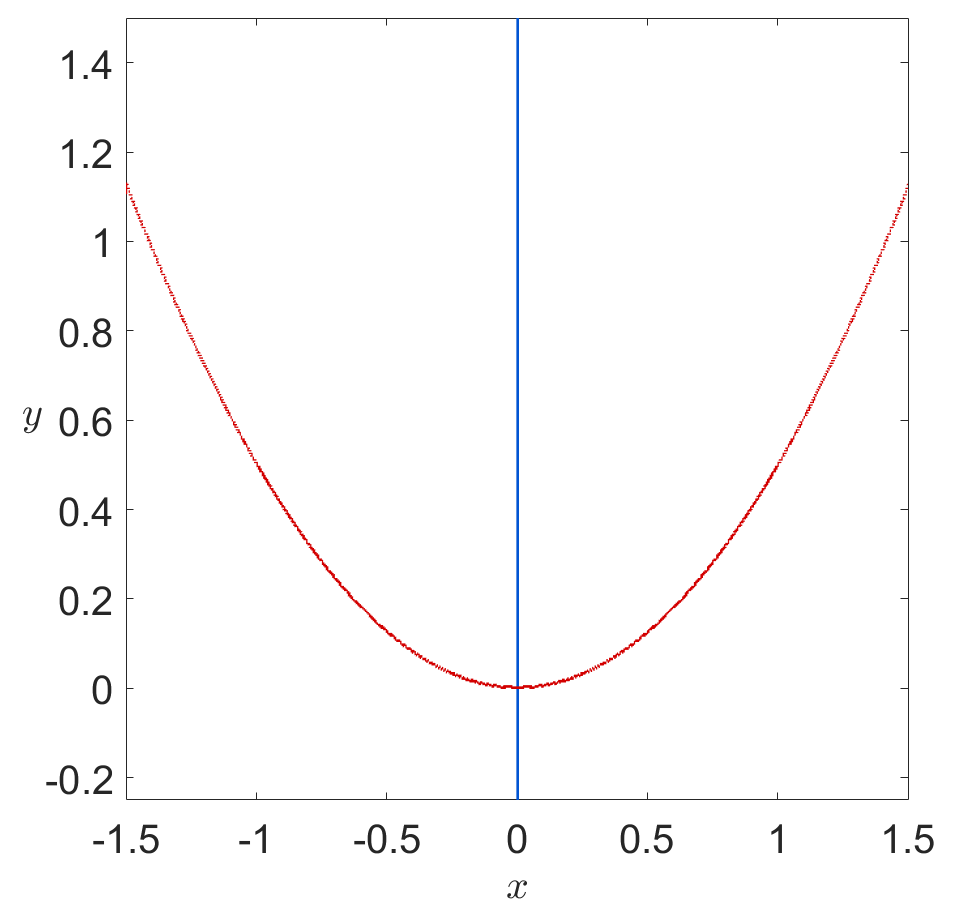}  
\end{center}
\caption{Phase space of system in Eq. \eqref{nonlin_saddle} for the model parameters $\lambda = -2$, $\mu = 1$ as revealed by LDs using the $p$-norm with $p = 1/2$. A) Total LD calculated by integrating forward and backward for $\tau_f = 26$ and $\tau_b = 25$, respectively. B) Analytical comparison of the LD scalar values with the stable (blue) and unstable (red) manifolds given by Eq. \eqref{nonlin_saddle_mani} for the saddle point at the origin of the system described in Eq. \eqref{nonlin_saddle}. C) Extraction of the stable and unstable manifolds from the LD values shown in panel A) by means applying a Sobel-type filter, that is, using the gradient of the LD function, $||\nabla \mathcal{L}_p||$, in order to detect the singular features (image edges) of the scalar output provided by LDs.}
	\label{ld_nonlin}
\end{figure}

We will show that the method of LDs detects the stable manifold when we integrate forward, whereas the unstable manifold is highlighted by the backward integration. For this computation, note that we can not obtain a closed-form expression for the integrals. Hence, we need to perform an asymptotic approximation similar to the one explained in \cite{lopesino2017}. Consider first the forward component of LDs:
\begin{equation}
	\mathcal{L}_p^{(f)}(x_0,y_0,t_0=0,\tau_f,0) = \int_{0}^{\tau_f}|\dot{x}|^p + |\dot{y}|^p \, dt \sim \dfrac{|x_0|^p}{p} \left(e^{p\tau_f} - 1\right) + \dfrac{|x_0|^{2p}}{2p} \left(e^{2p\tau_f} - 1\right)
\end{equation}
As we can see, this expression is non-differentiable when $x_0=0$, which coincides with the stable manifold of the system. Similarly, the backward component of LDs can be approximated by:
\begin{equation}
	\mathcal{L}_p^{(b)}(x_0,y_0,t_0=0,0,\tau_b) = \int_{-\tau_b}^{0}|\dot{x}|^p + |\dot{y}|^p \, dt \sim \dfrac{|x_0|^p}{p} \left(1-e^{-p\tau_b}\right) + \dfrac{2^{p-1}}{p}\left|y_0-\dfrac{1}{2}x_0^2\right|^p \left(e^{2p\tau_b} - 1\right)
\end{equation}
Since $\tau_b\gg 1$, the leading order singularity in the gradient of the backward component appears in the term that contains $|y_0-x_0^2/2|$. Therefore, it reveals the equation of the unstable manifold of the system. This analytic derivation demonstrates that the method successfully highlights the stable and unstable manifolds of the saddle point at the origin. We illustrate numerically in Fig. \ref{ld_nonlin} that the manifolds are detected. Indeed, panel A) displays the total LD values, and in B) we overlay the LD field with with the stable (blue) and unstable (red) manifolds described in Eq. \eqref{nonlin_saddle_mani}. As a validation, we extract the manifold locations by calculating the gradient of the LD function, $||\nabla \mathcal{L}_p||$, which allows us to detect the singular features (image edges) of the scalar output provided by LDs.

\subsection{Detection of Limit Cycles}
\label{subsec:sec2}

\subsubsection{Andronov-Hopf Bifurcations}

Our next goal is to show how Lagrangian descriptors can be used to detect the presence of limit cycles. We look first at the system described by the 2D normal form for a supercritical Andronov-Hopf bifurcation \cite{marsden76}, which is given by the following equations:
\begin{equation}
	\begin{cases}
	\dot{x} = \beta x - y - \sigma x \left(x^2 + y^2\right) \\[.2cm]
	\dot{y} = x + \beta y - \sigma y \left(x^2 + y^2\right)
	\end{cases}
\label{hopf_nf}
\end{equation}
We will consider here the case where $\sigma > 0$. If one writes Eq. \eqref{hopf_nf} in polar coordinates, it is straightforward to show that the system becomes:
\begin{equation}
	\begin{cases}
	\dot{r} = r(\beta - \sigma r^2) \\[.2cm]
    \dot{\theta} = 1
	\end{cases}
\label{hopf_polar}
\end{equation}
and therefore, it is clear that when $\beta > 0$ the system has a stable limit cycle at the circle of radius $\sqrt{\beta/\sigma}$, while for $\beta < 0$ the limit cycle disappears, and the origin becomes an asymptotically stable equilibrium point. One needs to be careful about the integration of initial conditions, because in this case, many trajectories will blow up in finite time. This is nicely exemplified if we solve analytically the system in Eq. \eqref{hopf_polar} for the situation where $\beta = 0$. For this parameter value, the analytical solution of the system \eqref{hopf_polar} is:
\begin{equation}
	r(t) = \frac{r_0}{\sqrt{2\sigma tr^2_0+1}} \quad,\quad \theta(t) = \theta_0 + t
\label{analSolLimCy}
\end{equation}
where $(r_0,\theta_0)$ is the initial condition of the trajectory at time $t_0 = 0$. We can clearly see that this solution blows up at time $t = -1/\left(2\sigma r_0^2\right) < 0$, so that one would run into trouble when computing trajectories backward in time. This issue is also present in the solutions to the dynamical system in Eq. \eqref{hopf_polar} for $\beta \neq 0$, and can be easily checked by integrating the ODE using separation of variables.

\begin{figure}[htbp]
	\begin{center}
	A)\includegraphics[scale=0.19]{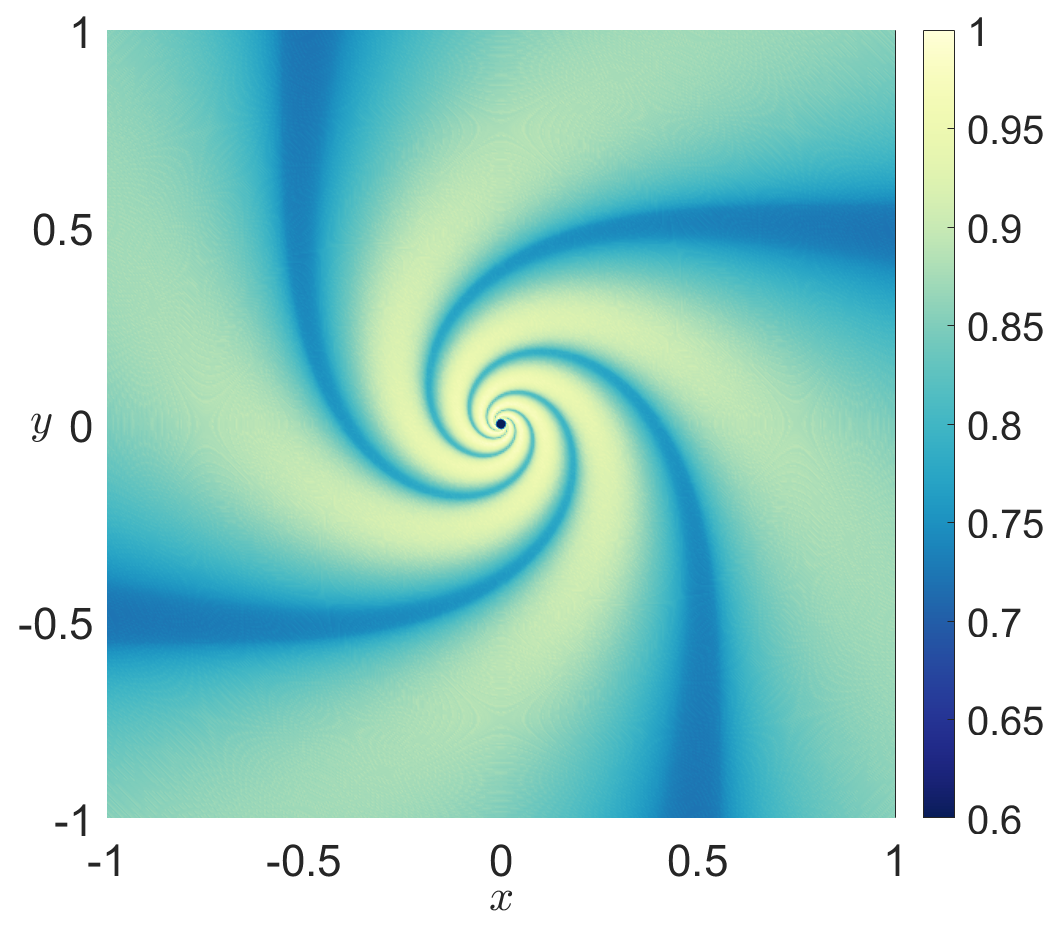} 
	B)\includegraphics[scale=0.19]{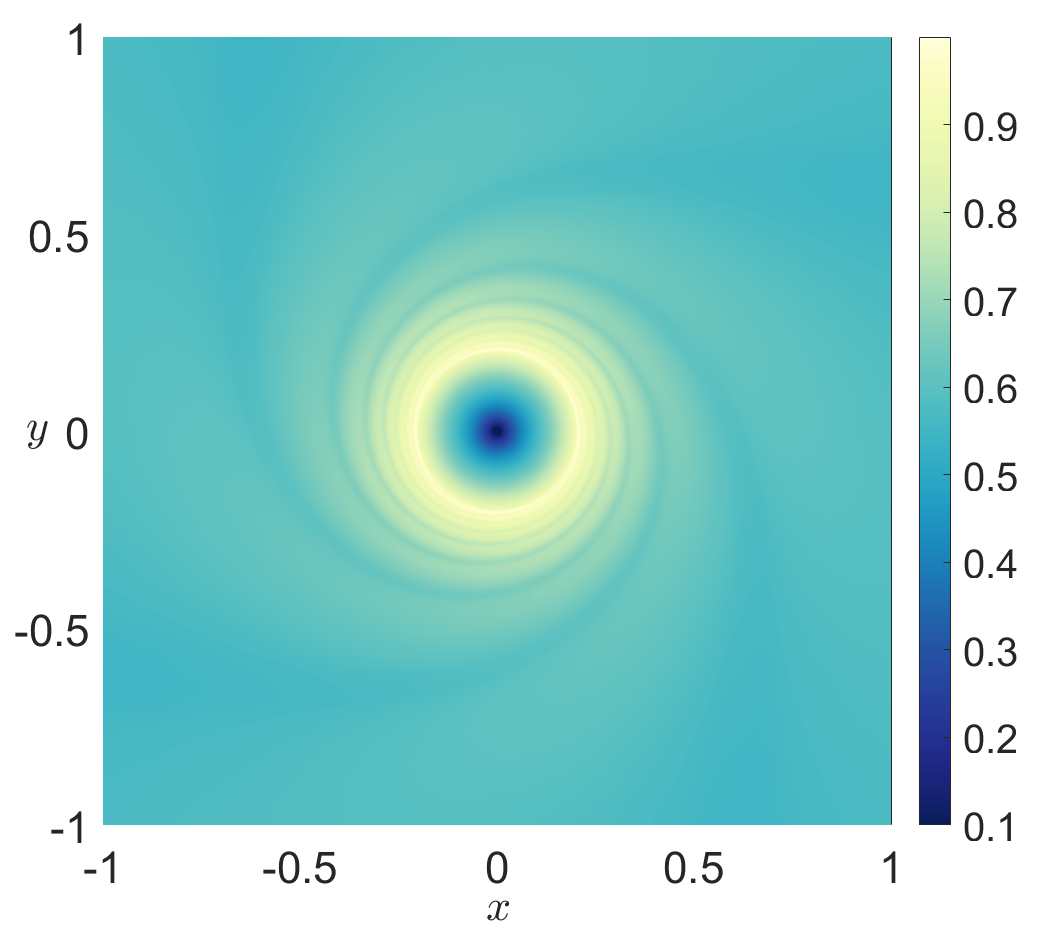} 
	C)\includegraphics[scale=0.19]{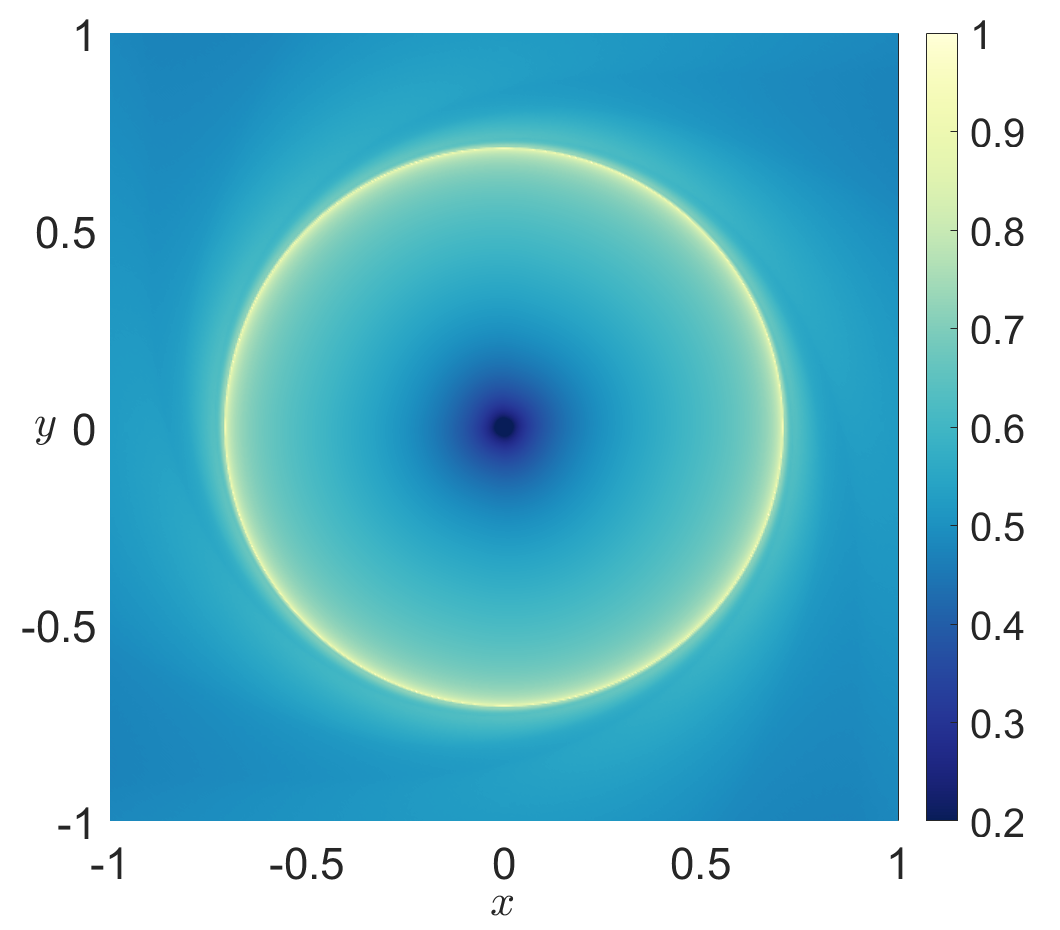} 
	\end{center}
	\caption{Phase space of the dynamical system in Eq. \eqref{hopf_nf}, as revealed by applying the $p$-norm definition of LDs with $p = 1/2$ using an integration time of $\tau_f = \tau_b = 8$, for different values of the model parameters. A) $\beta = -0.5$ and $\sigma = 1$. B) $\beta = 0$ and $\sigma = 1$. C) $\beta = 0.5$ and $\sigma = 1$. Notice the false positive displayed by the method in panel B), since in this case ($\beta = 0$) there is no limit cycle in the system. More details about this issue are given in the main text.}
	\label{ld_andronov}
\end{figure}

Our goal now is to apply LDs to the system in Eq. \eqref{hopf_nf} in order to reveal its phase space. Since there is an issue of trajectories blowing up in finite time we will use the strategy of stopping those trajectories that escape a circular region of radius $R = 4$ centered at the origin. In this way, we avoid problems when computing LDs. In order to reveal the phase space of this system we compute LDs for $\tau_f = \tau_b = 8$ using variable time integration, so that the initial conditions that escape the region of the phase space delimited by a circle of radius $R = 4$ centered about the origin are stopped. We run experiments for the mode parameter values $\sigma = 1$, and $\beta = -0.5,\, 0 ,\, 0.5$ and integrate trajectories forwards and backwards for a time $\tau_f = \tau_b = 8$ or until they leave the circle of radius $R = 4$, what happens first. We present the results of this simulations in Fig. \ref{ld_andronov}. In panel A), which corresponds to $\beta = -0.5$, the origin is a stable focus, and the LD values highlight that dynamical behavior. In B), for $\beta = 0$, this case is the critical value of the parameter where the system undergoes a supercritical Andronov-Hopf bifurcation. The origin is still a weakly asymptotically stable focus and a limit cycle is about to emerge. However, the LD scalar field seems to be indicating the existence of a limit cycle at points where its values are very large. This is not correct, and we will comment on this issue later. Lastly, in panel C) we show how LDs nicely capture the location of the limit cycle present in the phase space of the system for $\beta = 0.5$.

At this point we would like to discuss the false positive given by LDs for the case where $\beta = 0$. First, we calculate the forward component of LDs:
\begin{equation}
	\mathcal{L}_p^{(f)}(r_0,\theta_0,t_0=0,\tau_f,0) = \int_{0}^{\tau_f} |\dot{r}|^p + |\dot{\theta}|^p \, dt =
	\begin{cases}
	\dfrac{\sigma^{p-1}r_0^{3p-2}}{2-3p}\big[(2\sigma\tau_f r_0^2+1)^{1-3p/2}-1\big]+\tau_f \,,
	\quad\mbox{if }p\neq2/3 \\[0.6cm]
	\dfrac{\sigma^{-1/3}}{2}\ln(2\sigma\tau_f r_0^2+1) + \tau_f \,,
	\quad\mbox{if }p=2/3
	\end{cases}
\end{equation}
Thus, it is clear that there are no singularities in the forward integration. However, if we integrate backward, we obtain that:
\begin{equation}
	\mathcal{L}_p^{(b)}(r_0,\theta_0,t_0=0,0,\tau_b) = \int_{-\tau_b}^{0} |\dot{r}|^p + |\dot{\theta}|^p \, dt =
	\begin{cases}
	\dfrac{\sigma^{p-1}r_0^{3p-2}}{2-3p}\big[1-(-2\sigma\tau_b r_0^2+1)^{1-3p/2}\big] + \tau_b \,,
	\quad\mbox{if }p\neq2/3 \\[0.6cm]
	-\dfrac{\sigma^{-1/3}}{2}\ln(-2\sigma\tau_b r_0^2+1) + \tau_b \,,
	\quad\mbox{if }p=2/3
	\end{cases}
\end{equation}
and these expressions are non-differentiable when $r_0 = 1/\sqrt{2\sigma\tau_b}$. It is important to observe that this is a false positive, since we know that the system in Eq. \eqref{hopf_polar} does not have a limit cycle for $\beta=0$. Notice that this false positive shrinks as the integration time gets large, and it will tend to disappear when $\tau_b$ goes to infinity. This situation could be interpreted in the sense that the method is identifying the critical value of $\beta$ for which the system's phase space undergoes a bifurcation, and a limit cycle is about to be born. This potential issue of a false positive shows the importance of carefully checking the output provided by the method when analyzing the dynamics of a system. In summary, if a structure tends to disappear as the integration time gets large, it is reasonable to expect that this structure does not have any dynamical significance, although this should be always verified by other means before drawing any further conclusions.  

\subsubsection{The van der Pol Oscillator}

We focus our attention now on the classical van der Pol oscillator system, which is a paradigmatic example of a system displaying limit cycle behavior \cite{vdp1922,vdp1926,strogatz,meiss2017}. This dynamical system is described by the following second order differential equation:
\begin{equation}
	\ddot{x} + \mu \left(x^2-1\right) \dot{x}  + x = 0
\end{equation}
where $\mu \geq 0$. As a system of first order ODEs, this can be rewritten as:
\begin{equation}
	\begin{cases}
		\dot{x} = y \\[.2cm]
		\dot{y} = - x + \mu  \left(1-x^2\right) y
	\end{cases}
	\label{vdp_sys}
\end{equation}

\begin{figure}[htbp]
	\begin{center}
	A)\includegraphics[scale=0.26]{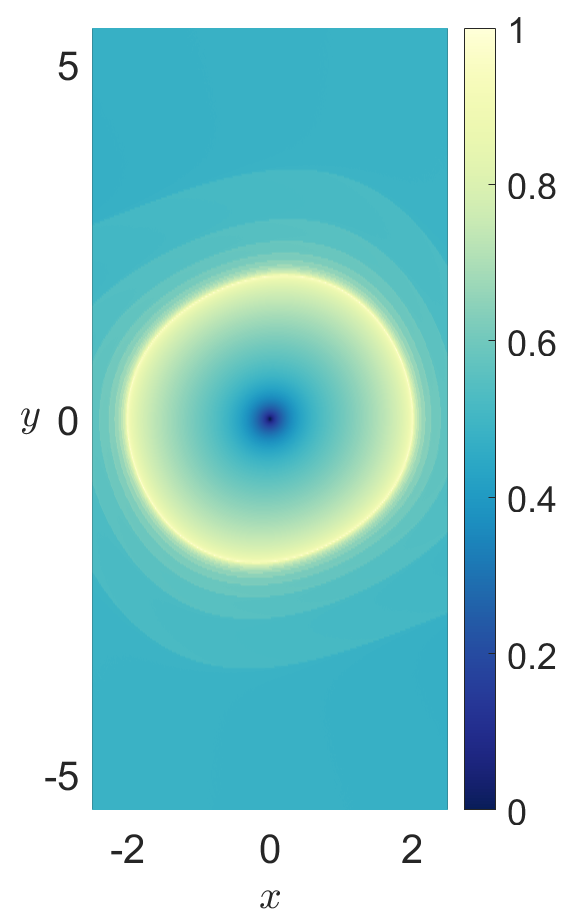} 
	B)\includegraphics[scale=0.26]{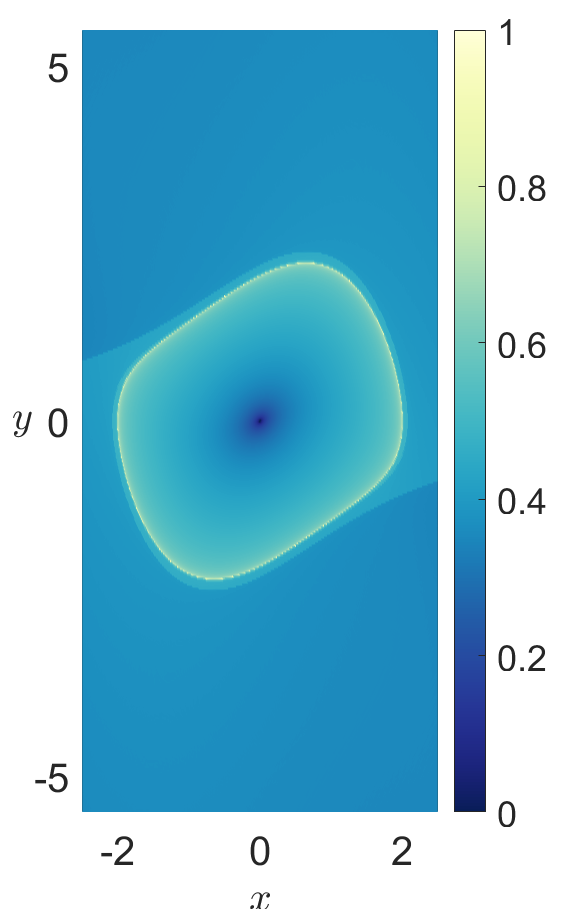} 
	C)\includegraphics[scale=0.26]{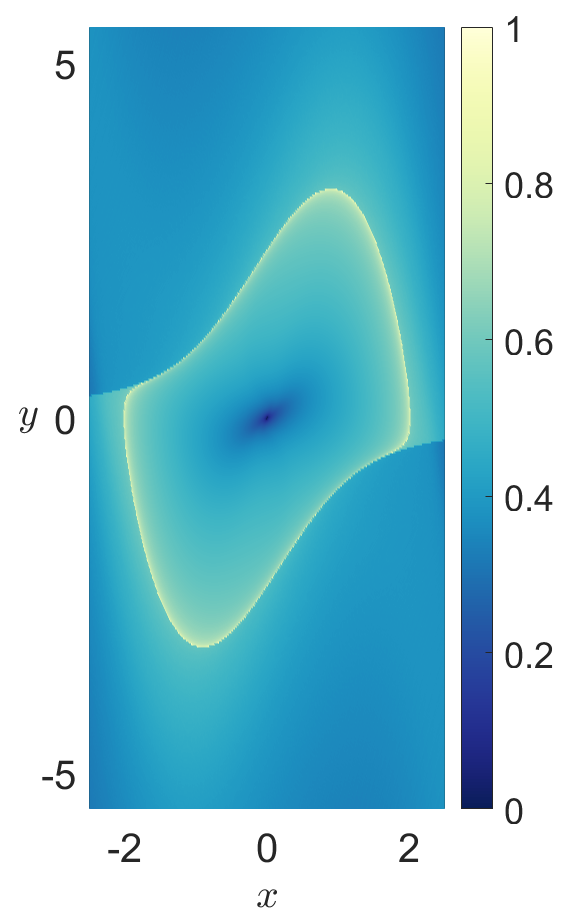} 
	D)\includegraphics[scale=0.26]{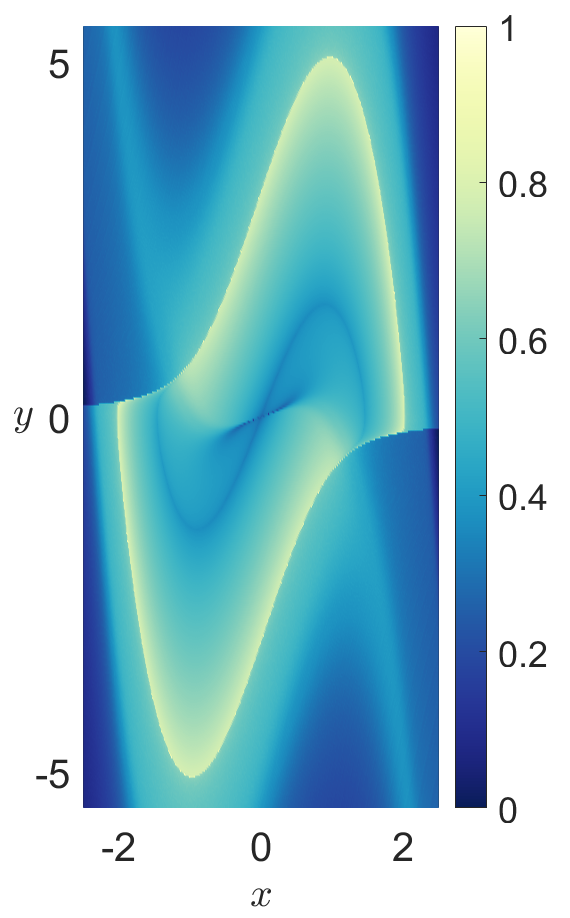} 
	\end{center}
	\caption{Phase space of the van der Pol oscillator in Eq. \eqref{vdp_sys}, as revealed by applying the $p$-norm definition of LDs with $p = 1/2$ using an integration time of $\tau_f = \tau_b = 50$, for different values of the model parameters. A) $\mu = 0.1$. B) $\mu = 0.5$. C) $\mu = 1.5$. D) $\mu = 3$. In order to avoid the issue of trajectories escaping in finite time, we stop them when they exit a circle of radius $R = 20$ centered at the origin.}
	\label{ld_vdp_lowmu}
\end{figure}

In order to compute LDs for this system we will use again the strategy of stopping trajectories if they escape a certain phase space region, since backward integration also has blow up issues in finite time. For this purpose we set a circle of radius $R = 20$ about the origin, and compute LDs forward and backward for $\tau_f = \tau_b = 50$, or until the trajectory leaves this circle, what happens first. To illustrate how the limit cycle of the system gets larger and its shape is distorted as $\mu$ is increased, we have simulated the system for $\mu = 0.1,\, 0.5,\, 1.5,\, 3$. Results of the total LD field are shown in Fig. \ref{ld_vdp_lowmu}, and the panels from left to right correspond to the cases from lower to higher values of $\mu$. It is clear from these plots that the method successfully recovers the geometry and location of the limit cycle for the van der Pol oscillator, and how it deforms as the model parameter value is increased.

\subsection{Analysis of Systems with Slow Manifolds}
\label{subsec:sec3}

\begin{figure}[htbp]
	\begin{center}
		A)\includegraphics[scale=0.19]{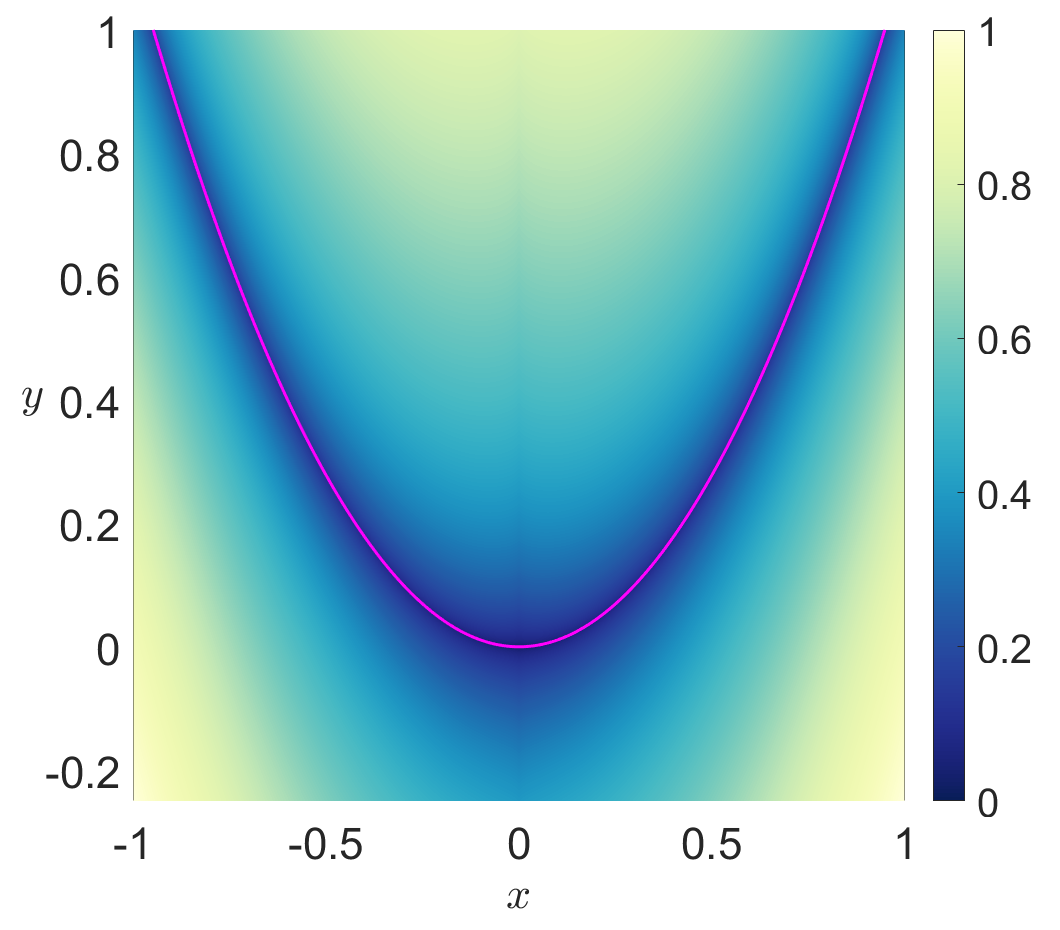} 
		B)\includegraphics[scale=0.19]{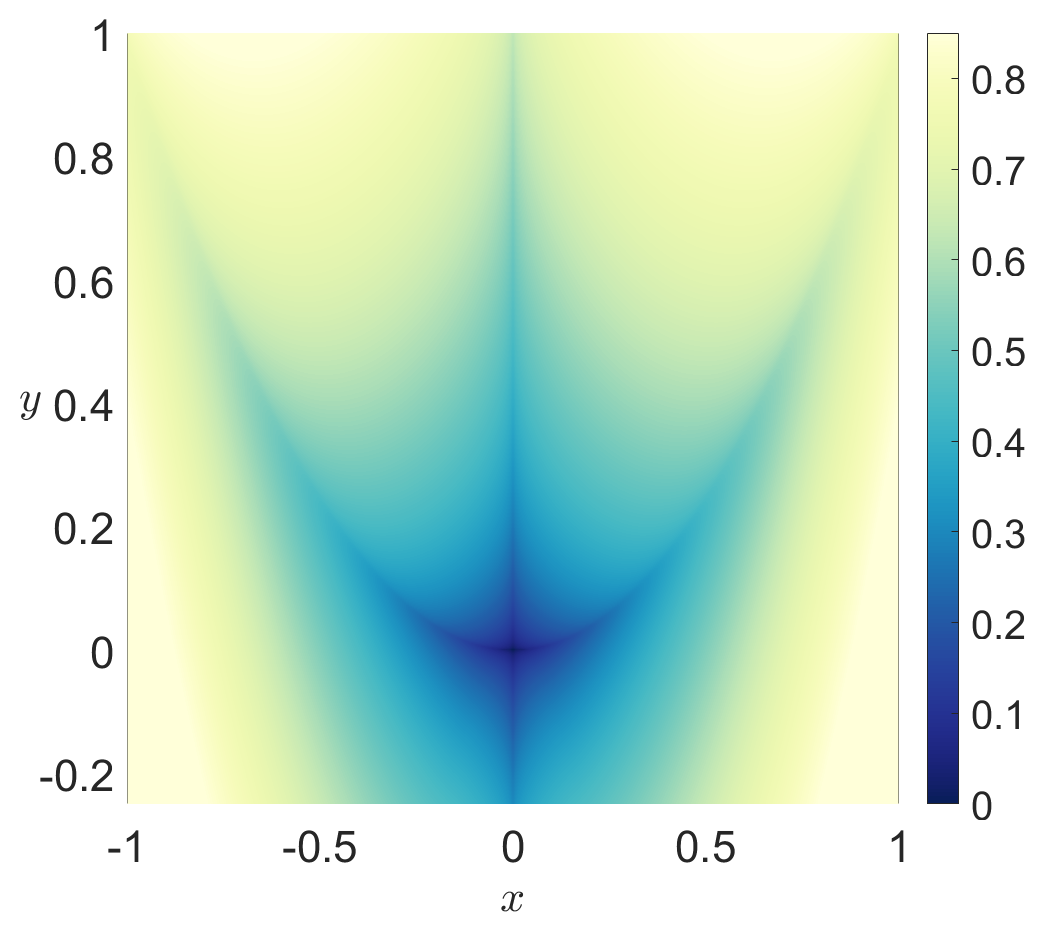} 
		C)\includegraphics[scale=0.19]{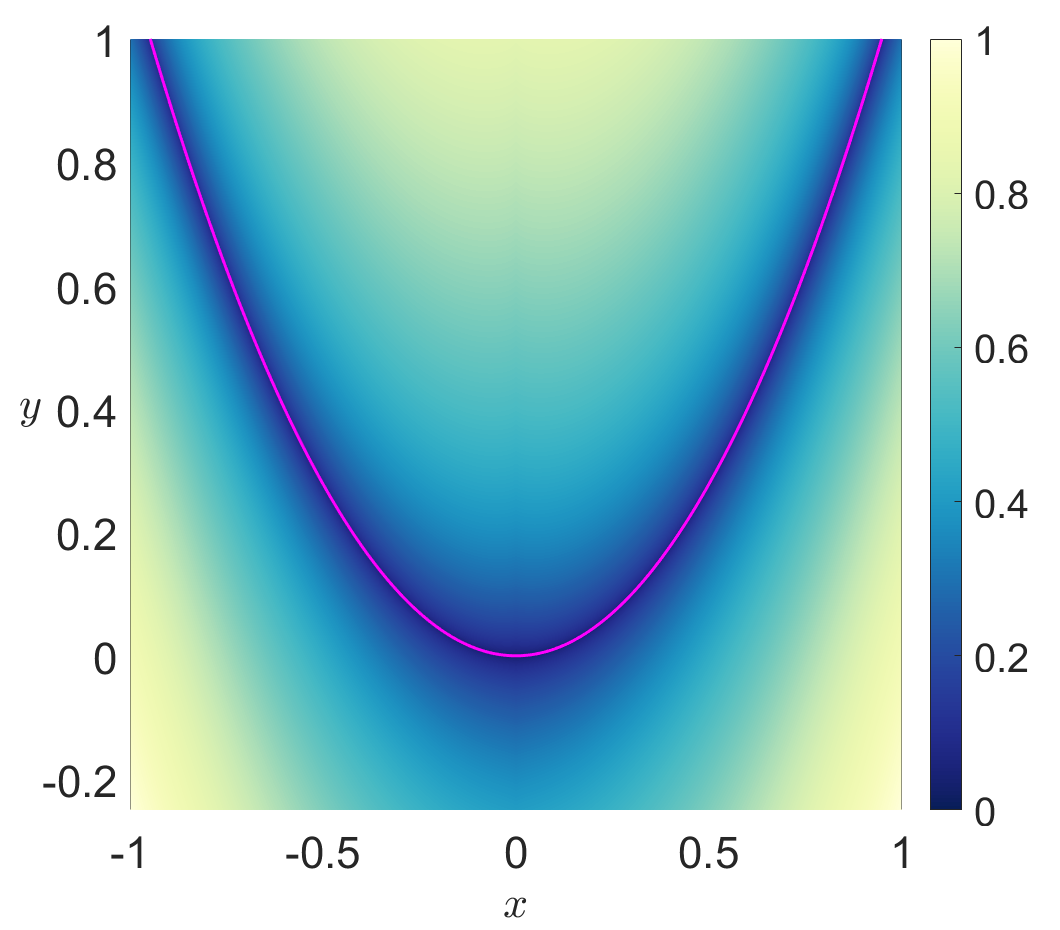}  
	\end{center}
	\caption{Phase space of the system in Eq. \eqref{nonlin_saddle} with model parameters $\lambda = -1$ and $\mu = -0.05$, as revealed by applying the $p$-norm definition of LDs with $p = 1/2$ using an integration time of $\tau_f = \tau_b = 5$. A) Total LDs. B) Forward LDs. C) Backward LDs. The magenta curve represents the slow manifold given by Eq. \eqref{slow_mani}.}
	\label{ld_slowmani}
\end{figure}

In this subsection we look at the capability of LDs to reveal slow manifolds in systems that display a separation of time scales \cite{kuehn2015}. Slow manifolds are relevant in many applications, such as those that involve the dynamical study of gliding animals \cite{nave2019}. We focus our attention first on the dynamical system we introduced in Eq. \eqref{nonlin_saddle} which has been widely studied in the literature as a basic model for slow-manifold dynamics \cite{brunton2016,lusch2018}. It is well known that for $\lambda \ll \mu < 0$ this system has a slow manifold at the curve:
\begin{equation}
    y = \dfrac{\lambda}{\lambda - 2\mu} x^2 \;.
    \label{slow_mani}
\end{equation}
We will study here the case where $\mu = -0.05$ and $\lambda = -1$. For this situation, the dynamics in the $y$-coordinate rapidly evolves towards making $\dot{y}$ very small, and at that point, trajectories follow the asymptotically attracting slow manifold so that dynamics on $x$ become dominant. To reveal the location of the slow manifold we apply LDs and integrate trajectories for $\tau_f = \tau_b = 5$. In Fig. \ref{ld_slowmani} we display from left to right the total, forward and backward LD scalar field respectively, and for comparison we overlay in magenta the location of the slow manifold given by Eq. \eqref{slow_mani}. This test illustrates that the method is correctly highlighting the slow manifold location in the phase space of the system.

\begin{figure}[htbp]
	\begin{center}
	A)\includegraphics[scale=0.28]{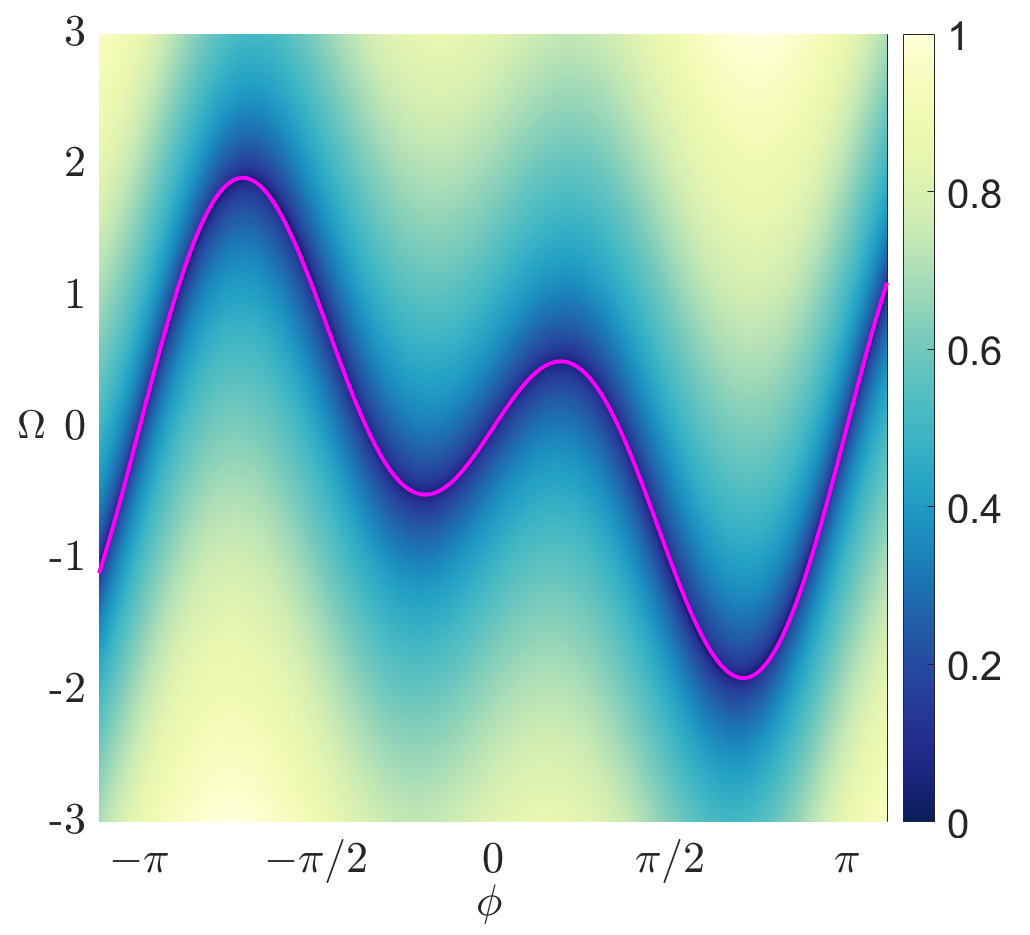} 
	B)\includegraphics[scale=0.28]{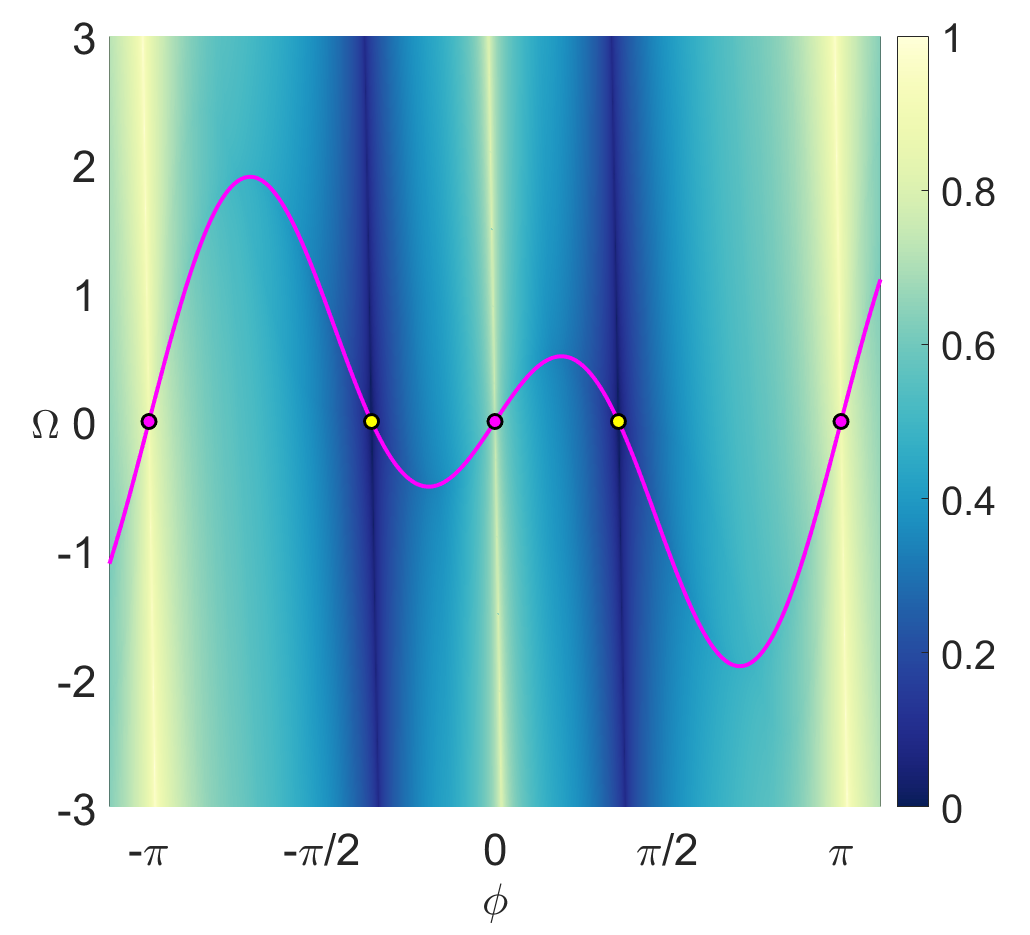} 
	C)\includegraphics[scale=0.28]{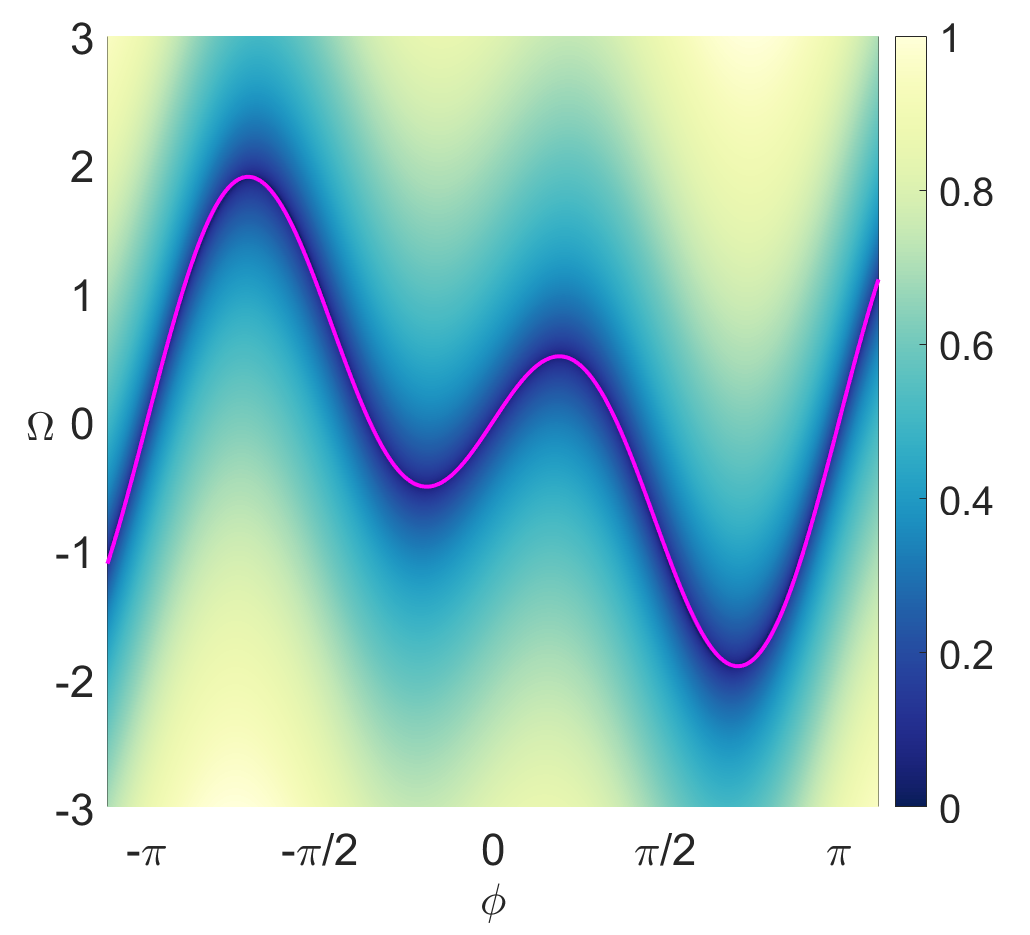}  
	D)\includegraphics[scale=0.28]{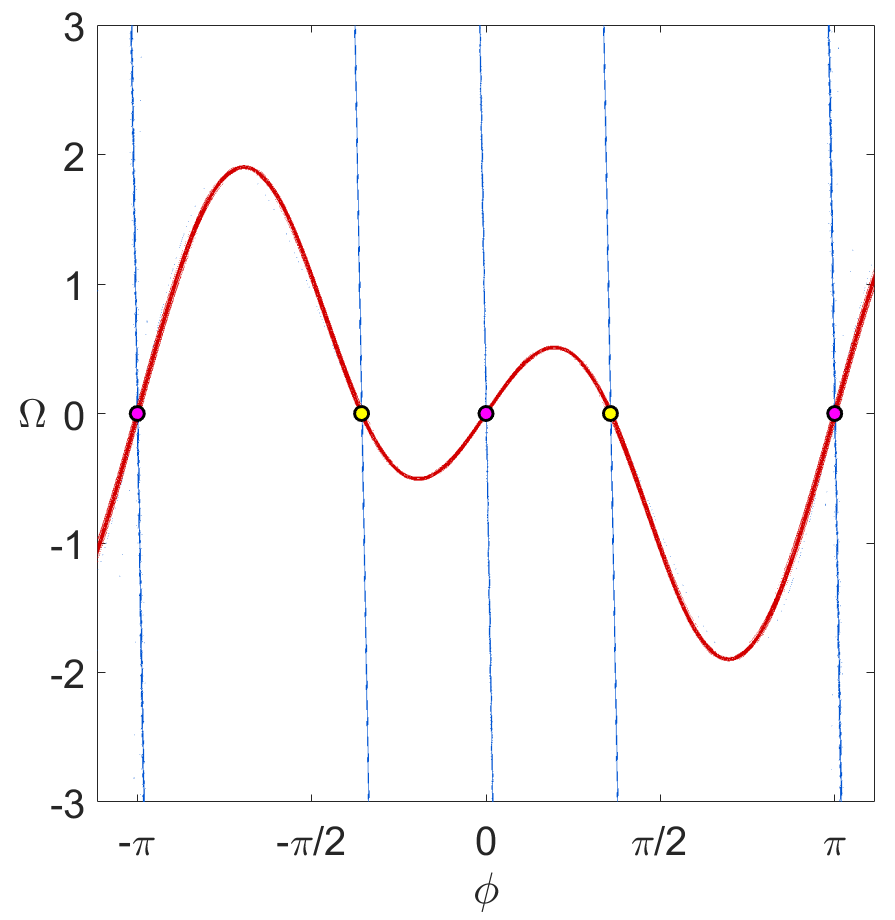}  
	\end{center}
	\caption{Phase space reconstruction of the bead on a rotating hoop problem in Eq. \eqref{bead} for the system parameters $\varepsilon = 0.02$ and $\mu = 2.3$, by means of LDs with $p = 1/2$. Trajectories are integrated forward and backward for a time $\tau_f = \tau_b = 10$. A) Total LD fuction; B) Forward component of LDs; C) Backward component of LDs; D) Stable (blue) and unstable (red) manifold extraction using the laplacian operator on the LD scalar output, that is, $\Delta \mathcal{L}_p$. In panels A)-C) we have depicted the slow manifold of the system in magenta. Yellow dots correspond to stable equilibria, and magenta dots represent saddle points.}
	\label{ld_bead}
\end{figure}

In order to provide further evidence of the capability of this tool for revealing slow manifolds, we focus next on two classical model dynamical systems that have been widely studied in the literature. First we analyze the classical mechanics problem of a bead in a rotating hoop. A detailed analysis of the stability of this system can be found in \cite{strogatz}. We are interested in solving the dynamical system which is given by:
\begin{equation}
	\begin{cases}
		\dot{\phi} = \Omega \\[.2cm]
		\dot{\Omega} = \dfrac{1}{\varepsilon} \left(f(\phi) - \Omega\right)
	\end{cases}
	\label{bead}
\end{equation}
where $f(\phi) = \left(\mu \cos(\phi) - 1\right) \sin(\phi)$. We will look at the case where $\mu > 1$ for which the system has two unstable equilibrium points at $\phi = 0$ and $\phi  = \pi$, and two stable equilibria at $\phi = \arccos\left(1/\mu\right)$. we will study the slow manifold that exists in the system for $\varepsilon = 0.02$ and $\mu = 2.3$. We compute LDs by integrating trajectories forward and backward for the same time $\tau_f = \tau_b = 10$. The output of these numerical experiments is presented in Fig. \ref{ld_bead}. As before, we have superimposed in the picture the slow manifold of the system in magenta, which is given by the curve $\Omega = f(\phi)$, in order to validate the results obtained. We have also depicted with yellow dots the stable equilibria, and magenta dots represent unstable equilibria. Notice that the method also reproduces for this problem the slow manifold, and it is also highlighting the stable and unstable manifolds of the saddle equilibrium points of the system. We check this fact by extracting the geometrical features of the phase space from the LD scalar field by means of applying the laplacian operator to it. Observe that with this approach we can also obtain an accurate approximation of the slow manifold from the ridges of the laplacian field $\Delta \mathcal{L}_p$.

\begin{figure}[htbp]
	\begin{center}
	A)\includegraphics[scale=0.28]{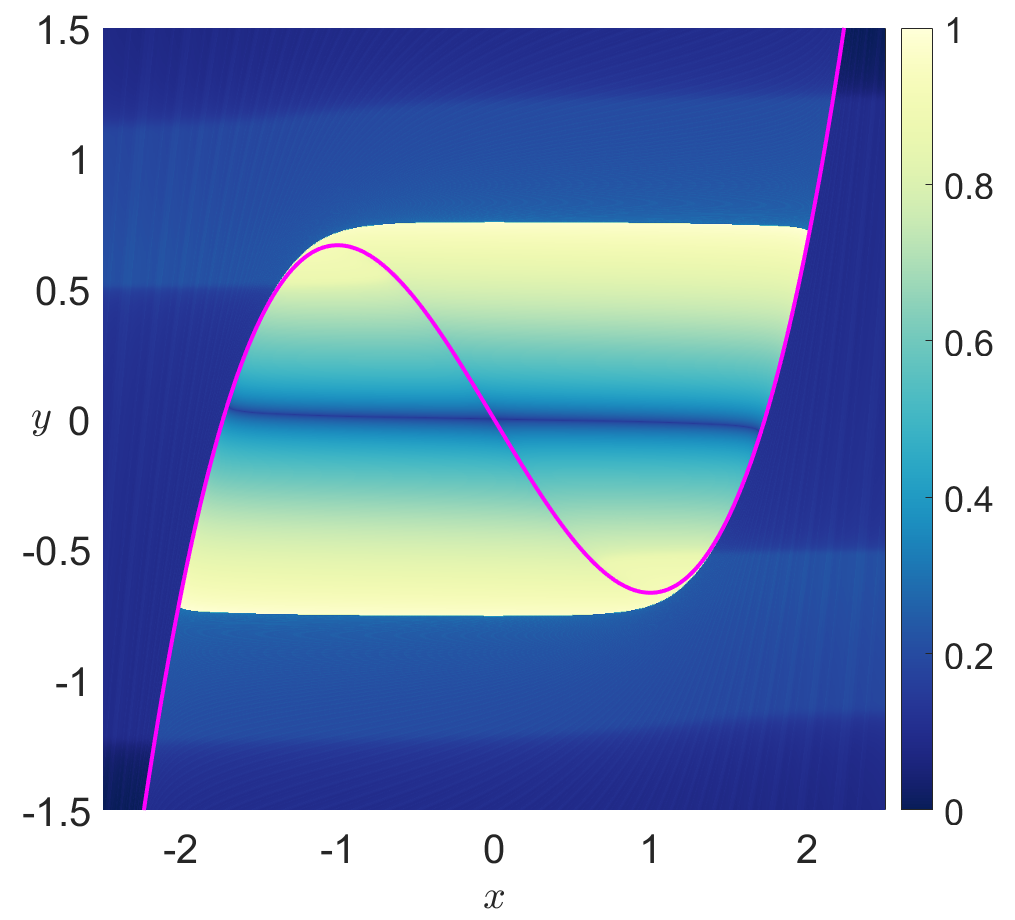} 
	B)\includegraphics[scale=0.28]{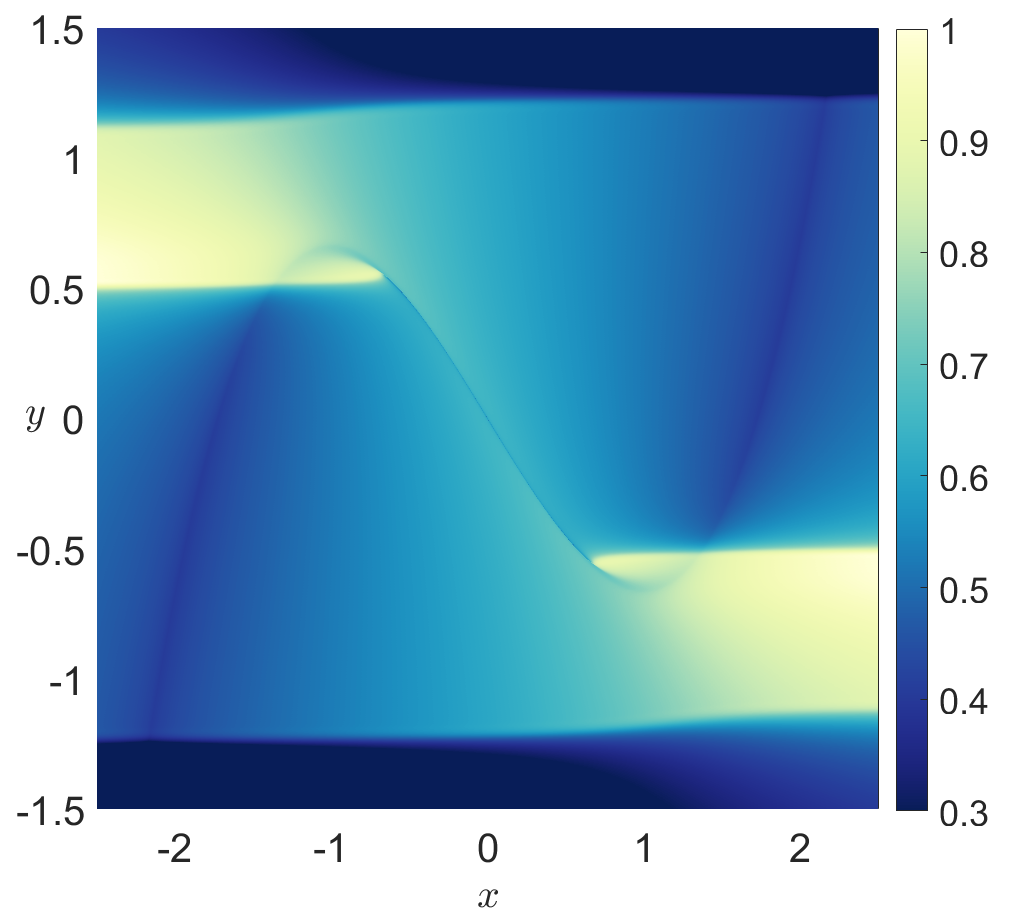} 
	C)\includegraphics[scale=0.28]{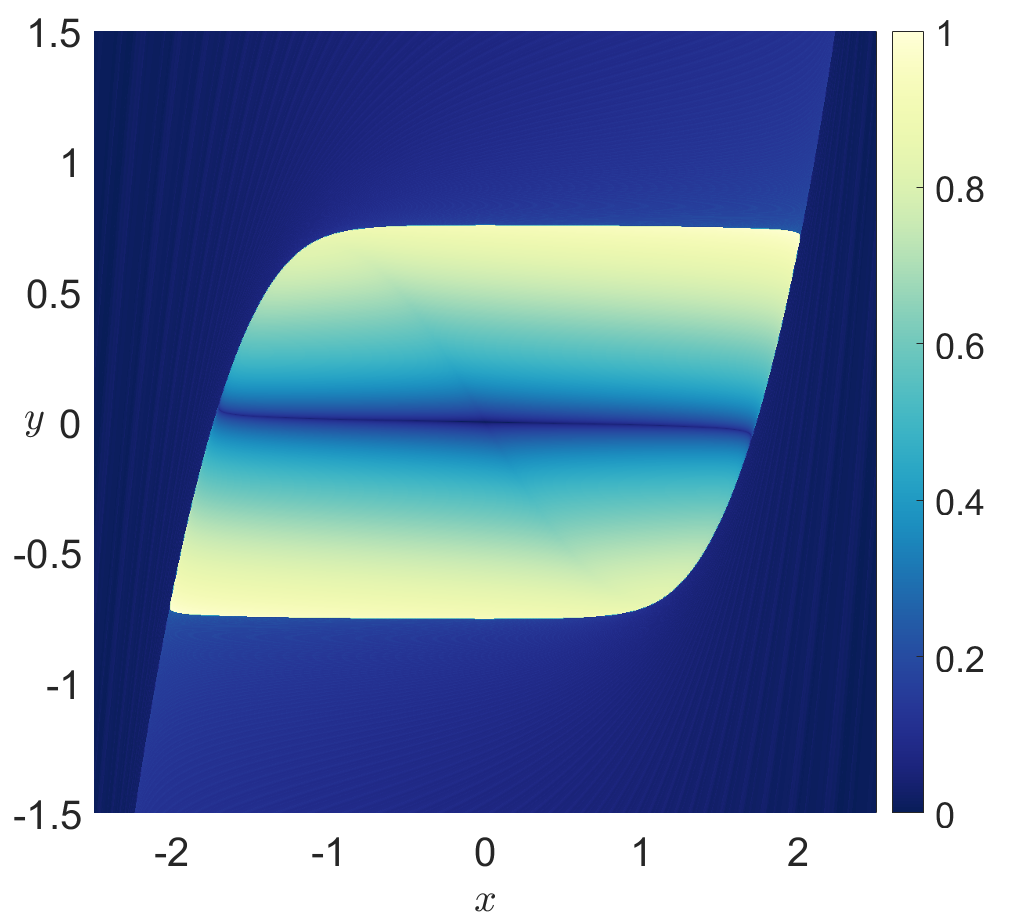}
	D)\includegraphics[scale=0.28]{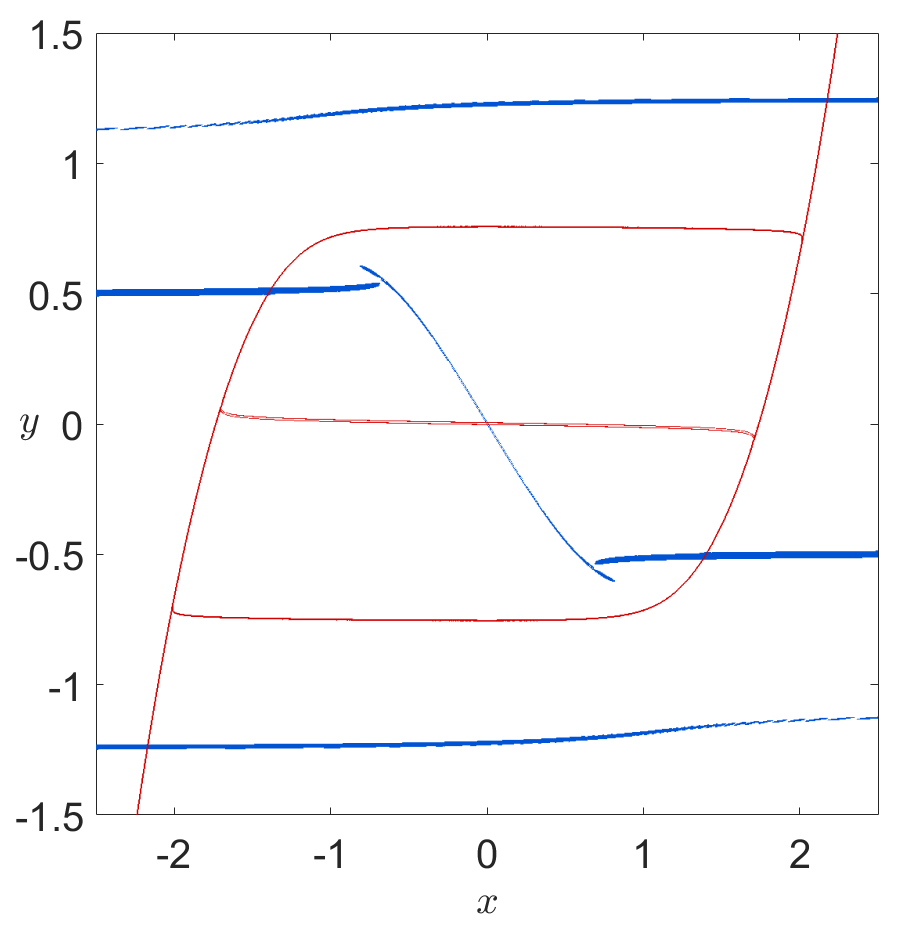}
	\end{center}
	\caption{Phase space visualization of the van der Pol oscillator under a Li\'enard transformation given by Eq. \eqref{lienard} for the system parameters $\mu = 10$, by means of LDs with $p = 1/2$. Trajectories are integrated forward and backward for a time $\tau_f = \tau_b = 50$. A) Total LD function; B) Forward component of LDs; C) Backward component of LDs; D) Approximation to the stable (blue) and unstable (red) branches of the slow manifold displayed by the system using the gradient operator on the LD scalar output. In panel A) we have depicted the slow manifold of the system in magenta.}
	\label{ld_vdp_lienard}
\end{figure}

We finish our analysis of the detection of slow manifolds with Lagrangian descriptors by revisiting the van der Pol oscillator problem in Eq. \eqref{vdp_sys} where $\mu \gg 1$. In this case, in order to visualize the phase space of the system we will carry out a Li\'{e}nard transformation. Consider the following definitions:
\begin{equation}
	F(x) = \dfrac{1}{3}x^3 - x \quad,\quad w = \dfrac{1}{\mu}y + F(x) = \dfrac{1}{\mu}\dot{x} + F(x) 
\end{equation}
then we can write Eq. \eqref{vdp_sys} in the new variables $x$ and $w$ as:
\begin{equation}
		\begin{cases}
		\dot{x} = \mu \left(w - F(x)\right) \\[.2cm]
		\dot{w} = -\dfrac{1}{\mu} x
	\end{cases}
	\label{lienard}
\end{equation}
This system has a slow manifold at the curve $w = F(x) = x^3/3 - x$, and in order to test how LDs perform when it comes to detecting this structure, we run a simulation using $\mu = 10$ and we integrate trajectories for $\tau_f = \tau_b = 50$ or until they leave a circle of radius $R = 6$ centered at the origin. Results are shown in Fig. \ref{ld_vdp_lienard}, and we have overlaid in magenta the curve that represents the true slow manifold of the system. This comparison between the output of LDs and the location of the slow manifold underlines the success of the method in highlighting this object. Furthermore, we can extract this information from the gradient of LDs, and obtain an approximation of the slow manifold from the scalar field values of the LD output. Notice that this analysis is also revealing the slow and fast branches of the underlying limit cycle of the system.

\subsection{Identification of Attractors}
\label{subsec:sec4}

\begin{figure}[htbp]
	\begin{center}
	A)\includegraphics[scale=0.28]{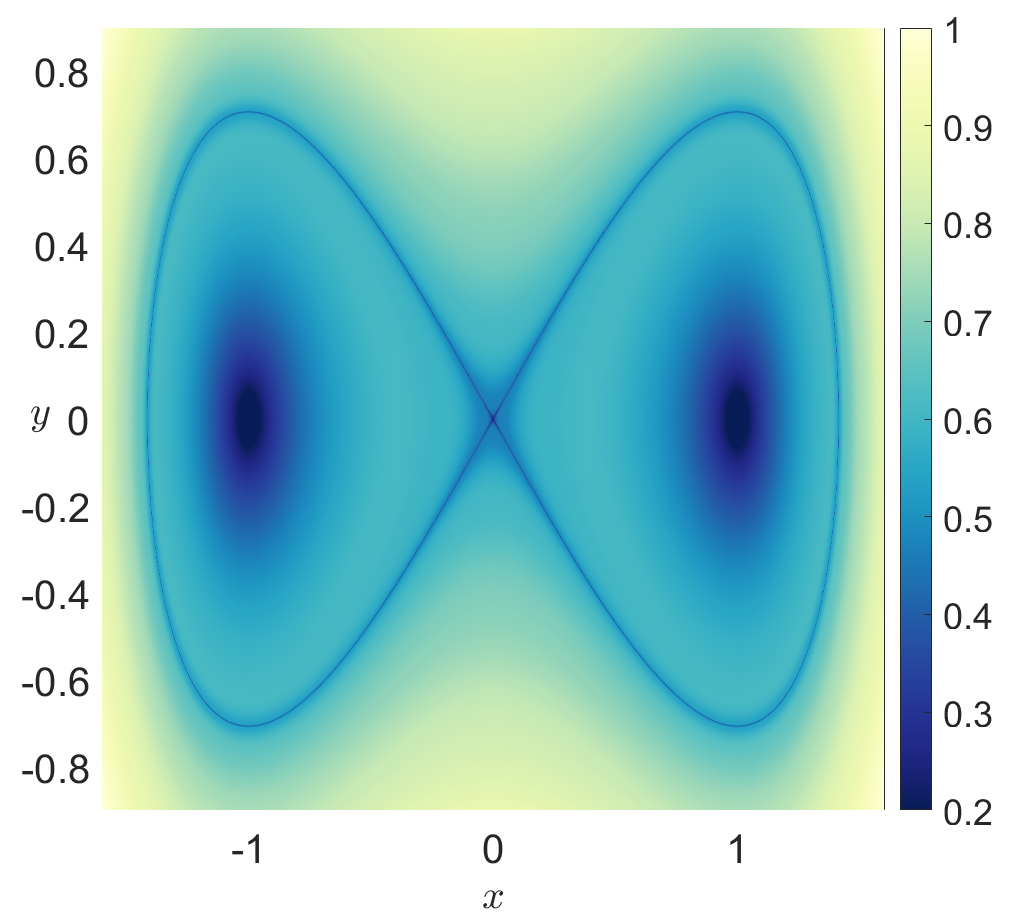} 
	B)\includegraphics[scale=0.28]{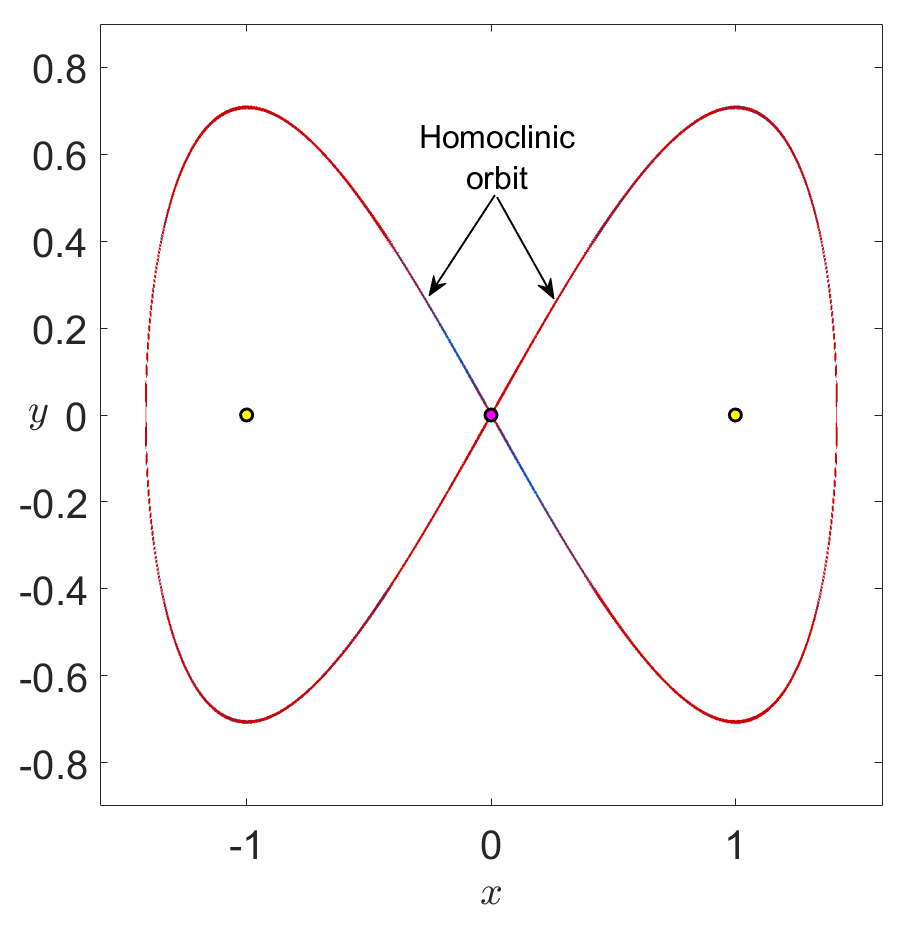}
	\end{center}
	\caption{Phase space of the Duffing system in Eq. \eqref{duffing} with no forcing and no damping, that is $\gamma = \delta = 0$, highlighted by LDs with $p = 1/2$ using a forward and backward integration tome of $\tau = 20$. A) Total LD values; B) Stable (blue) and unstable (red) manifolds extracted from the gradient of the LD scalar field. We have marked with yellow dots the stable equilibria, and with a magenta dot the saddle point at the origin.}
	\label{duffi_cons}
\end{figure}

The next dynamical topic that we would like to address is the capability of the method of Lagrangian descriptors to detect attractors and their basins of attraction. To illustrate this fact, we have chosen the Duffing equation as the model problem to study \cite{duffing1918,guck1983,korsch2008,kovacic2011}. The dynamics of the periodically forced and damped Duffing oscillator is described by Newton's second law:
\begin{equation}
	\ddot{x} + \delta \dot{x} - \alpha x + \beta x^3 = \gamma \cos(\omega t)
\end{equation}
which can be rewritten as a system of first order ODEs in the form:
\begin{equation}
	\begin{cases}
		\dot{x} = y \\[.2cm]
		\dot{y} = - \delta y + \alpha x - \beta x^3 + \gamma \cos(\omega t)
	\end{cases}
	\label{duffing}
\end{equation}
We will consider in our simulations that $\beta = 1$. We begin by illustrating how in the conservative unforced case, that is, $\gamma = \delta = 0$ we get the classical phase space picture of a double well system with a saddle point at the origin and two homoclinic orbits surrounding the center equilibria. To do so, we calculate LDs by using equal integration times forward and backward $\tau_f = \tau_b = 20$ and depict the results in Fig. \ref{duffi_cons}. In panel A) we present the total LD field, and in B) we have extracted the stable and unstable manifolds as ridges of the LD output by means of the gradient. The location of the centers is marked with yellow dots, and the magenta dot at the origin corresponds to the saddle point. Notice how LDs nicely highlight the figure eight shape formed by the homoclinic orbits that enclose the center equilibria of the system.

\begin{figure}[htbp]
	\begin{center}
	A)\includegraphics[scale=0.28]{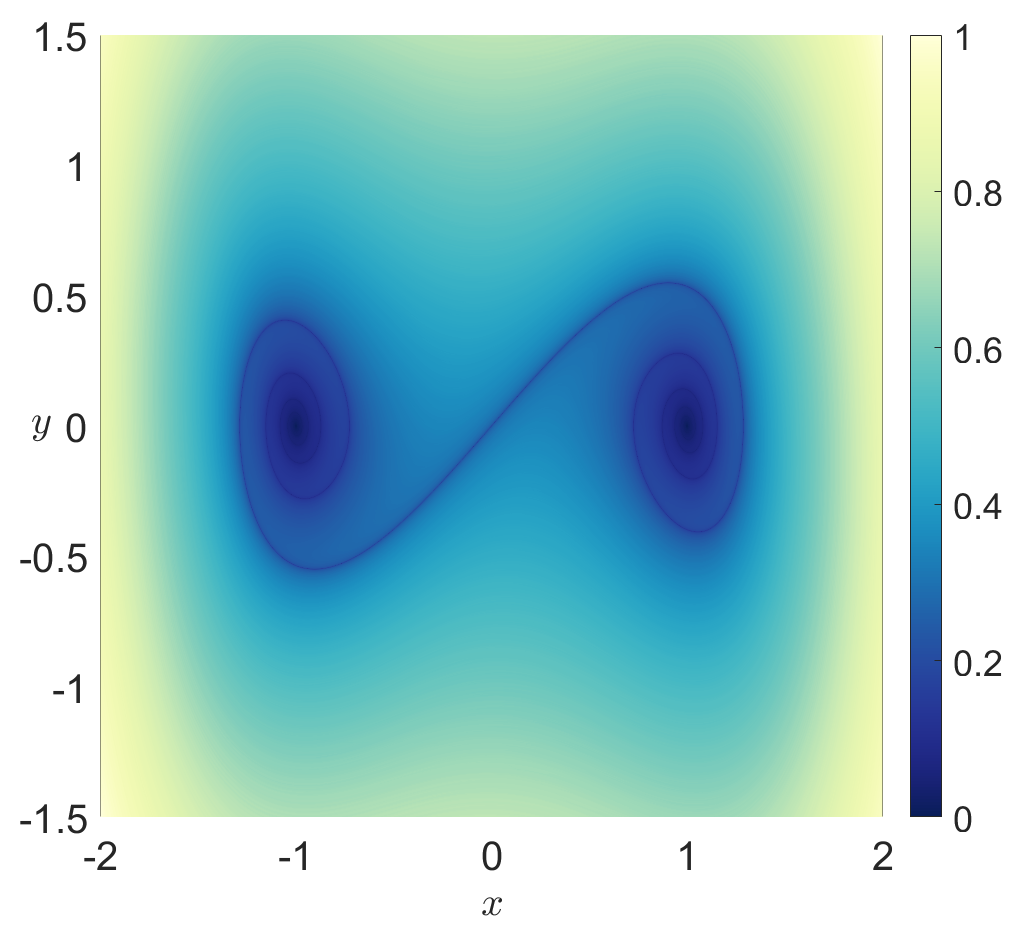} 
	B)\includegraphics[scale=0.28]{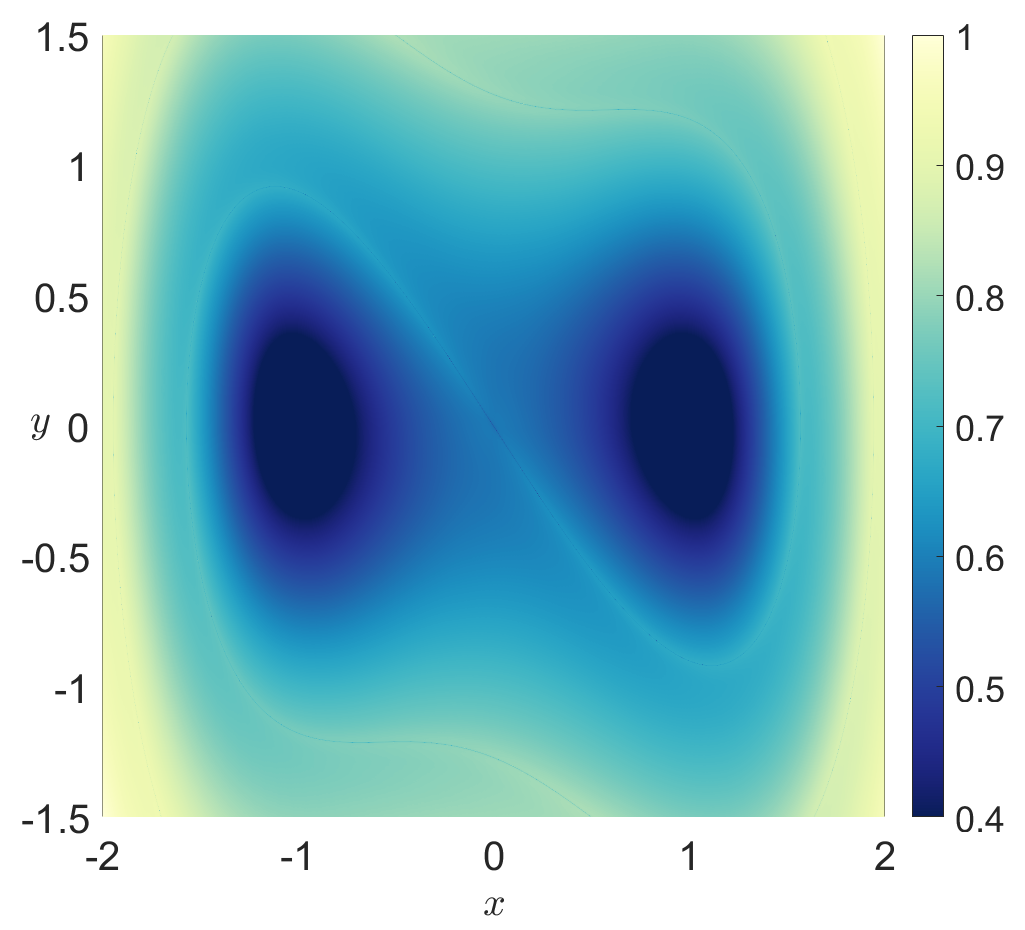} 
	C)\includegraphics[scale=0.28]{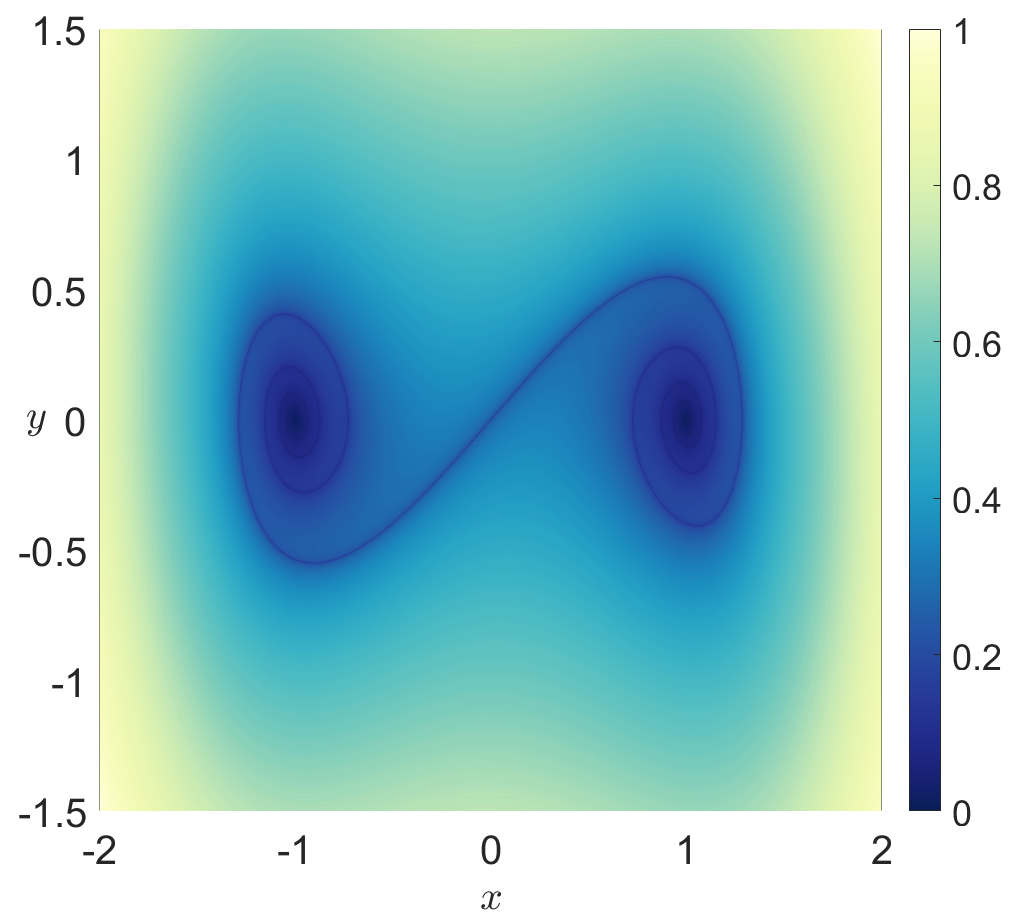} 
	B)\includegraphics[scale=0.28]{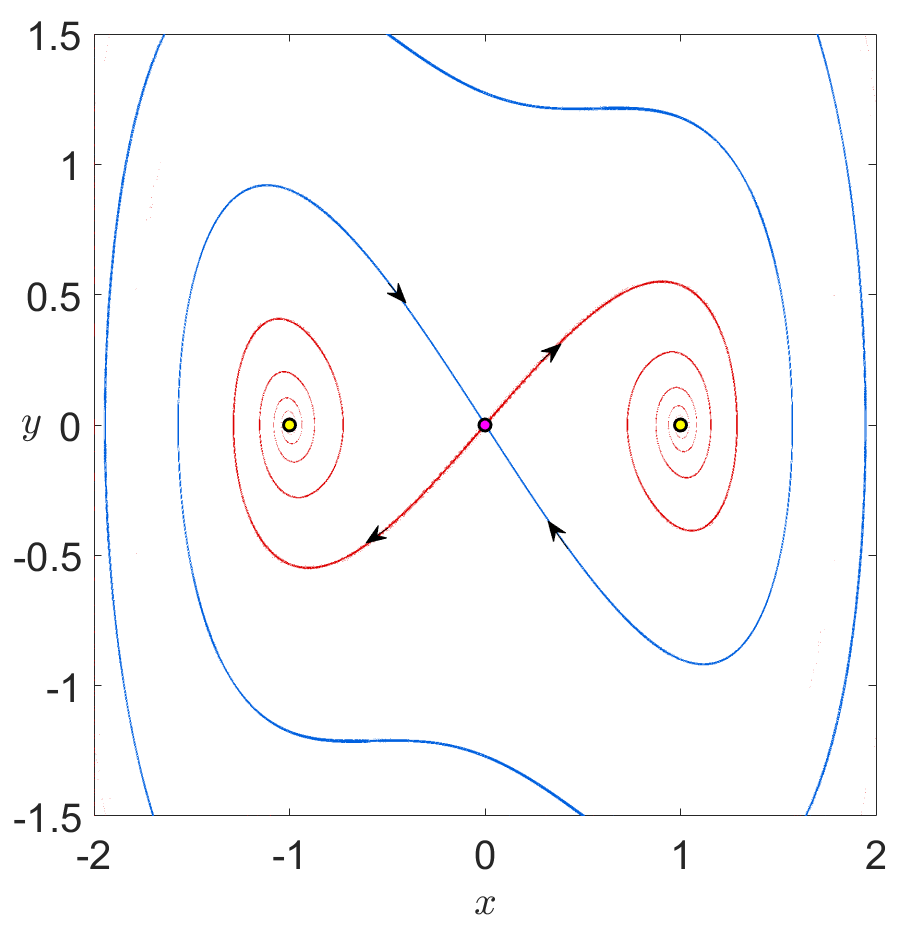} 
	\end{center}
	\caption{Phase space of the Duffing system in Eq. \eqref{duffing} with no forcing ($\gamma = 0$), and a damping coefficient of $\delta = 0.3$, as displayed by LDs with $p = 1/2$ using a forward and backward integration tome of $\tau = 25$. A) Total LD values; B) Forward LD; C) Backward LD; D) Stable (blue) and unstable (red) manifolds extracted from the gradient of the LD scalar field. We have marked with yellow dots the stable equilibria, and with a magenta dot the saddle point at the origin. For guidance, we have included arrows to mark the flow directions on the manifolds.}
	\label{duffi_diss}
\end{figure}

In the next experiment we continue working with the unforced Duffing oscillator ($\gamma = 0$) but switch on the damping in the system in order to show that LDs are capable of detecting the boundaries of the basins of attraction of both centers, that now turn into stable foci. We set the damping coefficient to $\delta = 0.3$, and for this parameter value, the eigenvalues at the origin are $\lambda_1 \approx 0.8612$ and $\lambda_2 \approx -1.1612$. Notice that this would imply that we have two different time scales at play in the neighborhood of the origin, and this could difficult the simultaneous visualization of the manifolds when we plot the scalar output of the total LD (forward plus backward), as we discussed for the example of the dissipative saddle. This becomes evident in Fig. \ref{duffi_diss} A) when we compute the total LDs (forward plus backward contributions) for $\tau_f = \tau_b = 25$. Panels B) and C) show the forward and backward components of LD respectively, and in them we can see that the method reveals the location of the stable and unstable manifolds of the system. These manifolds form the boundaries of the basins of attraction of the stable foci depicted with yellow dots in panel D). Despite the fact that it is difficult to visualize both manifold structures in the total LD field because of the different time scales at play, we can easily extract their locations by analyzing the gradient of the forward and backward components of LD separately and plot these information together as we illustrate in D). 

\begin{figure}[htbp]
	\begin{center}
	A)\includegraphics[scale=0.26]{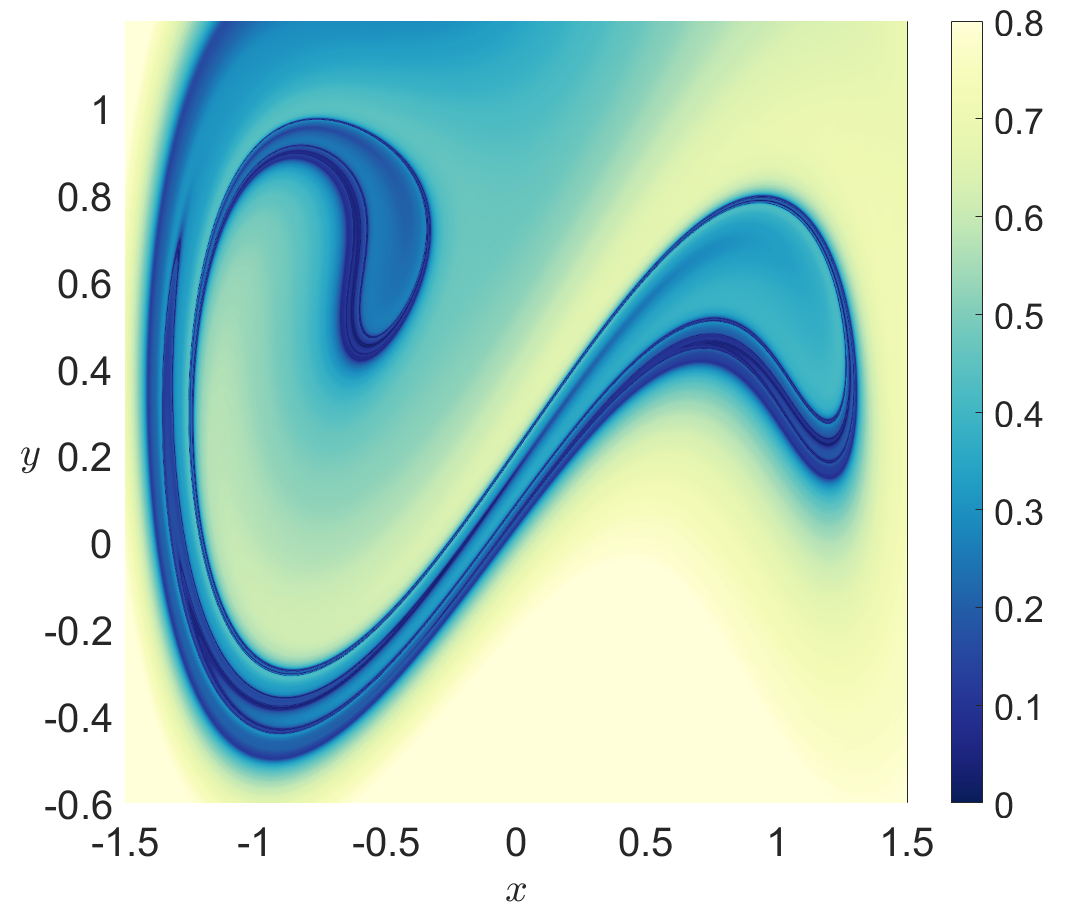} 
	B)\includegraphics[scale=0.26]{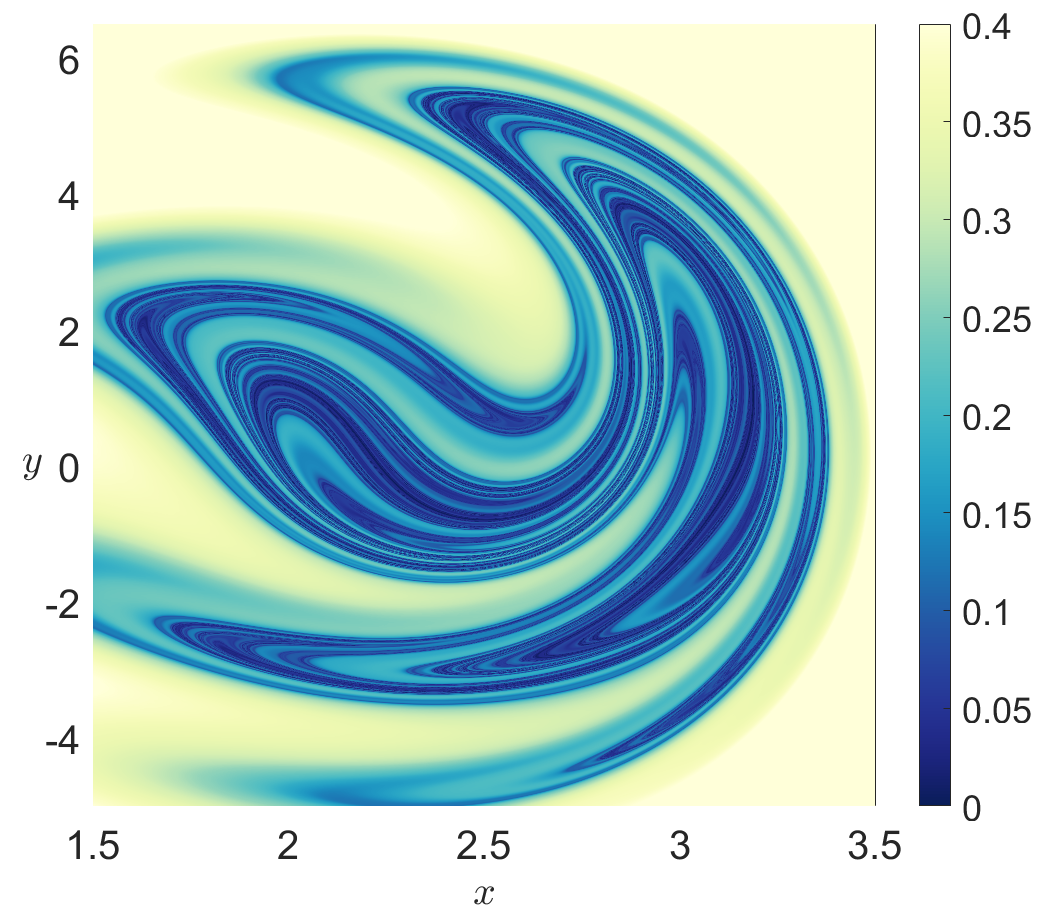} 
	C)\includegraphics[scale=0.28]{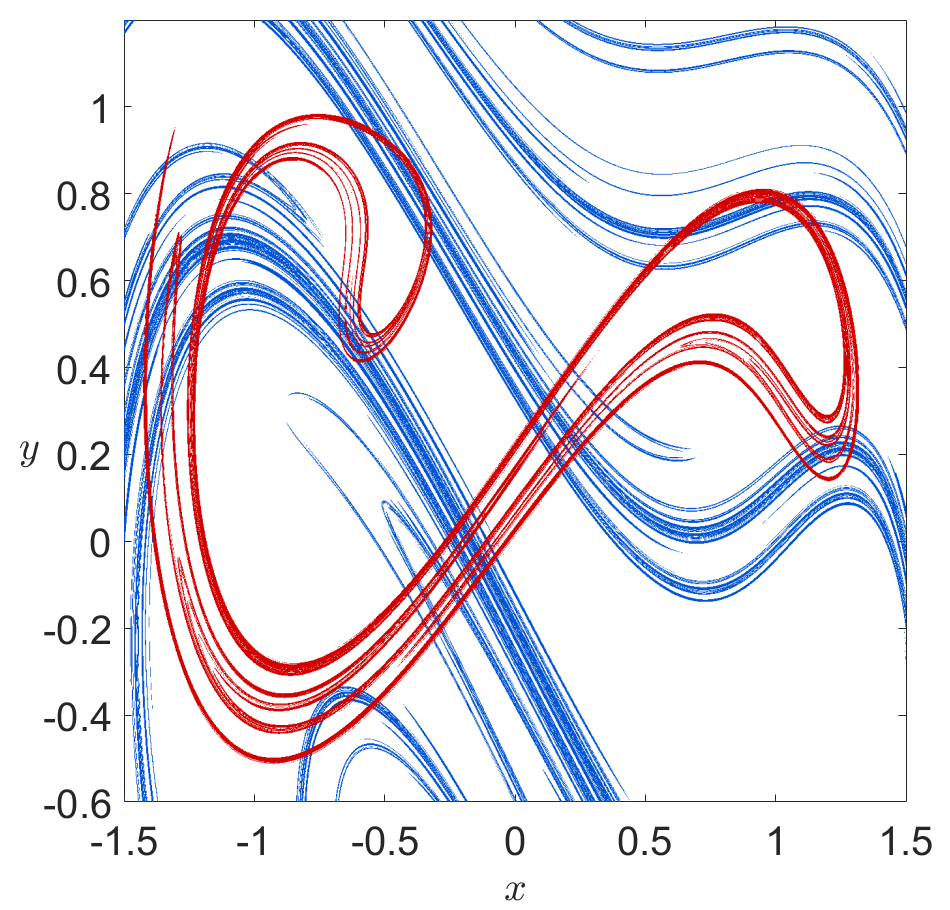} 
	D)\includegraphics[scale=0.28]{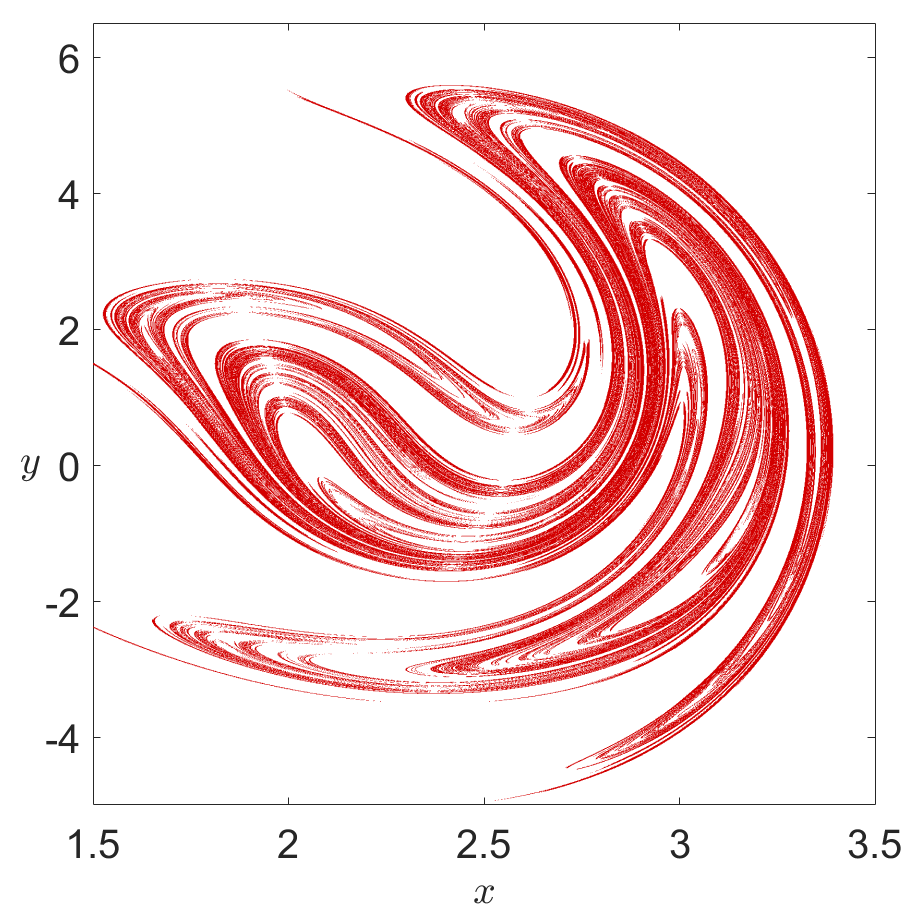} 
	E)\includegraphics[scale=0.28]{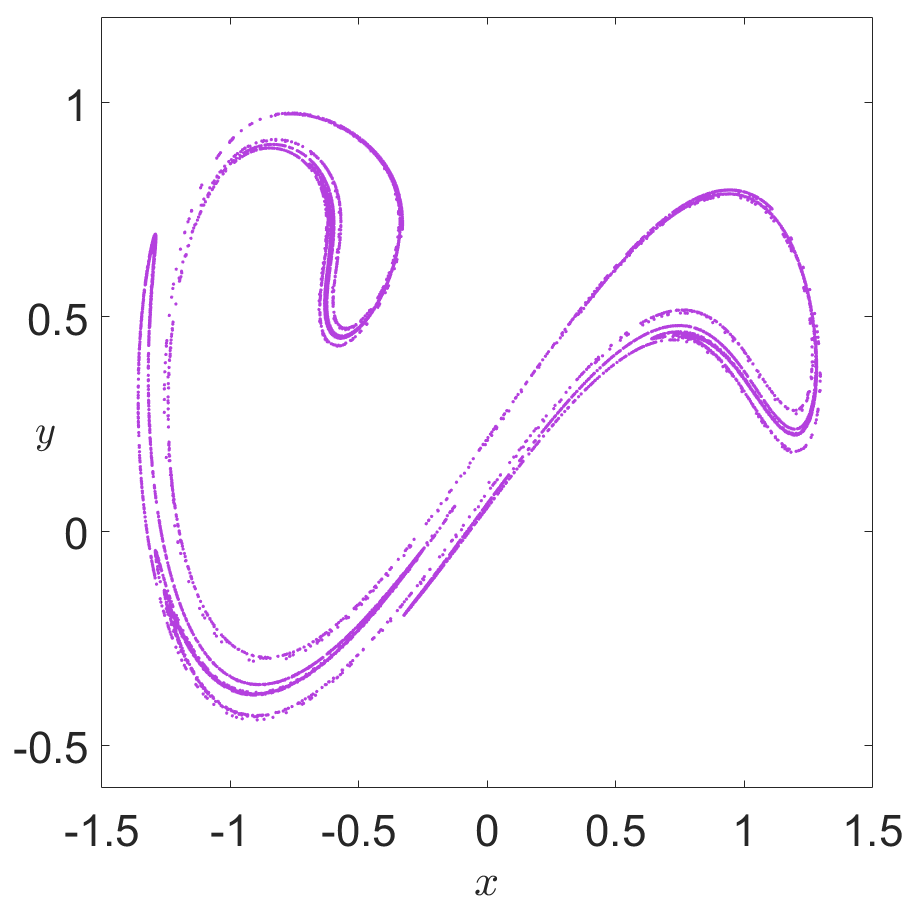} 
	F)\includegraphics[scale=0.28]{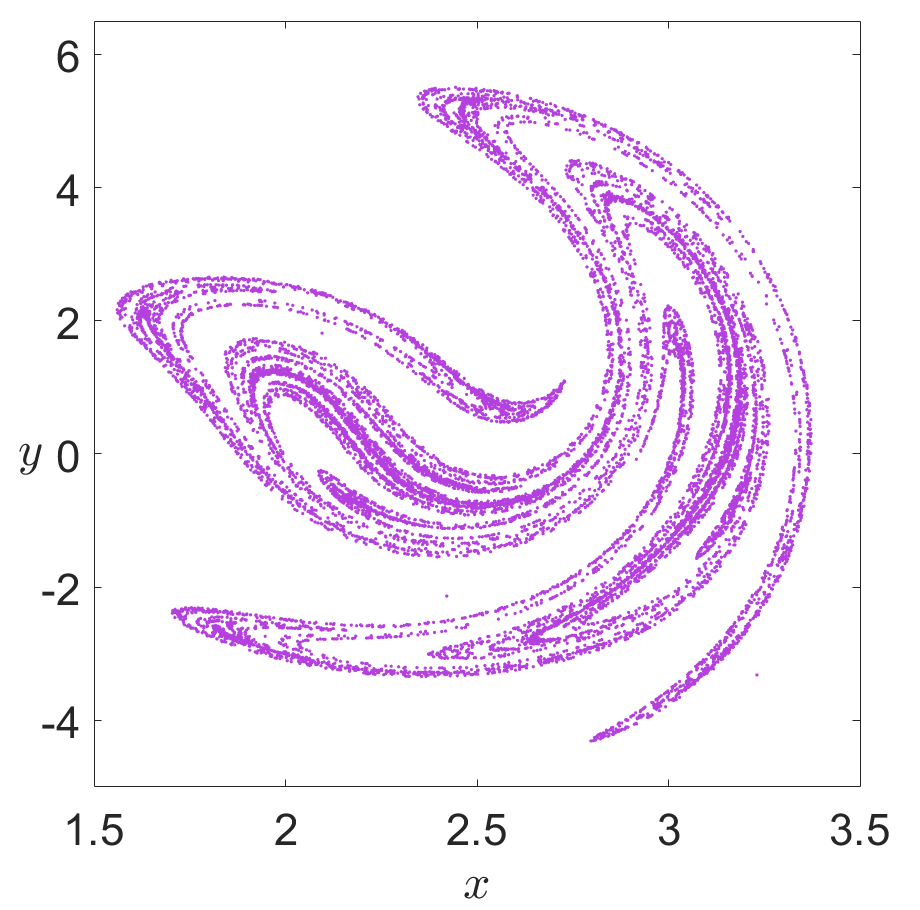} 
	\end{center}
	\caption{Strange attractors for the Duffing system in Eq. \eqref{duffing} at time $t = 0$ as revealed by LDs, and comparison with Poincar\'e sections. The first column corresponds to the model parameters $\alpha = 1$, $\delta = 0.3$, $\gamma = 0.5$ and $\omega = 1.2$, whereas the second column uses the values $\alpha = 0$, $\delta = 0.05$, $\gamma = 7.5$ and $\omega = 1$ (Ueda attractor). A) and B) depict the backward coponent of LDs with $p = 1/2$ using an integration time of $\tau = 20$ and $\tau = 50$ respectively. In panel C) we show the stable (blue) and unstable (red) manifolds extracted by means of the gradient of the LD scalar field. D) displays the unstable manifold extracted from the laplacian operator. E) and F) illustrate the Poincar\'e section obtained by strobing the location of a trajectory, starting from the initial condition $(1,0)$, every period of the forcing term for $15000$ periods.}
	\label{str_att}
\end{figure}

We continue the analysis of attractors with LDs by demonstrating that this tool can be used to reconstruct the intricate geometry of strange attractors. To do so, we explore two different cases for the model parameters and introduce forcing into the system, making it nonautonomous. We look at two different situations, one where $\alpha = 1$, $\delta = 0.3$, $\gamma = 0.5$ and $\omega = 1.2$ (displayed in the left column of Fig. \ref{str_att}), and for the other we select $\alpha = 0$, $\delta = 0.05$, $\gamma = 7.5$ and $\omega = 1$ (right column of Fig. \ref{str_att}). The second case  corresponds to the setting for the Ueda attractor \cite{ueda1979,ueda1980}. The existence and properties of these strange attractors for the Duffing equation have been widely studied in the literature, see e.g.  \cite{moon1979,holmes1980,guck1983,kovacic2011} and references therein. In order to reveal the strange attractor and its stable and unstable manifolds at time $t = 0$ we compute LDs with an integration time $\tau_f = \tau_b = 20$. Since the system is nonautonomous, the attractor is going to evolve and change shape with time, so our analysis is going to provide a snapshot of it at the initial time $t = 0$. If we want to obtain a picture of the attractor at any other time, say $t = t_1$, one just needs to compute LDs using that instant $t = t_1$ as the initial time. We compare our results with the classical technique used to reconstruct strange attractors that consists in calculating a Poincar\'e section by strobing the trajectory of an initial condition at integer multiples of the forcing period. For this test we have chosen $(1,0)$ as the initial condition at time $t = 0$, and we have recorded its location every period of the forcing term for $15000$ periods. In Fig. \ref{str_att} A) and B) we depict the forward component of LDs, while in the panels below, C) and D), the geometry of the stable (blue) and unstable (red) manifolds has been recovered from the LD values with a gradient-type filter. Finally, and as a validation of the results obtained with LDs, we present in E) and F) the corresponding Poincar\'e map in purple.

It is important to remark here that the reconstruction of the attractor by means of a Poincar\'e section approach requires a extremely large integration time, which would be prone to the accumulation of numerical error in the simulation. On the other hand, with the two examples we have discussed above, LDs provide a clear advantage for recovering the structure and location of the attractor, since our experiments were carried out with a small integration time of $\tau_f = \tau_b = 20$, which roughly corresponds to a couple of periods of the forcing, a number that by no means is comparable to the large time interval needed for the Poincar\'e map, Moreover, Poincar\'e maps can be only constructed for periodically forced systems, whereas LDs can be applied in a straightforward manner to both aperiodic and periodic forcing scenarios.

\subsection{Transition Ellipsoids in Hamiltonian systems with two degrees of freedom subject to dissipation}
\label{subsec:sec5}

We would like to finish this work by exploring the capability that the method of Lagrangian descriptors brings for the dynamical analysis of Hamiltonian systems with two degrees of freedom (DoF) subject to damping effects. In many applications, such as astrodynamics \cite{jaffe2002,koon2011}, structural mechanics \cite{collins2012,zhong2018snap}, capsize ship motion \cite{naik2017}, chemical reactions \cite{Uzer2002}, etc., it is of crucial importance to address transition events across an index-1 saddle critical point of a potential energy surface (PES). The phase space objects responsible for controlling this transport mechanism are the stable and unstable manifolds of the normally hyperbolic invariant manifold (NHIM) that exists in the vicinity of the index-1 saddle \cite{Wiggins94}. For 2 DoF Hamiltonian systems, the NHIM is an unstable periodic orbit, that is, it has the geometry of $S^1$, and its stable and unstable manifolds have the topology of cylinders ($S^1\times \mathbb{R}$). These 'tube manifolds' are known in the literature as spherical cylinders \cite{almeida1990,deleon1991} and act as conduits, i.e. a transportation network, that trajectories moving across the index-1 saddle follow along their evolution. Moreover, they are codimension-1 objects in the energy hypersurface of the Hamiltonian system, and thus, they are barriers to transport that trajectories can not cross \cite{wiggins2001}. These cylinders partition the energy surface into 'reactive' and 'non-reactive' trajectories. All trajectories starting from an initial condition inside the tube manifolds are termed as 'reactive', which means that they will cross the index-1 saddle region of the PES, whereas the initial conditions outside the tubes will not move across the index-1 saddle. 
However, when damping is added to the Hamiltonian system, it has been recently shown that the topology of the manifolds of the NHIM changes significantly and transforms from a cylindrical geometry into an ellipsoid \cite{zhong2020a,zhong2020b,zhong_thesis}. In particular, the stable cylinder of a conservative Hamiltonian system, which characterizes the trajectories crossing from one side to the other of the index-1 saddle in forward time, becomes an ellipsoid, known as a transition ellipsoid. This phase space structure determines all transition events in forward time in the dissipative setup. Furthermore, the unstable periodic orbit associated to the index-1 saddle equilibrium point that existed in the conservative system disappears and converts into a focus-type asymptotic orbit tending to the equilibrium point itself.

In order to demonstrate that LDs can be used as a simple mathematical tool to detect the location and geometry of this ellipsoid, we will analyze an uncoupled Hamiltonian system with 2 DoF defined from a double well potential in the $x$ DoF, and a harmonic oscillator in $y$. To this model, we will add a linear damping into Hamilton's equations of motion. Consider the Hamiltonian function:
\begin{equation}
H(x,y,p_x,p_y) = \dfrac{p_x^2}{2m_1} + \dfrac{p_y^2}{2m_2} + \dfrac{a}{4}x^4 - \dfrac{b}{2}x^2 + \dfrac{\omega^2}{2}y^2
\label{ham2dof}
\end{equation}
where $m_1$ and $m_2$ are the masses of the $x$ and $y$ DoF respectively, and $\omega$ is the angular frequency of the harmonic oscillator. For our analysis we will set $m_1 = m_2 = a = b = \omega = 1$. We can easily include the effect of dissipation in the system by writing Hamilton's equations in the form:
\begin{equation}
	\begin{cases}
	\dot{x} = \dfrac{\partial H}{\partial p_x} = \dfrac{p_x}{m_1} \\[.3cm]
	\dot{p}_x = -\dfrac{\partial H}{\partial x} - \gamma_x p_x = bx - ax^3 - \gamma_x p_x \\[.3cm]
	\dot{y} = \dfrac{\partial H}{\partial p_y} = \dfrac{p_y}{m_2} \\[.3cm]
	\dot{p}_y = -\dfrac{\partial H}{\partial y} - \gamma_y p_y = \omega y - \gamma_y p_y
	\end{cases}
	\label{ham_eqs}
\end{equation}
where $\gamma_x,\gamma_y > 0$ measure the dissipation strength in each DoF. We will look in this work at the case where dissipation has the same influence on both DoF, that is, $\gamma_x = \gamma_y = \gamma$, and we fix for all our calculations $m_1 = m_2 = a = b = \omega = 1$.

In order to analyze the dynamics of the system in Eq. \eqref{ham_eqs}, we define the following Poincar\'e surfaces of section (PSOS):
\begin{equation}
\begin{split}
    \Sigma_1 &= \left\lbrace (x,y,p_x,p_y) \in \mathbb{R}^4 \; \Big| \; y = 0 \; ,\; p_y > 0 \; \right\rbrace \\[.2cm]
    \Sigma_2 &= \left\lbrace (x,y,p_x,p_y) \in \mathbb{R}^4 \; \Big| \; p_y = 0 \; ,\; p_x > 0 \; \right\rbrace \\[.2cm]
    \Sigma_3 &= \left\lbrace (x,y,p_x,p_y) \in \mathbb{R}^4 \; \Big| \; x = -0.4 \; ,\; p_x > 0 \; \right\rbrace
\end{split}
\label{psecs}
\end{equation}
The main reasons for choosing these phase space slices are the following. First, $\Sigma_1$ serves to obtain a global picture of the basins of attraction of the stable equilibrium points located at the wells of the system. Second, we will use the section $\Sigma_2$ to probe the ellipsoidal geometry of the stable manifolds, since it cuts the ellipsoid longitudinally across its equator. Finally, $\Sigma_3$ is a transverse plane to the ellipsoid that would yield a circular cross-section, which gives us the classical picture known as a reactive islands \cite{almeida1990}. This name comes from the fact that all initial conditions chosen inside the closed curve obtained in this cut are reactive trajectories that cross the index-1 saddle region on forward time, and any other initial condition outside it never crosses the saddle along its evolution.

\begin{figure}[htbp]
	\begin{center}
    A)\includegraphics[scale=0.27]{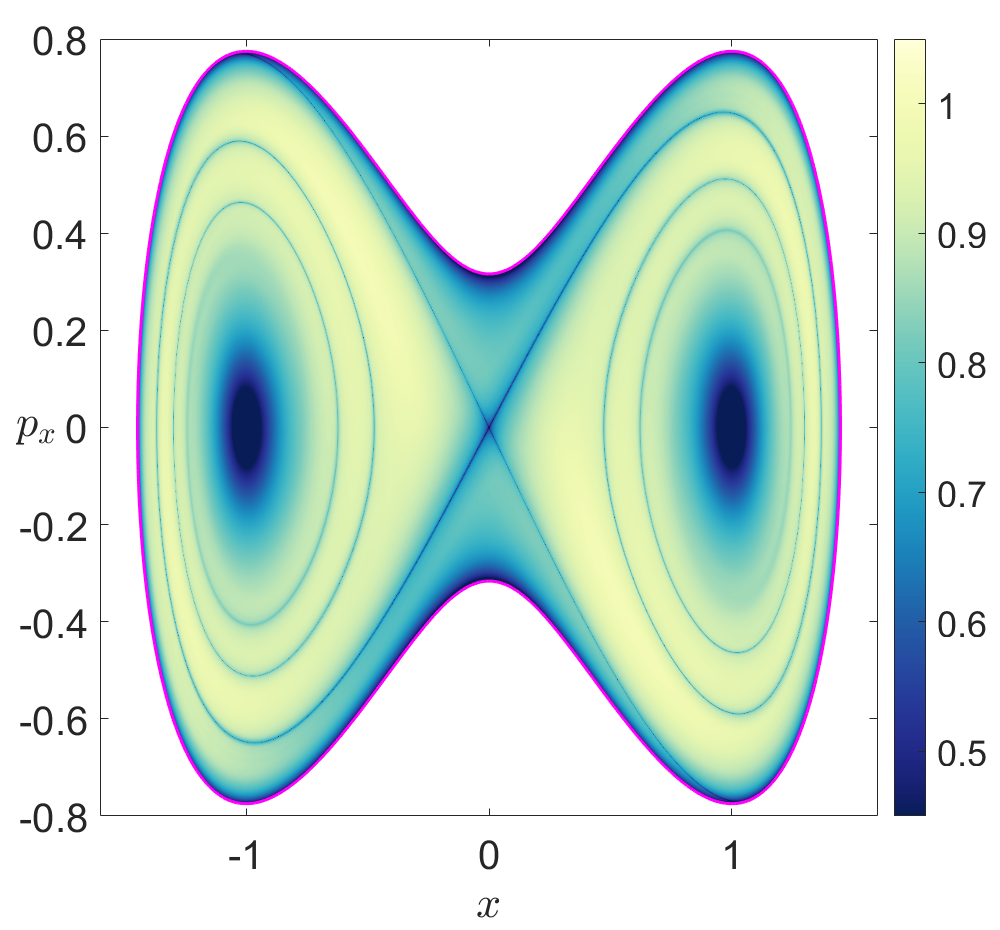} 
    B)\includegraphics[scale=0.27]{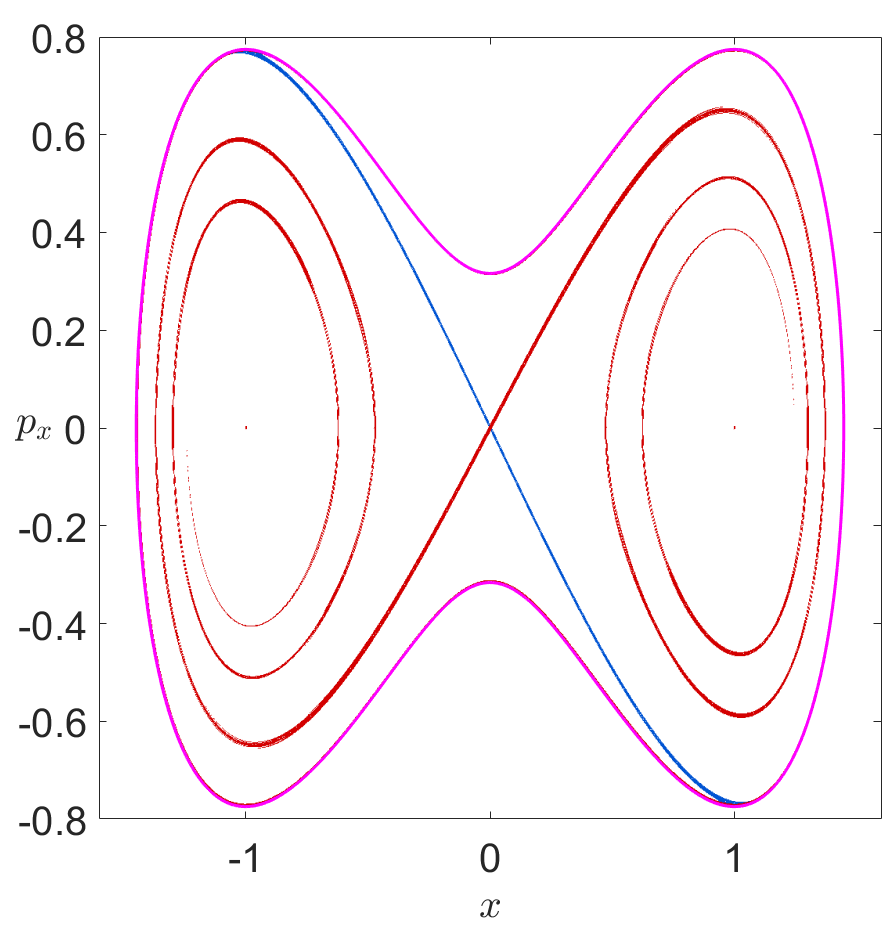} 
    C)\includegraphics[scale=0.27]{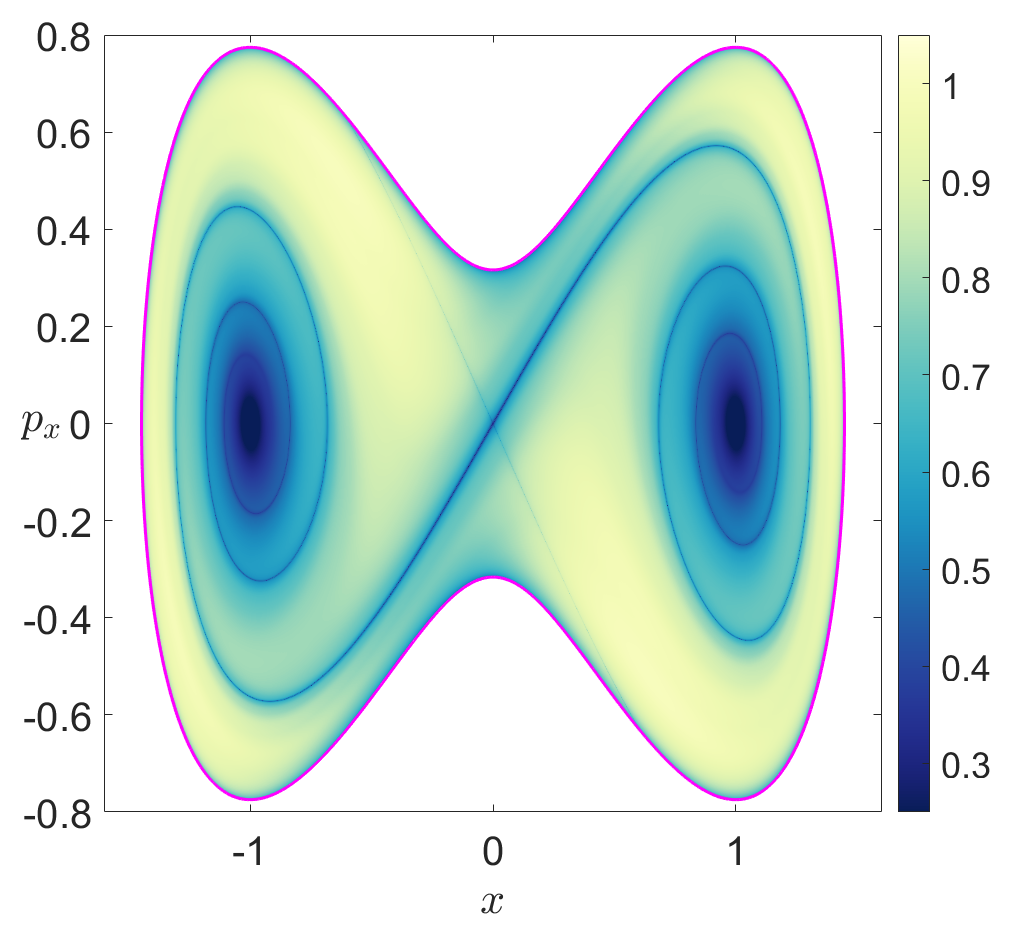} 
    D)\includegraphics[scale=0.27]{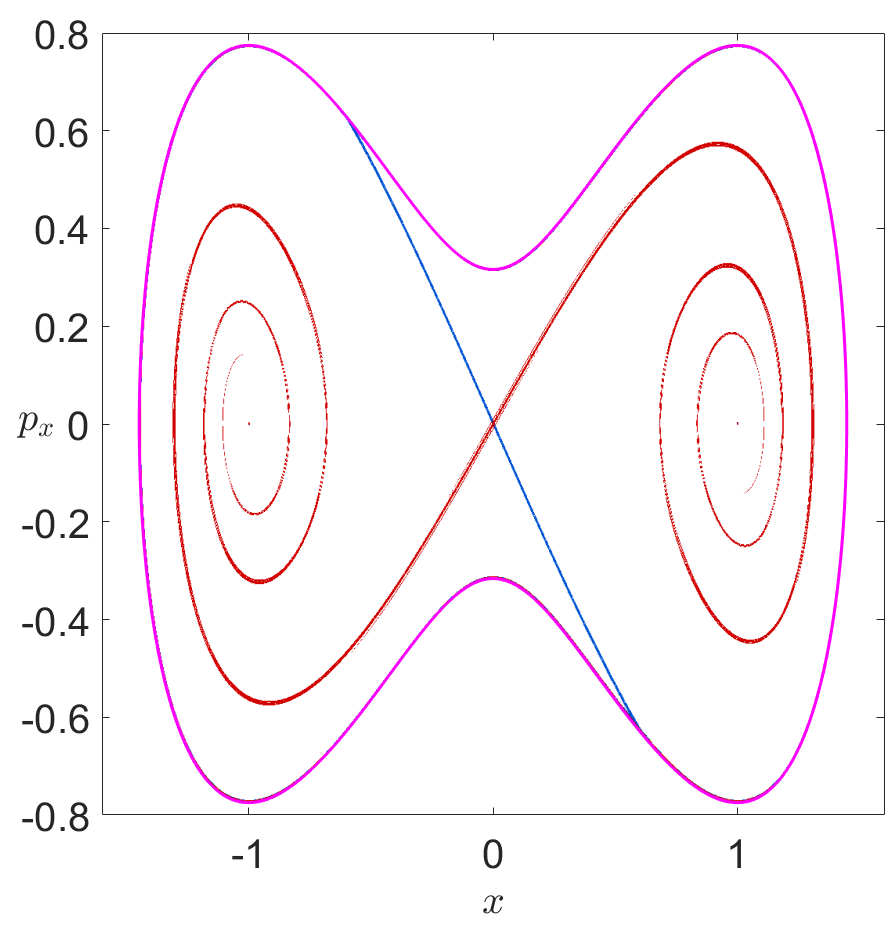}
    E)\includegraphics[scale=0.27]{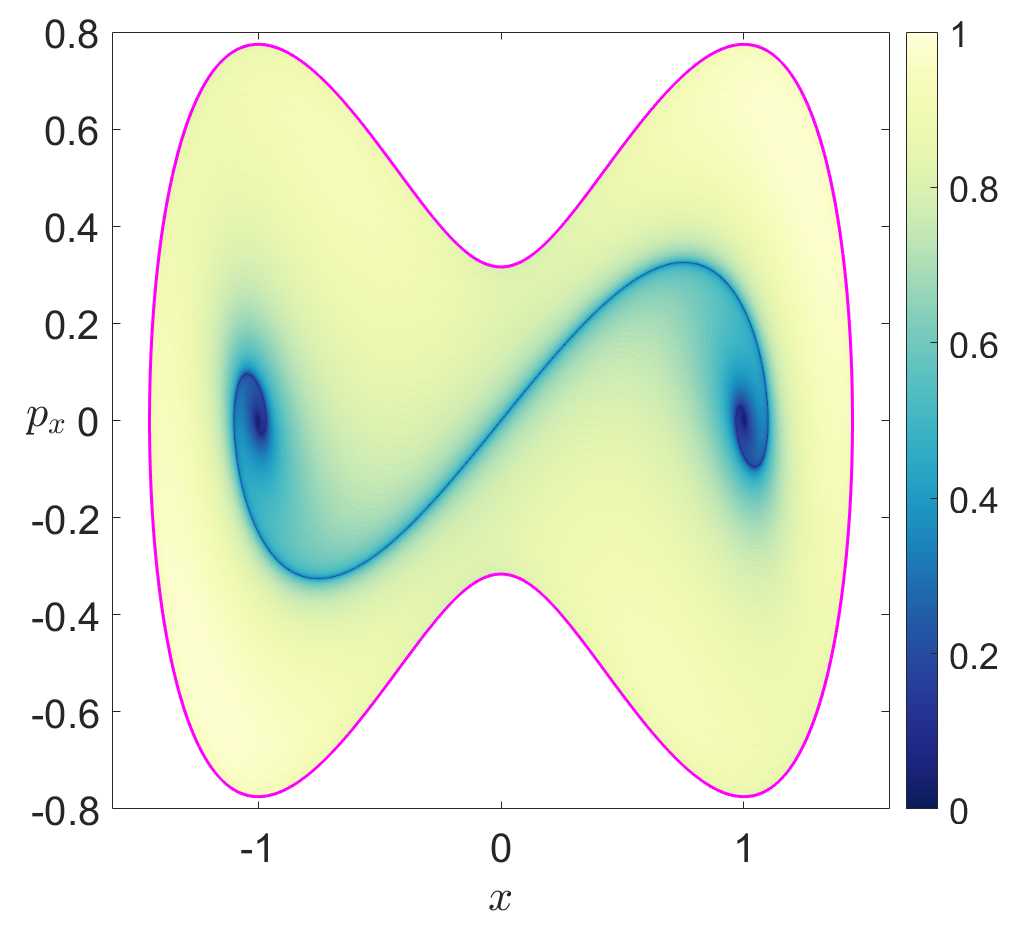} 
    F)\includegraphics[scale=0.27]{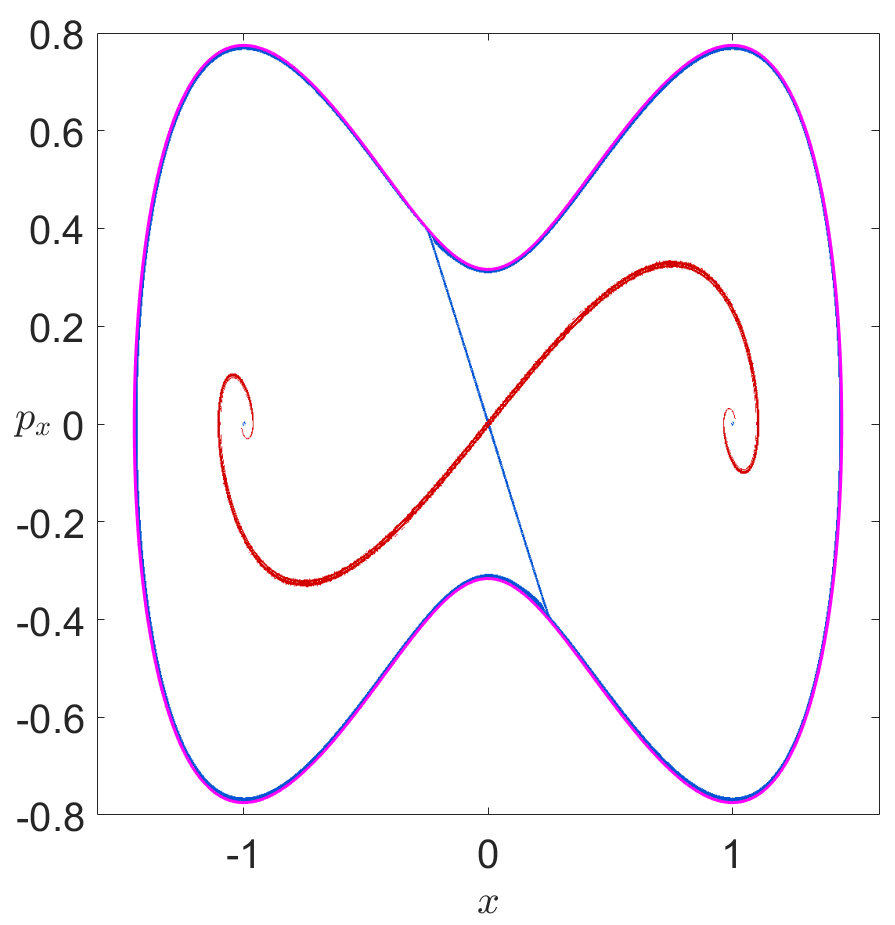}
	\end{center}
	\caption{Lagrangian descriptors calculated for the double well system in Eq. \eqref{ham_eqs} on the Poincar\'e section $\Sigma_1$ defined in Eq. \eqref{psecs}, using $p = 1/2$ and an integration time of $\tau_f = \tau_b = 15$. The initial energy of the system is $H_0 = 0.05$. Panels in each row correspond, from top to bottom, to dissipation values of $\gamma = 0.1, \, 0.25, \, 1$, respectively. We have depicted the initial energy boundary of the system as a magenta curve. The first column displays the total LD values, and the column on the right shows the sable (blue) and unstable (red) manifolds extracted from the LD scalar field by means of applying the laplacian operator.}
	\label{ld_dwell_sec_y_0}
\end{figure}

We start looking at the effect of increasing the strength of the damping in the geometry of the invariant manifolds. It is straightforward to show that the conservative Hamiltonian system ($\gamma = 0$) has an equilibrium point of saddle$\times$center stability type at the origin. When dissipation is included, its stability changes, and converts into a saddle$\times$focus. We consider here three different values for the damping coefficient, $\gamma = 0.1,\, 0.25,\, 1$. If we compute LDs in the PSOS $\Sigma_1$, see Fig. \ref{ld_dwell_sec_y_0}, the basins of attraction of the wells in the system are revealed. The spiraling of the unstable manifold that emanates from the equilibrium point at the origin and ends in the equilibrium point in the left and right wells reduces as the damping increases, and this explains why it takes less time for trajectories to approach the stable equilibria. Notice also how the stable manifolds get closer to the origin for large values of the friction, and this represents that the amount of trajectories that can cross from well to well across the saddle point of the PES reduces significantly with large damping values. This implies that the transition ellipsoid becomes smaller in size. Notice that the stable manifold in the total LD scalar field displayed in the left column of Fig. \ref{ld_dwell_sec_y_0} becomes harder to locate as dissipation is increased, but this is expected in the output of LDs, as we have discussed for other examples in this work, since the time scales of the system separate. If we would like to fix this visualization issue, one could integrate initial conditions forward and backward for different times. However, this is not necessary because the manifolds can be directly extracted from the forward and backward components of LDs by means of applying the gradient or laplacian of the LD scalar field, or other edge detection algorithms used in image processing.

\begin{figure}[htbp]
	\begin{center}
A)\includegraphics[scale=0.27]{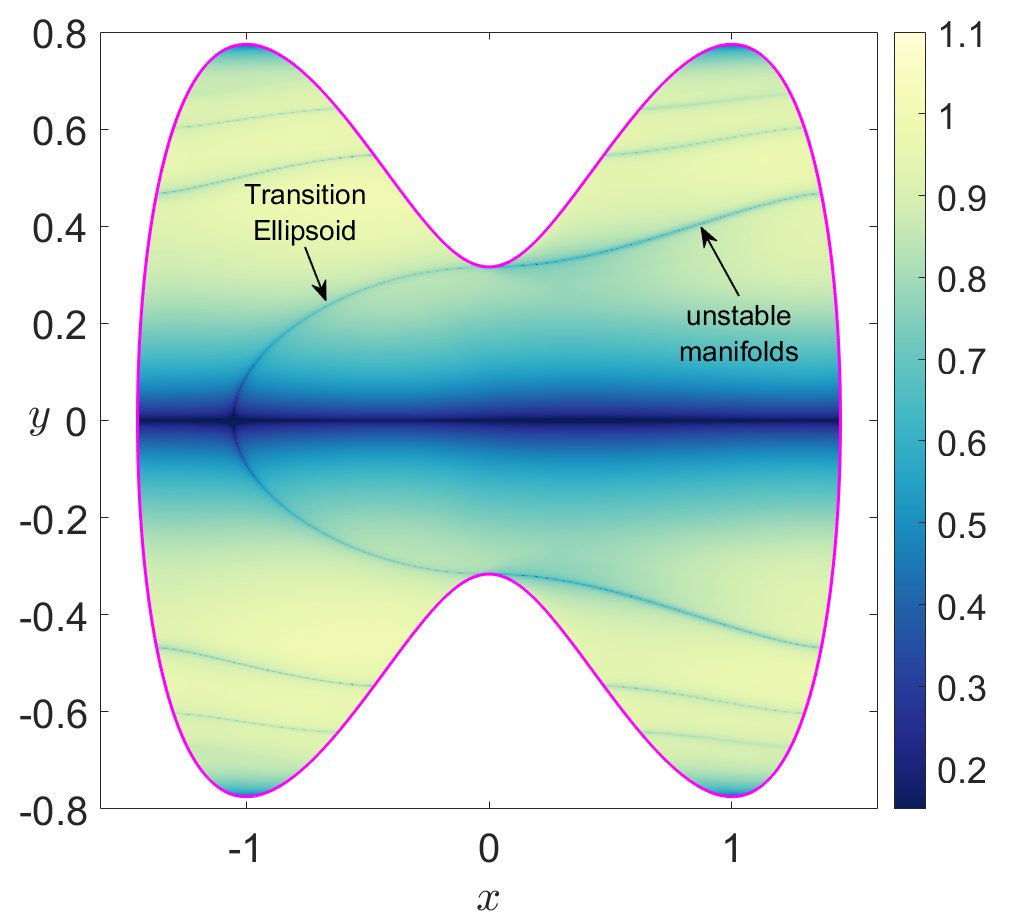}
B)\includegraphics[scale=0.27]{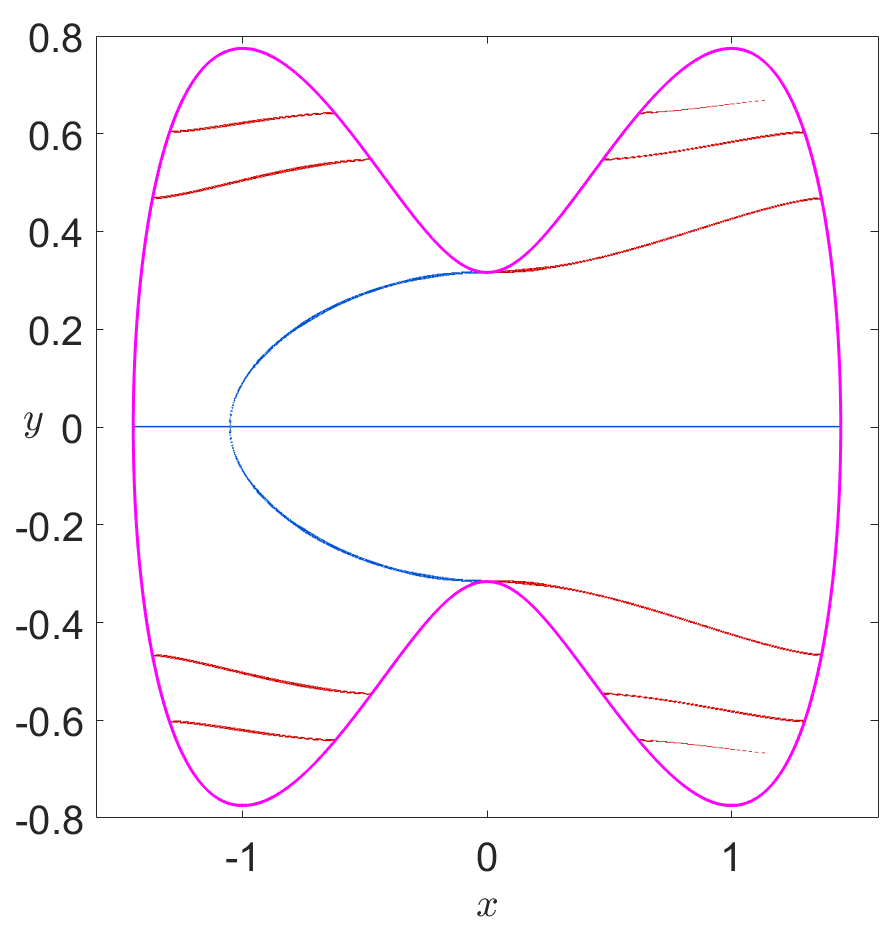} 
C)\includegraphics[scale=0.27]{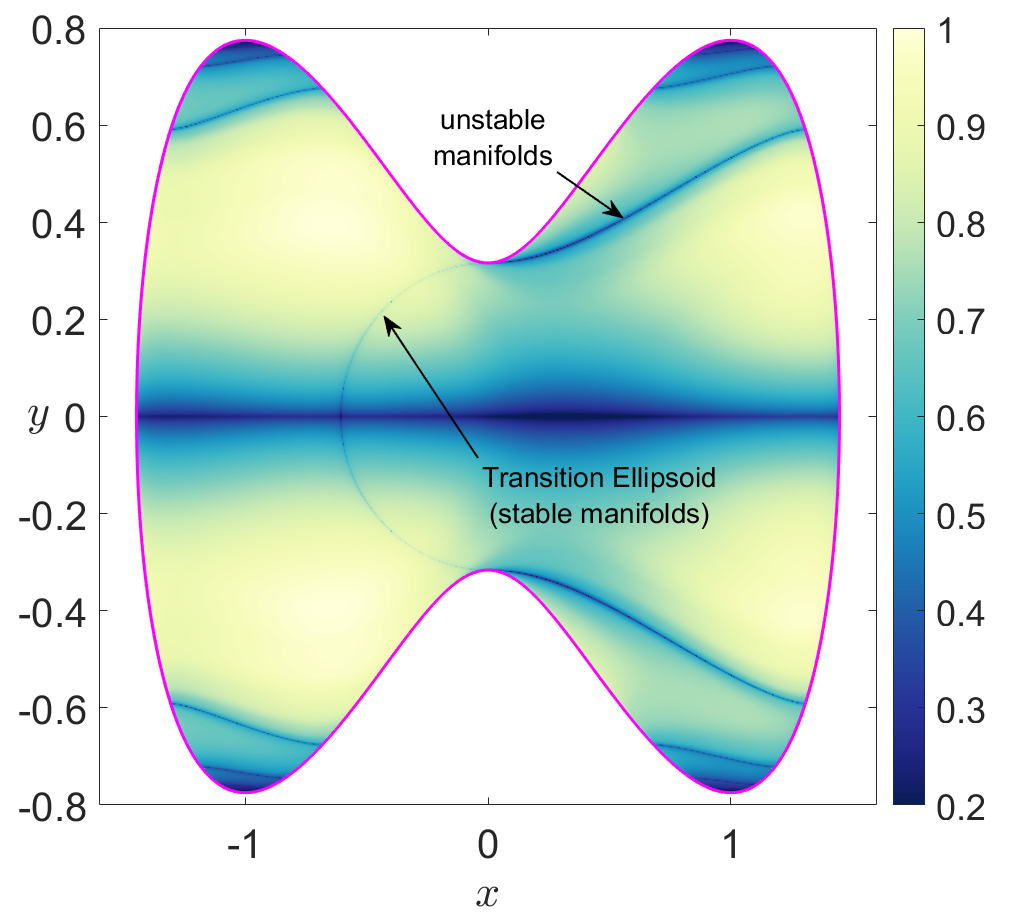} 
D)\includegraphics[scale=0.27]{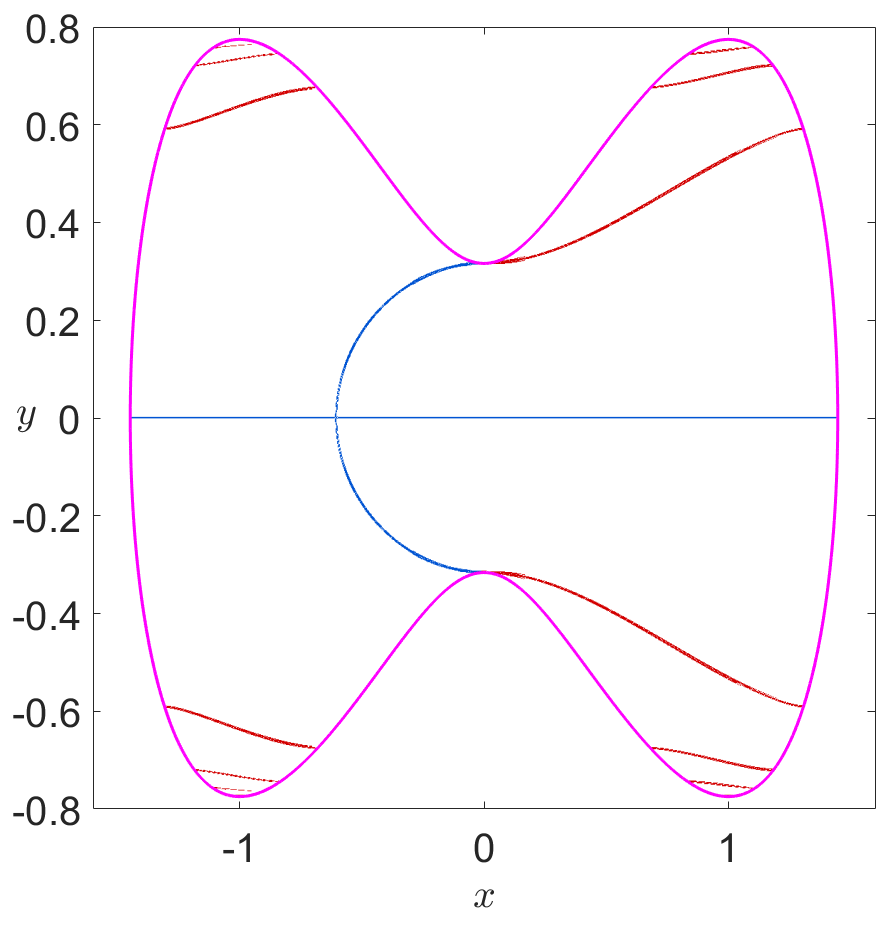}
E)\includegraphics[scale=0.27]{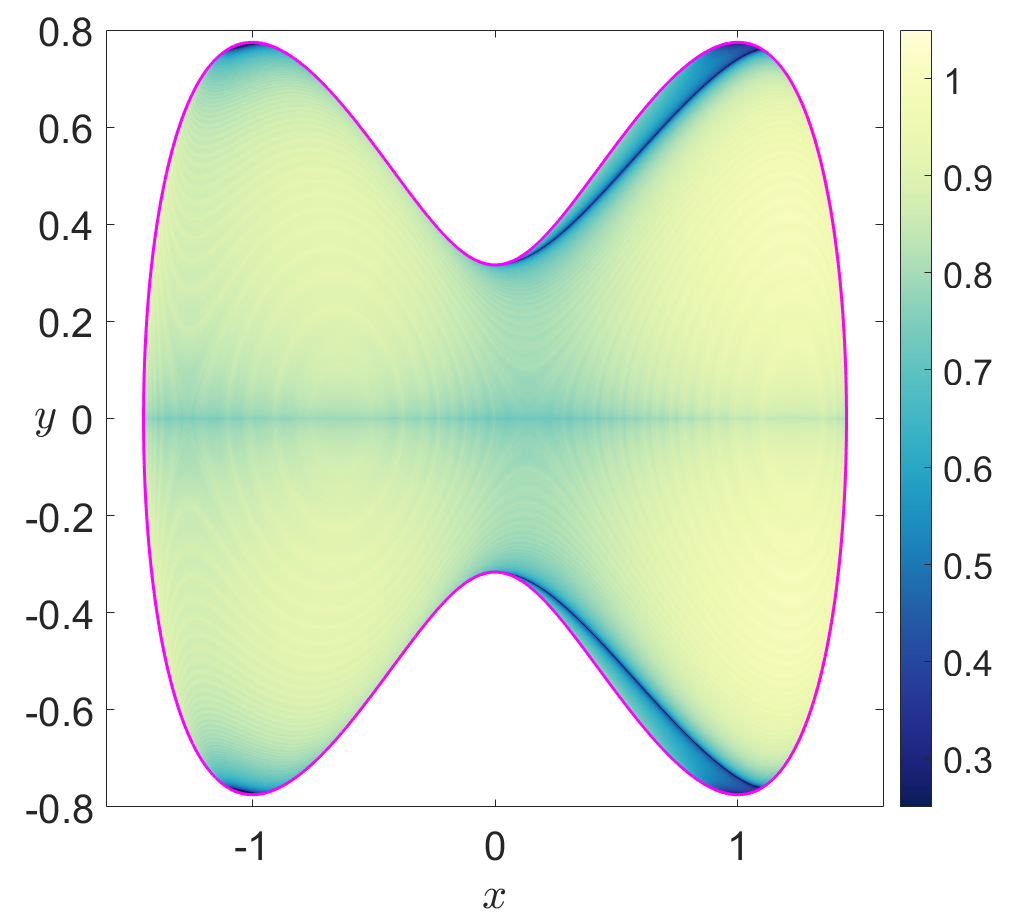} 
F)\includegraphics[scale=0.27]{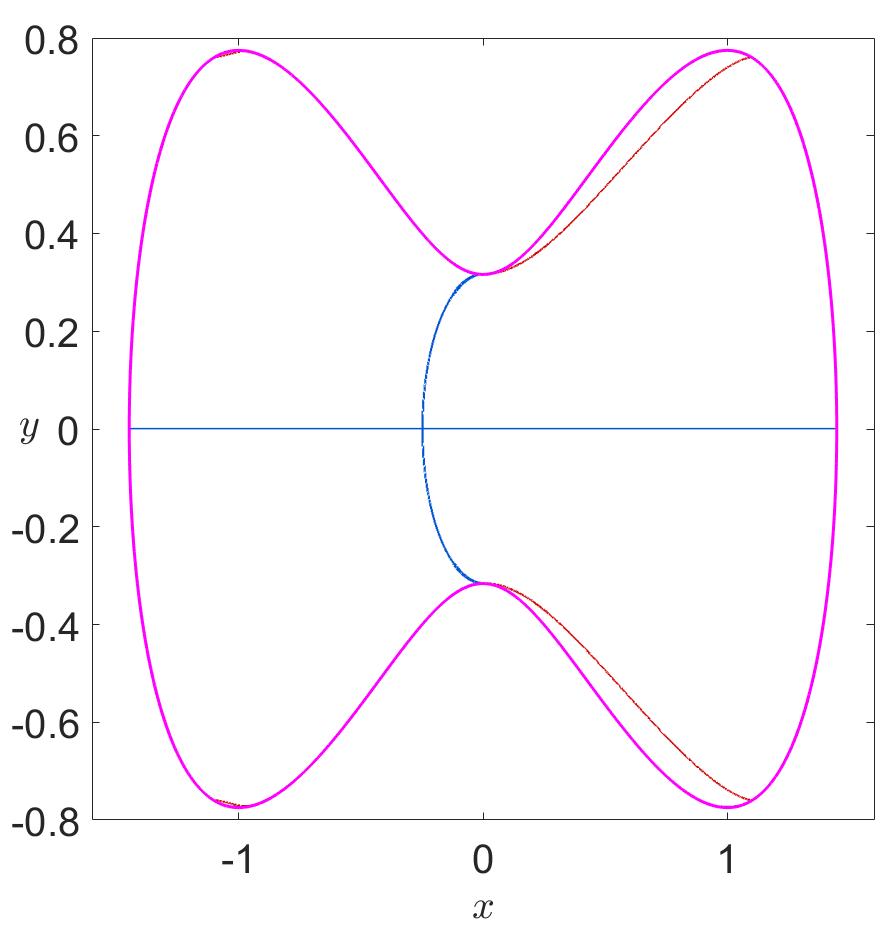}
	\end{center}
	\caption{Lagrangian descriptors calculated for the double well system in Eq. \eqref{ham_eqs} on the Poincar\'e section $\Sigma_2$ defined in Eq. \eqref{psecs}, using $p = 1/2$ and an integration time of $\tau_f = \tau_b = 15$. The initial energy of the system is $H_0 = 0.05$. Panels in each row correspond, from top to bottom, to dissipation values of $\gamma = 0.1, \, 0.25, \, 1$, respectively. We have depicted the initial energy boundary of the system as a magenta curve. The first column displays the total LD values, and the column on the right shows the sable (blue) and unstable (red) manifolds extracted from the LD scalar field by means of applying the laplacian operator.}
	\label{ld_dwell_sec_py_0}
\end{figure}

\begin{figure}[htbp]
	\begin{center}
A)\includegraphics[scale=0.27]{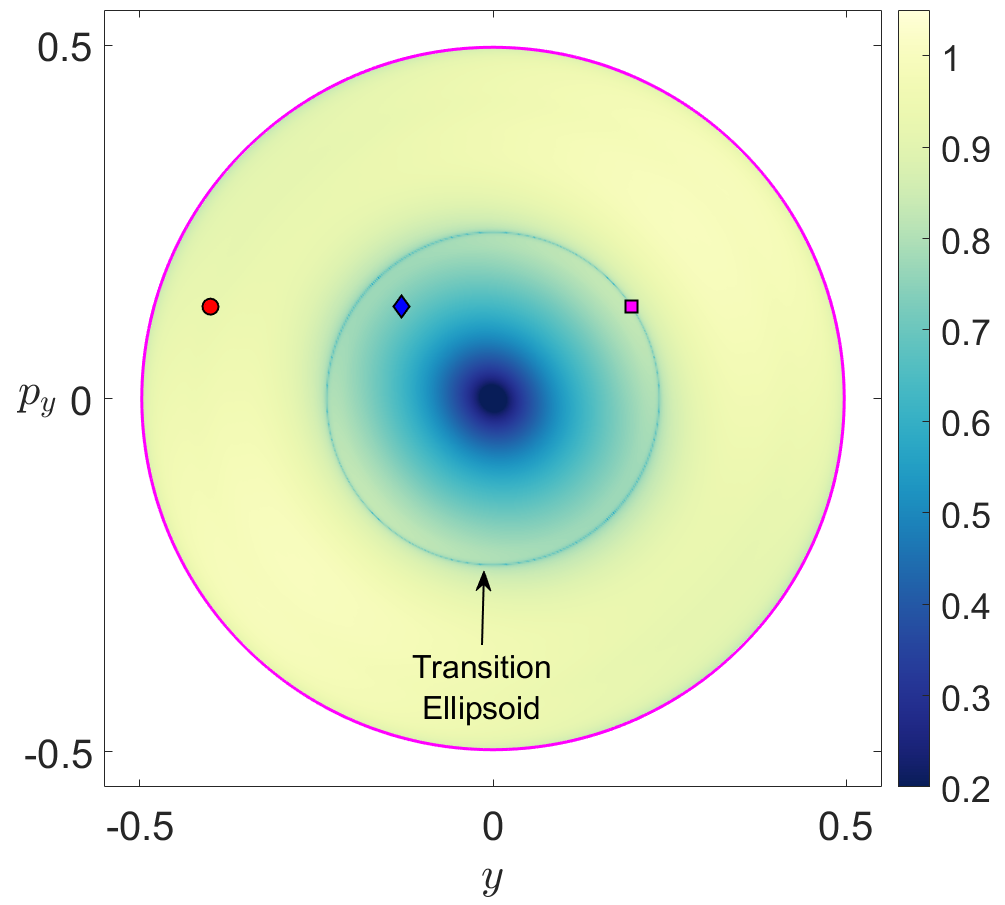}
B)\includegraphics[scale=0.28]{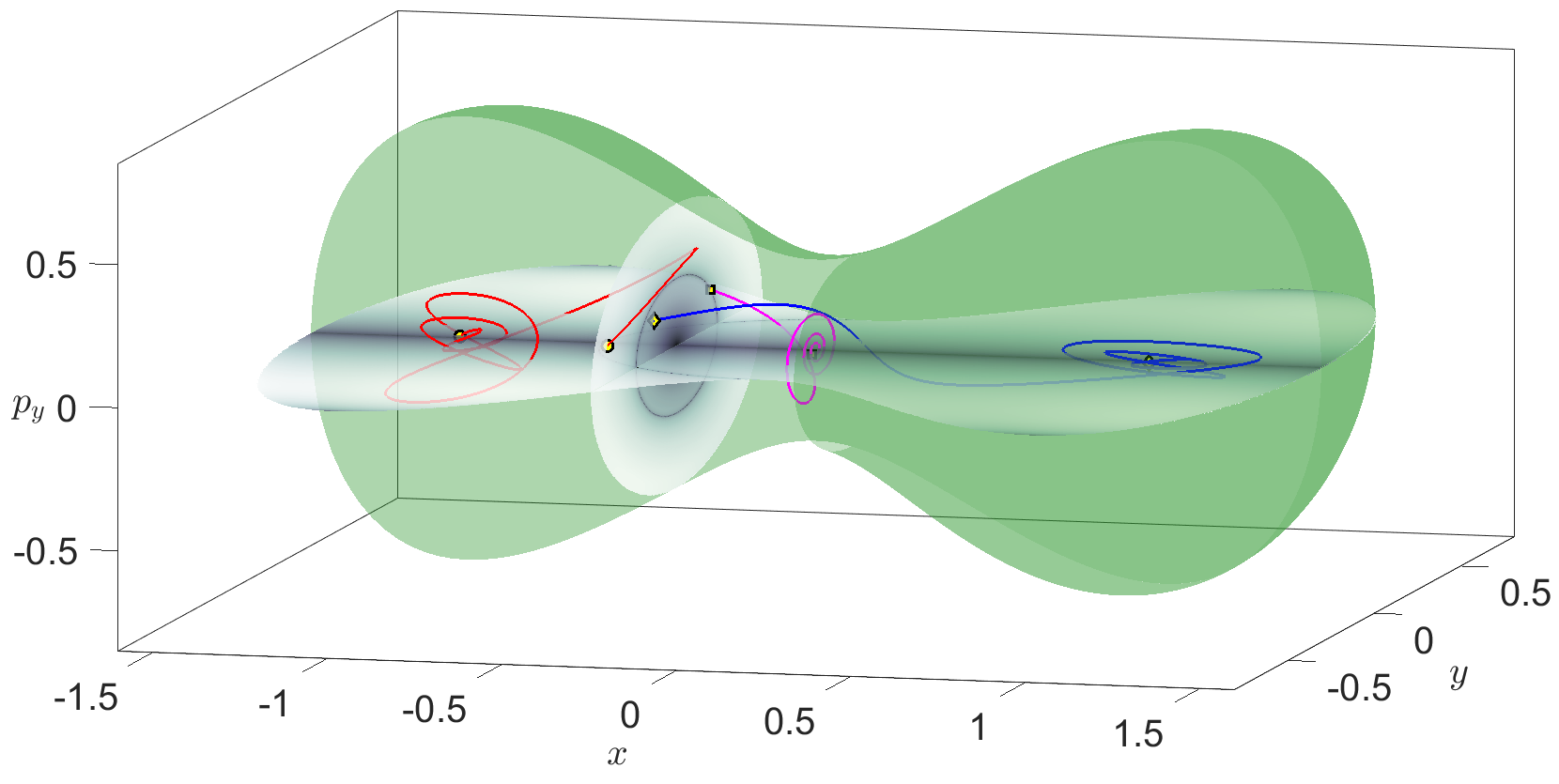} 
C)\includegraphics[scale=0.27]{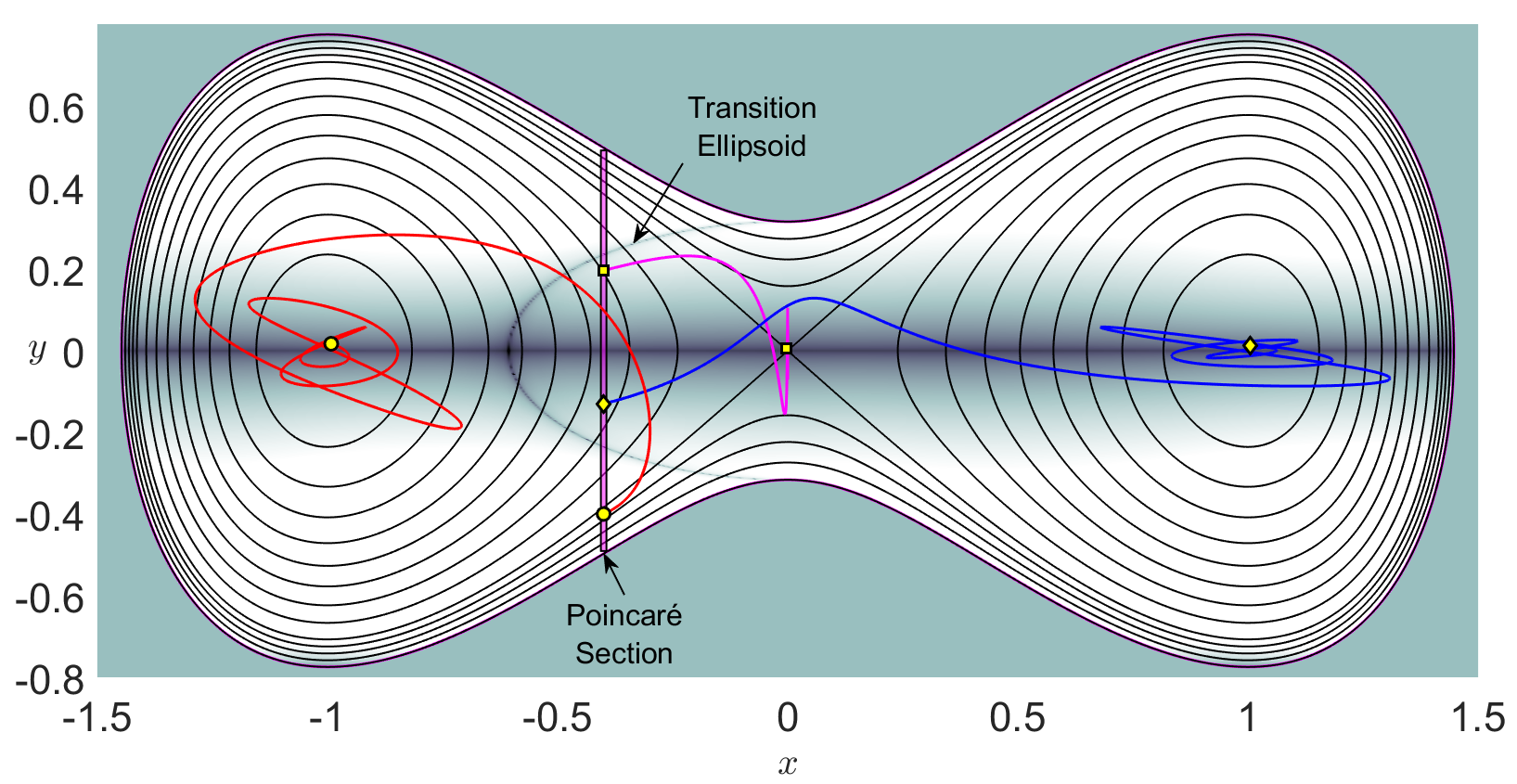} 
	\end{center}
	\caption{Computation of LDs and visualization of three-dimensional dynamics for the double well system in Eq. \eqref{ham_eqs} with initial energy $H_0 = 0.05$ and dissipation strength $\gamma = 0.25$. A) LDs on the Poincar\'e section $\Sigma_3$ defined in Eq. \eqref{psecs}, using $p = 1/2$ and an integration time of $\tau_f = \tau_b = 15$. We have marked three different initial conditions, the red circle is outside the transition ellipsoid, the blue diamond is inside, and the magenta square is on its boundary. The magenta curve depicts the initial energy boundary of the system. B) Three dimensional dynamics in forward time of the initial conditions selected in panel A) as they evolve from $t_0 = 0$ to $t = 20$. We have superposed LDs calculated on the phase space slices $\Sigma_2$ and $\Sigma_3$, and the green surface represents the initial energy shell. C) Configuration space projection of the trajectories in panel B), superposed with LDs calculated on the section $\Sigma_2$. Potential energy contours are depicted as black cuves and the outer grey area corresponds the forbidden region for the initial energy of the system.}
	\label{3d_ld_trajs}
\end{figure}

In order to clearly visualize the ellipsoidal geometry of the stable manifold we calculate LDs on the configuration space PSOS $\Sigma_2$. We depict in Fig. \ref{ld_dwell_sec_py_0} the LD values (left column) and the stable (blue) and unstable (red) manifolds obtained from the laplacian operator of the LD function (right column). In this phase space slice, the ellipsoidal geometry is nicely captured by LDs, confirming its existence in the system. Also, it shrinks in size as the damping gets large. To finish, we analyze the system further to validate that the structure highlighted by LDs is in fact the stable manifold, and that this transition ellipsoid controls the transport of trajectories across the index-1 saddle. To do so, we calculate LDs on the PSOS $\Sigma_3$ using a dissipation of $\gamma = 0.25$, and show the results of the forward LD in Fig. \ref{3d_ld_trajs} A). This slice intersects the ellipsoid transversely into a circle. To test if this closed curve corresponds to the ellipsoid boundary, we select three different initial conditions, the red circle outside the ellipsoid, the blue diamond inside the manifold, and the magenta square on the boundary of the stable manifold: We evolve them forward in time and plot their trajectories in 3D in panel B) and projected onto configuration space in C). This numerical experiment convincingly shows that LDs is correctly detecting the boundary of the transition ellipsoid, since the red initial condition does not cross the index-1 saddle and is attracted by the stable equilibrium point of the left well. On the other hand, the blue trajectory moves across the saddle point, and then it ends its evolution at the right well equilibrium point. Finally, the magenta trajectory which starts exactly at the boundary of the ellipsoid, behaves as expected, and asymptotically tends towards the saddle point at the origin. It does so following a spiral-like trajectory that evidences the saddle-focus stability nature of the saddle equilibrium.

\section{Conclusions}
\label{sec:conc}

In this work we have demonstrated how the method of Lagrangian descriptors can be successfully implemented to reveal the relevant phase space structures of dissipative dynamical systems. This study substantially broadens the applicability spectrum of this tool and exemplifies how a simple scalar diagnostic technique can provide relevant insights in order to understand dynamical phenomena in real problems. 

Our analysis of different classical models from nonlinear dynamics has shown that this technique has the capability to detect the presence of limit cycles, basins of attraction, strange attractors and slow manifolds. Moreover, by including the possibility of integrating trajectories forward and backward for different times we have generalized the definition of LDs, and this allows us to account for the separation of time scales displayed in these type of problems. On the other hand, in the context of Hamiltonian systems subject to dissipative forces, where the topology of the stable and unstable manifolds changes due to the presence of friction, we have explored this effect on a double well system with two degrees of freedom. The study of this problem by means of LDs convincingly illustrates the advantage that this tool brings for the location of transition ellipsoids which characterize the phase space trajectories that are allowed to move across the index-1 saddle separating two well regions of the underlying potential energy surface.

In the future, our work will focus on applying Lagrangian descriptors to develop a more detailed understanding on the impact of dissipation on the geometry of invariant manifolds and their interaction, and the role that these phase space structures play on governing the dynamics in many Hamiltonian systems subject to friction that arise in real problems and applications.

\section*{Acknowledgments} 

VJGG would like to acknowledge the financial support received from the EPSRC Grant No. EP/P021123/1 and the Office of Naval Research Grant No. N00014-01-1-0769 for his research visits over the past years to the School of Mathematics, University of Bristol. Many fruitful discussions with Prof. Stephen Wiggins during this period have inspired the development of this work.

\bibliographystyle{natbib}
\bibliography{SNreac}

\end{document}